\documentclass{amsart}

\numberwithin{equation}{section}

\usepackage{amssymb}
\usepackage{enumerate, xspace}
\usepackage{esint}
\usepackage[normalem]{ulem}
\usepackage{marvosym} 
\usepackage{bbm}       
\usepackage[backgroundcolor=blue!20!white, linecolor=blue!20!white,
textsize=footnotesize]{todonotes}
\usepackage[colorlinks]{hyperref}
\usepackage{geometry}
\newtheorem{thm}{Theorem}[section]
\newtheorem{lem}[thm]{Lemma}
\newtheorem{sub-lem}[thm]{Sub-Lemma}
\newtheorem{cor}[thm]{Corollary}
\newtheorem*{notation}{Notation}
\newtheorem{sublem}[thm]{Sublemma}

\newtheorem{prop}[thm]{Proposition}

\newtheorem{defin}[thm]{Definition}

\newtheorem{rem}[thm]{Remark}

\newcommand\cA{{\mathcal A}}
\newcommand\cB{{\mathcal B}}
\newcommand\cC{{\mathcal C}}
\newcommand\cD{{\mathcal D}}
\newcommand\cE{{\mathcal E}}
\newcommand\cF{{\mathcal F}}
\newcommand\cG{{\mathcal G}}

\newcommand\cK{{\mathcal K}}
\newcommand\cL{{\mathcal L}}
\newcommand\cO{{\mathcal O}}
\newcommand\cM{{\mathcal M}}

\newcommand\cQ{{\mathcal Q}}
\newcommand\cS{{\mathcal S}}

\newcommand\cW{{\mathcal W}}

\newcommand\bB{{\mathbb B}}
\newcommand\bC{{\mathbb C}}

\newcommand\bE{{\mathbb E}}

\newcommand\bH{{\mathbb H}}
\newcommand\bL{{\mathbb L}}
\newcommand\bN{{\mathbb N}}

\newcommand\bP{{\mathbb P}}

\newcommand\bR{{\mathbb R}}

\newcommand\bZ{{\mathbb Z}}

\newcommand\ve{\varepsilon}

\newcommand\vf{\varphi}

\newcommand{\tB}{\tilde{B}}
\newcommand{\tg}{\tilde{g}}
\newcommand\tpsi{\tilde{\psi}}
\newcommand\tJ{\widetilde{J}}
\newcommand\tT{\widetilde{T}}

\newcommand{\bbe}{\mathbbm e}
\newcommand{\bbm}{\mathbbm m}

\newcommand{\bbg}{\mathbbm g}

\newcommand\Id{{\mathbbm{1}}}
\newcommand{\Pe}{\bP_{\! e}}
\newcommand{\tri}{{|\:\!\!|\:\!\!|}}
\newcommand{\hS}{\hat{S}}
\newcommand{\cop}{ \upsilon}     
\newcommand{\zm}{z_-}               
\newcommand{\zp}{z_+}                
\newcommand{\fA}{\mathfrak A}      
\newcommand{\fa}{\mathfrak a}         
\newcommand{\hT}{\widehat{T}}
\newcommand{\wcL}{\widehat \cL}

\newcommand{\Const}{C_{\#}}
\newcommand{\const}{c_{\#}}

\newcommand{\sign}{\operatorname{sign}}
\newcommand{\musrb}{\mu_{\mbox{\tiny SRB}}}
\newcommand{\diam}{\operatorname{diam}}

\begin{document}

\title{Central Limit Theorem for Sequential Dynamical Systems}
\author{Mark F. Demers}
\address{Mark F. Demers\\
Department of Mathematics\\
 Fairfield University, Fairfield CT 06824, USA.}
\email{{\tt  mdemers@fairfield.edu}}
\author{Carlangelo Liverani}
\address{Carlangelo Liverani\\
Dipartimento di Matematica\\
II Universit\`{a} di Roma (Tor Vergata)\\
Via della Ricerca Scientifica, 00133 Roma, Italy.
}
\address{\vskip-16pt \noindent Mathematics Department\\
University of Maryland\\
4176 Campus Drive - William E. Kirwan Hall\\
College Park, MD 20742-4015}
\email{{\tt liverani@mat.uniroma2.it}}
\date{\today}
\vspace*{-.5\baselineskip}
\begin{abstract}
We present a general approach to establish the Central Limit Theorem with error bounds for sequential dynamical systems.  The main tool we develop is the application to this setting of a projective metric on complex cones, following
the ideas introduced by Rugh and Dubois.
To demonstrate the power of the proposed setting, we apply it to both sequential expanding maps, where similar results are known, and to sequential dispersing billiards, for which no such results are currently known.
\end{abstract}
\thanks{This work was supported by the Grant (PRIN 2017S35EHN) and by the MIUR Excellence Department Project MatMod@TOV. C.L.~acknowledges membership to the GNFM/INDAM and support from the Brin Center, University of Maryland.  M.D.~was partially supported by National Science Foundation Grant DMS 2350079. Finally, we would like to thank Dmitry Dolgopyat for many valuable discussions.}
\maketitle

\tableofcontents

\section{Introduction}
This paper is concerned with obtaining a self-normed Central Limit Theorem with error terms for 
sequential dynamical systems. 
Such results are relevant for many open problems; let us just mention random 
Lorentz gases of the type initially presented in \cite{Le03}, in which a particle moves between cells formed by randomly placed obstacles
in an approximately $\mathbb{Z}^2$ lattice.  In each cell, the particle finds a different arrangement of scatterers.  
From the particle's point of
view, it undergoes a sequence of different billiard maps corresponding to the cells it visits.  Such a model was 
considered in several works \cite{Le03, Le06, AL20} and most recently in \cite{DL22}, which introduced the application
of Birkhoff cones to dispersing billiards. 
More broadly, many results regarding the Central Limit Theorem for sequential systems already exist, for example, \cite{CR07,NMTV, NTV18, HS20, H20, LS20, KL21, NFT21, S22, DK23} in which are investigated either random dynamical systems or sequential expanding maps, possibly non-uniformly expanding. 

In the present paper, we present a general setting to establish a self-normed Central Limit Theorem with error terms for sequential dynamical systems, which allows us to treat both observables that are not bounded and discontinuous hyperbolic dynamical systems, including dispersing billiards. The results are obtained via the study of the properties of the corresponding transfer operators. As we are interested in presenting the general strategy, we do not push for optimal results. Since some estimates rest on a Taylor expansion, better results can be obtained by computing more terms. Note, however, that to obtain optimal results, it is necessary to have an estimate of the
 higher momenta in terms of the second momenta, and it is not clear how to achieve this in the present generality (see Remark \ref{rem:momenta} for a more detailed explanation
of this issue).

In Section \ref{sec:setting}, we present two sets of abstract hypotheses:  the first hypothesis assumes only that the real transfer operator is contracting in a projective metric, together with a uniform bound on the complex twisted transfer operator; the second hypothesis assumes that the complex transfer operator is also contracting. The second possibility yields much sharper results. We explain in Section~\ref{sec:complex} under which conditions the latter extra information on the transfer operator can be obtained via the theory of complex cones introduced by Rugh \cite{Ru10} and further developed by Dubois \cite{Du09, Du11}. Since this theory is of independent interest and is spread among several articles, for the reader's convenience, we provide a short, self-contained account of the theory in this section. Our presentation of complex cones is not as 
general as the one that can be found in \cite{Ru10, Du09, Du11}, but it is more than sufficient for all the applications 
to the Central Limit Theorem that we can think of.

To illustrate the broad applicability of our result, we apply it to two emblematic classes of examples; indeed more than half of the paper is devoted to these examples. The first application is to smooth expanding maps 
(the simplest possibility and one for which the most results already exist); the second application is to dispersing billiards (a technically difficult case motivated by the Lorentz gas and for which no sequential results of this type currently exist).  We remark that other cases, such as Anosov and Axiom A diffeomorphisms, lie in some sense between the applications contained in this paper.  The analysis of hyperbolic systems in dimensions greater than two is not an obstruction to the present technique.  Indeed, anisotropic spaces for Anosov and Axiom A diffeomorphisms (including for those with singular SRB measures) in any dimension have been constructed by a variety of authors \cite{GL1, BT1, GL2} and the connection between anisotropic spaces and projective cones is described in \cite{DKL}.  While the details of such applications would make the present article excessively long, yet the
technique presented here is applicable to such systems as well.  Indeed, the expanding maps example shows that it is not necessary that the map have a common invariant measure. On the other hand, the billiards examples are studied using Banach spaces of distributions, which implicitly shows that the existence of an absolutely continuous invariant measure is also not necessary.

The application to expanding maps in Section~\ref{sec:expanding} is done to compare our present 
Theorems~\ref{thm:main} and \ref{thm:main_bis} with existing results. We show that our technique reproduces all existing 
results apart from the ones in \cite{DH24}, which just appeared. In addition, in Proposition \ref{prop:variance} and Corollary \ref{cor:verify} we provide a new criterion to verify the growth of the variance in the sequential case, while Theorem \ref{thm:random_var} provides a sharper criterion for the random case. This is an issue often not addressed in the literature but fundamental in applications.  The strategy used in \cite{DH24} is similar to ours when applied to the case in which one can use complex cones but takes advantage of the fact that it deals only with bounded observables and expanding maps. In particular, it uses a martingale decomposition trick, adapted from \cite[Chapter 3]{Merlevede}, in order to estimate higher momenta via the second momenta. It is unclear how to extend such a trick to the level of generality we are addressing in the present work.  Possibly, something can be done using ideas in \cite[Theorem 1.2]{Li95c} applied to the stable foliation of the sequential maps, but as the present paper is already rather long, we refrain from exploring this issue.
The application to billiards in Section~\ref{sec:billiard} is brand new and shows the generality and power of the theory put forward in this paper. We expect this application to be relevant for the study of random Lorentz gases, a long-standing open problem.

For our motivating example of a particle moving through cells in a $\mathbb{Z}^2$ lattice \cite{Le03},
exponential loss of memory is proved in \cite{DL22} using projective cones;  however, the
present technique does not yet yield the Central Limit Theorem in this case due to a
lack of control of mixing in the $\mathbb{Z}^2$-coordinate.
It is, however, worth being precise about what is, and what is not, covered. Seen from the
moving particle, the model of \cite{Le03, Le06, AL20} is exactly a sequential billiard drawn
from a fixed family $\cF(\tau_*, \cK_*, E_*)$: at each collision, the particle meets the
scatterers of the cell it currently occupies, and the configuration changes only when it
leaves the cell. The sequential systems studied in Section~\ref{sec:billiard} are therefore
precisely the ones arising in that model, and Theorem~\ref{thm:billiard CLT} applies to them.
What is still missing is not the billiard dynamics, but the control of the $\bZ^2$
displacement. This, as remarked at the beginning, is the problem of the Lorentz gas in a random environment, which has been open for many decades.
So the present paper is a step forward in that direction, at least in the case in which the system has a drift, a case for which a technique to attack the problem has already been developed, e.g. see \cite{DK23}, provided one can control the inner dynamics, which is what the present paper shows how to do in a class of models.


\subsection*{Notation}
For the reader's convenience, here are the most relevant symbols used throughout the paper.

\smallskip\noindent\emph{Sequential dynamics.}
\nopagebreak
\begin{center}\footnotesize
\begin{tabular}{@{}p{3.4cm}@{}p{0.68\textwidth}@{}}
$M_k$, $f_k$ & compact Riemannian manifolds and measurable maps $f_k:M_k\to M_{k+1}$, $k\in\bN$; Section~\ref{sec:abstract}\\
$f_{k,j}$ & the composition $f_{k-1}\circ\cdots\circ f_j$, with $f_{j,j}(x)=x$\\
$\mu_k$ & reference probability measures on $M_k$, with $(f_k)_*\mu_k\ll\mu_{k+1}$\\
$\cL_k$ & transfer operator of $f_k$ with respect to the $\mu_k$, \eqref{eq:transfer_op}\\
$\cL_{k,\lambda,n}$ & twisted transfer operator $\cL_k(e^{i\sigma_n^{-1}\lambda\hat g_k}\,\cdot\,)$, \eqref{eq:twisted}; abbreviated $\cL_{k,\lambda}$ when $n$ is fixed, \eqref{eq:simpler}, with compositions $\cL_{k,j,\lambda}$, \eqref{eq:transfer_gen}\\
$\rho$, $\nu_0$, $\bE$ & initial density $\rho\in\cB_0$, initial measure $d\nu_0=\rho\, d\mu_0$, and expectation with respect to $\nu_0$\\
\end{tabular}
\end{center}

\smallskip\noindent\emph{Cones, spaces and norms.}
\nopagebreak
\begin{center}\footnotesize
\begin{tabular}{@{}p{3.4cm}@{}p{0.68\textwidth}@{}}
$C_k$, $V_k$ & the real cones and $V_k=\operatorname{span}(C_k)$, conditions (C-\ref{cond:positivity})--(C-\ref{cond:norm})\\
$\bbe_k$, $\preceq$ & order unit of $C_k$, normalised by \eqref{eq:normalize_e}, and the order it induces\\
$\|\cdot\|_k$, $\cB_k$ & the cone norm \eqref{eq:real cone norm} and the Banach space it completes\\
$\bB_k$ & complexification of $\cB_k$, Lemma~\ref{lem:complexification}; its norm is still written $\|\cdot\|_k$\\
$\cC_\bR$, $\cC_\bC$ & the real and complex cones of Section~\ref{sec:complex}, the latter defined by \eqref{eq:com_cone_def}\\
$\cS$, $\bbm_k$ & the functionals defining $\cC_\bR$ and the norm-dominating functional $\bbm_k\in\cS_*$, \eqref{eq:maybetoomuch}\\
$d_H$, $\delta_\cC$, $\Delta$ & Hilbert metric on $\cC_\bR$, its complex analogue \eqref{eq:gauge}, and the diameter of an image cone\\
$C_*$ & the uniform bound on compositions of the $\cL_k$, condition (C-\ref{cond:powerbounded})\\
\end{tabular}
\end{center}

\smallskip\noindent\emph{Observables, sums and the characteristic function.}
\nopagebreak
\begin{center}\footnotesize
\begin{tabular}{@{}p{3.4cm}@{}p{0.68\textwidth}@{}}
$g_k$, $\hat g_k$ & the observables and their centred versions, \eqref{eq:centering0}\\
$S_n$, $\hS_n$ & the Birkhoff sums $\sum_{k<n}g_k\circ f_{k,0}$ and their centred version, \eqref{eq:center sum}\\
$\sigma_n^2$ & the variance $\bE(\hS_n^2)$, \eqref{eq:variance0}; it is assumed to diverge in all the main statements\\
$F_n$ & the distribution function of $\hS_n/\sigma_n$, \eqref{eq:Fn def}\\
$\lambda$ & the Fourier variable; $\lambda_0$, $\lambda_1$ the thresholds below which $|\lambda|\sigma_n^{-1}$ must lie\\
$L_n$, $C_L$ & the window $L_n=C_L\ln\sigma_n$ separating near from far indices, \eqref{eq:Ln}\\
\end{tabular}
\end{center}

\smallskip\noindent\emph{Constants and conventions.}
\nopagebreak
\begin{center}\footnotesize
\begin{tabular}{@{}p{3.4cm}@{}p{0.68\textwidth}@{}}
$K$, $j_0$ & the constant and the moment order of (O-\ref{cond:multiplier}), also the constant of (O-\ref{cond:lambda_multiplier}), (O-\ref{cond:lambda_multiplier2})\\
$c$ & the exponential rate of \eqref{eq:decay_real} and (O-\ref{cond:lambda_multiplier2}); in Section~\ref{sec:billiard} we write $\nu=e^{-c}$, and there $c$ also denotes the cone constant of \eqref{eq:cone 1}\\
$\alpha_{k,j,\lambda},\ h_{k,j,\lambda},$ $\ell_{k,j,\lambda}$ & the leading eigenvalue, eigenvector and functional of the twisted operators, (O-\ref{cond:lambda_multiplier2})\\
$\const$, $\Const$, $\cO(\cdot)$ & generic constants and order notation, Notation~\ref{not:constants}\\
\end{tabular}
\end{center}

\section{Setting and Main Results}
\label{sec:setting}

Before stating the abstract framework in which we are able to study the Central Limit Theorem for sequential systems, we briefly describe two applications:
smooth expanding maps and dispersing billiards.  This gives the reader an indication of the broad applicability of the present
approach.


\subsection{Sequential Smooth Expanding Maps: A Preview}\label{sec:expanding_p}\ \\
\label{sec:exp preview}
The case of expanding maps has been extensively studied, and our general theory does not provide any new results, apart from the discussion of the growth of the variance. We discuss this case in Section \ref{sec:expanding} only so the reader can compare the present results with the literature.  However, even in this case, we are able to present new results for the problem of studying the variance.

Let $M$ be a smooth compact connected Riemannian manifold.
Let $(f_k)_{k \in \bN} \subset C^2(M,M)$ be a uniform family of expanding maps in the following sense:  $\exists A>0$, $\vartheta >1$ such that for all $k \in \bN$,
\begin{itemize}
  \item $\|Df_k^{-1}\|_\infty \leq \vartheta^{-1}$;
  \item $\|D^2f_k\|_\infty\leq A$.
\end{itemize}
Let $g_k\in C^1(M,\bR)$ with $\sup_{k \in \bN} \| g_k \|_{C^1} \le C$ for some $C>0$.
Given an initial density $\rho\in L^1$, we have that $g_k\circ f_{k-1}\circ \cdot\circ f_0$ is a random variable. We designate its expectation by
\[
\bE(g_k\circ f_{k-1}\circ \cdots\circ f_0)=\int_M g_k\circ f_{k-1}\circ \cdots\circ f_0(x) \rho(x) dx,
\]
and call $\bP$ the associated probability distribution.
 Define the centered observables,
\begin{equation}\label{eq:centered}
\begin{split}
\hat{g}_k & = g_k - \int_M g_k \circ f_{k-1} \circ \cdots \circ f_0\cdot  \rho  \\
\hat{S}_n & = \sum_{k=0}^{n-1} \hat{g}_k \circ f_{k-1} \circ \cdots \circ f_0 \,\, ; \qquad \sigma_n^2 = \mathbb{E}(\hat{S}_n^2)\, ,
\end{split}
\end{equation}
and the cumulative distribution function
\begin{equation}\label{eq:cumulative}
F_n(x) = \mathbb{P}\left( \left\{ \frac{ \hat{S}_n }{\sigma_n} \le x \right\} \right) .
\end{equation}
The above, apart from some minor quirks, is the same model discussed in \cite{CR07,NMTV, H20}.
In Section~\ref{sec:expanding intro} we will prove the following.
\begin{thm}
\label{thm:expanding CLT}
Let  $\rho \in \cC^1(M)$. Then for all $c_\star\in (0,1)$ and all $n\in\bN$ such that $\sigma_n\geq \max\{1,c_\star n^{\frac 13}\ln (n+1)\}$, there exists $\bar C>0$ such that
\[
\left|F_n(x)-\frac{1}{\sqrt{2\pi}}\int_{-\infty}^xe^{-\frac {y^2}2}dy\right|\leq \bar C \sigma_n^{-3}(\ln \sigma_n)^2 n.
\]
\end{thm}

\noindent The above Theorem can be applied only if $\sigma_n$ grows fast enough. As mentioned, in the special case of expanding maps, it is possible to strengthen the result, see \cite{DH24}. However, it is always necessary to meet some requirements on the growth of $\sigma_n$.

 Unfortunately, contrary to the case of a single map, we are not aware of an existing result that provides general checkable criteria to ensure some variance growth. By checkable, we mean that one can ensure the wanted variance growth by computing only averages on a finite number of finite time trajectories.
 
 In Proposition \ref{prop:variance}, we provide such a criterion for the case of smooth expanding maps. Our criterion naturally generalizes the usual criterion for a single map. While the argument is rather general, some technical problems do not allow, at the moment, to extend it to the general setting in which we obtain the CLT.
 
\subsection{Sequential Dispersing Billiards: A Preview}\ \\
We consider a class of dispersing billiard tables formed by a fixed number $\ell$ of closed, convex sets, which we call scatterers, on the
torus $\mathbb{T}^2$.  The boundaries of the scatterers are assumed to have strictly positive curvature.  
The dynamics of the billiard map are defined by the motion of a point particle traveling at unit speed between collisions
and reflecting elastically at collisions with the scatterers.
We consider 
the family of billiard maps $\cF(\tau_*, \cK_*, E_*)$ corresponding to configurations of $\ell$ scatterers which satisfy
\begin{itemize}
  \item The minimum and maximum free flight times between collisions are bounded by $\tau_*$ and $\tau_*^{-1}$, respectively.  
  \item The minimum and maximum curvature of the boundary of the scatterers are bounded by $\cK_*$ and
  $\cK_*^{-1}$, respectively.
  \item The maximum of the $C^3$ norm of the scatterers is bounded by $E_*$.
\end{itemize}
All billiard maps in $\cF(\tau_*, \cK_*, E_*)$ act on the same phase space $M$ and preserve the same invariant
measure $d\musrb = c_0 \cos \vf \, drd\vf$, where $c_0$ is the normalizing constant.

Given this family $\cF(\tau_*, \cK_*, E_*)$, we  prove a sequential Central Limit Theorem
along any {\em $N_{\cF}$-admissible sequence} of maps $( T_i )_{i \in \mathbb{N}} \subset \cF(\tau_*, \cK_*, E_*)$.
See Definition~\ref{def:admissible} for the definition of $N_{\cF}$-admissible sequence.\footnote{Informally,
given $N_{\cF} \in \mathbb{N}$,
an $N_{\cF}$-admissible sequence is a sequence of billiard maps drawn from the family $\cF(\tau_*, \cK_*, E_*)$
which comprise blocks of length $N_{\cF}$ in which all maps in each block are close to one another, but maps in different blocks are not necessarily close.}

For the Central Limit Theorem we consider observables sampled at \emph{each} collision,
rather than once every $N_\cF$ collisions: given $\tilde g_{k,j}\in C^\alpha(M)$ for some
$\alpha\in(0,1/3)$, with $|\tilde g_{k,j}|_{C^\alpha}\leq \tilde K$ for all $k\in\bN$ and
$j<N_{\cF}$, we set $g_k=\sum_{j=0}^{N_{\cF}-1}\tilde g_{k,j}\circ T_j\circ\cdots\circ T_1$,
see \eqref{eq:observe}. Note that the $g_k$ themselves are not H\"older, the billiard maps
being discontinuous; that one can nevertheless sample at every collision rests on
Proposition~\ref{lem:bounded L}, which is not contained in \cite{DL22}.  Given a strictly positive density $\rho \in C^1$,\footnote{The results holds more generally for $\rho\in C_{\bR}$; see
Section~\ref{sec:cone def} for the precise definition of $C_{\bR}$.} $d\nu_0 = \rho d\musrb$, and denote
by $\mathbb{E}$ the expectation with respect to $\nu_0$. We define the centered variables exactly as in \eqref{eq:centered} (where the integrals are w.r.t. $\musrb$) and the cumulative distribution function as in \eqref{eq:cumulative}. Our results can be summarized as follows (see Theorem \ref{thm:billiard CLT} for a precise statement, while the beginning of Section~\ref{sec:billiard} contains a precise definition of the sequential system and Section~\ref{sec:check O1_a} a detailed description of the observables).
\begin{thm}
\label{thm:sample}
Let $\cF(\tau_*, \cK_*, E_*)$ be as defined above and let $(T_i)_{i \in \mathbb{N}} \subset \cF(\tau_*, \cK_*, E_*)$
be an $N_{\cF}$-admissible sequence, and $\ln \rho \in \cC^1$, $\int \rho \, d\musrb=1$. 

Then for all $c_\star \in (0,1)$ there exists $\bar C >0$ such that, for all $n \in \mathbb{N}$ such that $\sigma_n \ge \max \{ 1, c_\star n^{1/3} \ln(n+1) \}$,

\[
\left| F_n(x) - \frac{1}{\sqrt{2\pi}} \int_{-\infty}^x e^{-\frac{y^2}{2} } \, dy \right| \le \bar{C} \sigma_n^{-3} (\ln \sigma_n)^2 n .
\]
\end{thm}
Note that in the non-sequential case, when $T_k=T$ and $g_k=g$, if $g$ is not a coboundary, then $\sigma_n^2\sim n$. Hence the error term is $\cO(n^{-\frac 12} (\ln \sigma_n)^2)$ which is optimal, apart from the logarithm.

In the sequential case we do not have general criteria to check that $\sigma_n \ge \max \{ 1, c_\star n^{1/3} \ln(n+1) \}$ (apart from the case in which the billiards are extremely hyperbolic, that is the operators $\cL_{T}$ have a large spectral gap\footnote{We thank Dmitry Dolgopyat for pointing this out to us.}).
However, sharper results can be obtained in the random case, see Section \ref{sec:random} where we improve on the usual condition (cf. Remark  \ref{rem:L2cob} for details).

We prove Theorem \ref{thm:sample} in Section~\ref{sec:billiard} (that is, we prove its more precise formulation Theorem \ref{thm:billiard CLT}) as an application of the abstract theory described in the next section. However, checking the conditions of the abstract theory entails a non-trivial amount of work, which is carried out in  Section~\ref{sec:billiard}.


\subsection{Abstract Framework and Central Limit Theorems}
\label{sec:abstract}

Consider a sequence of dynamical systems $f_k:M_k\to M_{k+1}$, $k\in\bN$,\footnote{ We use the convention that $0\in\bN$.} where $M_k$ are compact Riemannian manifolds of uniform diameter, and the $f_k$ are measurable functions with respect to the Borel $\sigma$-algebra. Let $\mu_k$ be a sequence of probability measures such that $(f_k)_*\mu_k\ll\mu_{k+1}$.
For all $g\in L^\infty(M_{k}, \mu_k)$ and $\vf\in L^1(M_{k+1},\mu_{k+1})$ we define the Transfer Operator $\cL_k$ as
\begin{equation}\label{eq:transfer_op}
\int_{M_{k}} g\vf \circ f_k \, d\mu_k=\int_{M_{k+1}} \vf \cL_kg \, d\mu_{k+1}.
\end{equation}
A direct computation shows that, $\mu_{k+1}$-a.s., we have
\[
\cL_k g(x) = \sum_{y \in f_k^{-1}(x)} g(y)/J_k(y),
\]
where $J_k(y) = \frac{d\mu_{k+1}(f_k y)}{d\mu_k(y)}$. Note that, by definition, $\cL_k$ is a positive operator and a contraction as an operator from $L^1(M_k,\mu_k)$ to $L^1(M_{k+1},\mu_{k+1})$.
Next, assume that there exist convex cones $C_k\subset L^1(M_k,\mu_k)$, $C_k\cap -C_k=\emptyset$, with the following properties
\begin{enumerate}[\bf (C-1)]
\item \label{cond:positivity} for each $h\in C_k$ we have $\int_{M_k} h  \, d\mu_k > 0$.
\item \label{cond:cone_cont} $\cL_k (C_{k})\subset C_{k+1}$ and the diameter (with respect to the Hilbert metric, see \cite[Appendix D]{DKL}) of $\cL_k (C_{k})$ in $C_{k+1}$ is uniformly bounded.
\item \label{cond:norm} there exists $\bbe_k\in C_k$ such that for each $h\in V_k:=\operatorname{span}(C_k)\subset L^1(M_k,\mu_k)$ there exists $\cop\in \bR_+$:
\[ 
\cop \bbe_k +h\in C_k.
\]
\end{enumerate}
Without loss of generality, by (C-\ref{cond:positivity}), we can choose the $\bbe_k$ such that
\begin{equation}\label{eq:normalize_e}
\int_{M_k}\bbe_k d\mu_k=1.
\end{equation}
\begin{rem}
Note that the above setting is not the most general possible: e.g. $\cL_k$ could be a more general transfer operator with some weight; condition (C-\ref{cond:cone_cont}) could be weakened to some diameter depending on $k$, provided that it does grow under appropriate control; in the conditions (O-\ref{cond:multiplier}), (O-\ref{cond:lambda_multiplier}) below one could let $K$ grow, moderately, with $k$; and so on. We refrain from such endless generalizations, which the reader can easily work out, if needed, to present the basic idea in its more straightforward form.
\end{rem} 

Due to condition (C-\ref{cond:norm}) we can associate to each $k$ a Banach space $\cB_k$ obtained by the completion of $V_k$ with respect to the norm $\|\cdot\|_k$ associated to the cone $C_k$:
\begin{equation}
\label{eq:real cone norm}
\|h\|_k=\inf\{\cop\in\bR_+\;:\; -\cop\bbe_k\preceq h\preceq \cop \bbe_k\}
\end{equation}
where $h\succeq g$ iff $h-g\in C_k$ (see \cite[Equation D.2.1]{DKL} for more details). By \cite[Lemma D.5]{DKL} it follows that $\cB_k$ is a Banach lattice with the order structure $\succeq$. In particular, if $-g\preceq h\preceq g$, then $\|g\|_k\geq \|h\|_k$. Also $\|\bbe_k\|_k=1$.

Next, note that if $h\succeq g$, i.e. $h-g\in C_k$, by condition (C-\ref{cond:cone_cont}), $\cL_k(h-g)\in C_{k+1}$, so  $\cL_kh \succeq\cL_k g$. That is, $\cL_k$ is a positive, order preserving operator.
Accordingly, for all $\alpha\in\cB_k$ we have $\|\cL_k\alpha\|_{k+1}\leq \|\alpha\|_k\|\cL_k\bbe_k\|_{k+1}$.\footnote{ In fact, by conditions  (C-\ref{cond:cone_cont}) and (C-\ref{cond:norm}), $\cL_k\bbe_k\in V_{k+1}$ so $\|\cL_{k}\bbe_k\|_{k+1}<\infty$.}
In other words, $\cL_k\in L(\cB_k,\cB_{k+1})$.  Moreover, by (C-\ref{cond:positivity}) we have 
\begin{equation}\label{eq:int_bound}
\left|\int_{M_k}\alpha d\mu_k\right|\leq \int_{M_k}\bbe_k \|\alpha\|_kd\mu_k=\|\alpha\|_k.
\end{equation}
In addition, by \cite[Lemma D.4]{DKL} (with the choice $\rho(g)=\int g$), there exists $c>0$ such that, for each $k\geq j$, $h\in\cB_j$, with $\int_{M_j}hd\mu_j=0$, we have
\begin{equation} \label{eq:decay_real}
\left\|\cL_k\cdots\cL_j h\right\|_{k+1}\leq 3 e^{-c(k-j +1 )}\|h\|_{j}.
\end{equation}
 We further assume
\begin{enumerate}[\bf (C-1)]\setcounter{enumi}{3}
\item \label{cond:powerbounded} there exists $C_*>0$ such that, for each $ k\geq j\in\bN$, and $h\in\cB_j$, we have $\|\cL_{k}\cdots \cL_j h\|_{k+1}\leq C_*\|h\|_j$.
 \end{enumerate}
 
To continue, we need to complexify the Banach spaces $\cB_k$ and extend the action of $\cL_k$ to such an extension (to simplify notation, we still call $\cL_k$ the extension).
\begin{lem}\label{lem:complexification} There exists a canonical complex extension $\bB_k$ of $\cB_k$ and the action of $\cL_k$ extends to a bounded operator in $L(\bB_k,\bB_{k+1})$.
Moreover, for  each $k,m\in\bN$, the norm of $\cL_{k+m}\cdots\cL_k:\bB_k\to\bB_{k+m+1}$ is bounded by $\sqrt 2 C_*$.
\end{lem}
The above fact is well known but, for the reader's convenience, we provide the proof in Appendix \ref{app:complexification}.

To keep the notation simple from now on we will use $\|\cdot\|_k$ also for the norm on the complex space, since no confusion can arise.
Note that if $f,g\in \cB_k$ then, by \eqref{eq:int_bound},
\begin{equation}\label{eq:int_bound2}
\left|\int_{M_k} (f+ig)d\mu_k\right|\leq \left|\int_{M_k} f d\mu_k\right|+\left|\int_{M_k} g d\mu_k\right|\leq \|f\|_k+\|g\|_k\leq \sqrt 2\|f+ig\|_k.
\end{equation}
Finally, we introduce the real functions (observables) $\{g_k\}_{k\in\bN}$, $g_k\in L^1(M_k,\mu_k)$ for which we shall prove a sequential Central Limit Theorem.
We define $f_{k,j}=f_{k-1}\circ \dots \circ f_j$, $f_{j,j}(x)=x$ and consider the sum
\[
S_n =\sum_{k=0}^{n-1} g_k\circ f_{k,0}.
\]
We let $x_k=g_k\circ f_{k,0}(x_0)$, and assume that $x_0$ is distributed according to a probability measure $d\nu_0 = \rho d\mu_0$, 
$\rho\in \cB_0$.
We will use the notation $\bE$ to denote the expectation with respect to the initial measure $d\nu_0= \rho d\mu_0$.

Then the $x_k$ are random variables and\footnote{ For now, the elements in the sum could also be infinite. Shortly, we will introduce a hypothesis that makes them finite.}
\[
\bE\left(S_n\right):=\sum_{k=0}^{n-1}\int_{M_0} \rho \cdot g_k\circ f_{k,0} d\mu_0
= \sum_{k=0}^{n-1}\int_{M_k} g_k \cL_{k-1}\cdots \cL_0 \rho\; d\mu_k.
\]
It is then natural to define
\begin{equation}\label{eq:centering0}
\hat g_k=g_k-\int_{M_k} g_k \cL_{k-1}\cdots \cL_0 \rho d\mu_k
\end{equation}
so that
\begin{equation}\label{eq:centering}
\int_{M_k} \hat g_k \cL_{k-1}\cdots \cL_0 \rho d\mu_k=0,
\end{equation}
and
\begin{equation}
\label{eq:center sum}
\hS_n  =\sum_{k=0}^{n-1} \hat g_k\circ f_{k,0}.
\end{equation}
\smallskip
By definition, $\bE\left(\hS_n\right)=0$. It is then natural to define
\begin{equation}\label{eq:variance0}
\sigma_n^2=\bE(\hS_n^2)= \sum_{k,j=0}^{n-1} \int_{M_{n}}\cL_{n-1}\cdots\cL_{k}(\hat g_{k}\cL_{k-1}\cdots \cL_{j}(\hat g_j\cL_{j-1}\cdots \cL_{0}\rho)) \, d\mu_n.
\end{equation}
Next, for all $n \in \bN$ with $\sigma_n>0$, $k \le n$,  $\lambda \in \mathbb{R}$, we define the operators $\cL_{k,\lambda, n}$: for each $h\in L^1(M_k,\bC)$,
\begin{equation}\label{eq:twisted}
\cL_{k,\lambda, n} h=\cL_k (e^{i\sigma_n^{-1}\lambda \hat g_k}h ).
\end{equation}
From now on, for simplicity, we assume $\sigma_n>0$ for all $n\in\bN$.
Note that, for $h\in L^1$ and $\lambda\in\bR$ we have $e^{i\sigma_n^{-1}\lambda\hat g_k}h\in L^1$, so $\cL_{k,\lambda,n}$ is well defined, and a contraction, as an operator in $L^1(M_k,\bC)$.\\
We can now state our condition on the observables $\{g_k\}$.
 \begin{enumerate}[\bf (O-1)] 
 \item \label{cond:multiplier} There exists $\lambda_0, K>0$ and $j_0\in \bN$, $j_0>3$, such that for all  $n\in\bN$, $k\leq n$, and $|\lambda|\leq \lambda_0\sigma_n$, $\cL_{k,\lambda, n}$ extends to an operator in $L(\bB_k,\bB_{k+1})$, and for each $h\in\bB_k$, and $0<j\leq j_0$, we have 
 \[
 \|\cL_{k,\lambda,n} ( g_k^j h)\|_{k+1}\leq K^j\|h\|_k.
 \] 
 \end{enumerate}
In addition, we assume
\begin{enumerate}[\bf(O-1)]\setcounter{enumi}{1}
\item \label{cond:lambda_multiplier} There exists $K,\lambda_0>0$, such that for all 
$j,k,n\in\bN$, $j \le k \le n$, 
$|\lambda|\leq \lambda_0\sigma_n$, and $h\in \cB_{j}$ we have 
$\|\cL_{k,\lambda,n}\cdots \cL_{j,\lambda,n} h\|_{k+1}\leq K \|h\|_j$.
\end{enumerate}

Note that, by \eqref{eq:int_bound}, (O-\ref{cond:multiplier}), with $\lambda=0$, and (C-\ref{cond:powerbounded}) we have
\[
\sup_{k\in\bN}\left|\int_{M_k} g_k \cL_{k-1}\cdots \cL_0 \rho\; d\mu_k\right|=\sup_{k\in\bN}\left|\int_{M_{k+1}} \cL_k g_k \cL_{k-1}\cdots \cL_0 \rho d\mu_{k+1}\right|\leq K C_*\|\rho\|_0.
\]
The above shows that \eqref{eq:centering0} is well defined.
In addition, by \eqref{eq:centering0}, for each $j\leq j_0$ the function $\hat g_k^j$ is a linear combination of the $g_k^i$, $i\leq j$, with coefficients bounded in terms of $K$, $C_*$ and $\|\rho\|_0$, by the above bound. Hence, by (O-\ref{cond:multiplier}) and (O-\ref{cond:lambda_multiplier}), there exists $\hat K>0$, depending only on $K$, $C_*$ and $\|\rho\|_0$, such that, for each $|\lambda|\leq\lambda_0\sigma_n$, $h\in\bB_k$ and $0<j\leq j_0$,
\begin{equation}\label{eq:ghat_bound}
\|\cL_{k,\lambda,n} (\hat g_k^j h)\|_{k+1}\leq \hat K^j\|h\|_k.
\end{equation}
Note also, that the difference between $\cL_{k,\lambda,n}$ defined with the centered 
observable $\hat g_k$ as opposed to the uncentered $g_k$ is a constant factor $e^{ia_k}$ for some
$a_k \in \mathbb{R}$. Thus, the verification of (O-\ref{cond:multiplier}) is independent of 
the choice of $\rho$.

Next, we show that the above assumptions imply $\sigma_n<\infty$.
\begin{lem}\label{lem:variance_up}
There exists $\tilde{C}>0$ such that, for all $n\in\bN$,
\[
\sigma_n^2\le \tilde{C} n.
\]
\end{lem}
\begin{proof}
Indeed, using that $\cL_{j-1} \cdots \cL_0 \rho \in \cB_j$, together with \eqref{eq:int_bound},
(C-\ref{cond:powerbounded})
and \eqref{eq:ghat_bound}, we may estimate
\begin{equation}
\label{eq:var_upper}
\begin{split}
\sigma_n^2 & = \sum_{k=0}^{n-1} \int_{M_{k+1}} \cL_k( {\hat g}^2_k \cL_{k-1} \cdots \cL_0 \rho) \, d\mu_{k+1} \\
& \quad +  2 \sum_{k=1}^{n-1} \sum_{j=0}^{k-1} \int_{M_{k+1}} \cL_{k}(\hat g_{k}\cL_{k-1}\cdots \cL_{j}(\hat g_j\cL_{j-1}\cdots \cL_{0}\rho)) \, d\mu_{k+1} \\
& \le n C_*  \hat K^2 \| \rho \|_0 + 6 C_*  \hat K^2 \| \rho \|_0 \sum_{k=1}^{n-1} \sum_{j=0}^{k-1} e^{-c(k-j-1)} 
\le \tilde{C} n \, ,
\end{split} 
\end{equation} 
where, in the last line, we have used \eqref{eq:decay_real} to estimate
$\| \cL_{k-1}\cdots \cL_{j}(\hat g_j\cL_{j-1}\cdots \cL_{0}\rho) \|_k$, since we have that $\int_{M_{j+1}} \cL_{j}(\hat g_j\cL_{j-1}\cdots \cL_{0}\rho) \, d\mu_{j+1} =0$ and $\cL_{j}(\hat g_j\cL_{j-1}\cdots \cL_{0}\rho) \in \cB_{j+1}$ by condition (O-\ref{cond:multiplier}).
\end{proof}

Finally, we state a stronger assumption (the composition of twisted transfer operators is of Perron-Frobenius type), which will allow us to obtain stronger results:

\begin{enumerate}[\bf(O-1)]\setcounter{enumi}{2}
\item \label{cond:lambda_multiplier2} There exist $\lambda_0, c, K>0$ and elements $h_{k,j,\lambda,n}\in\cB_k$, $\ell_{k,j,\lambda,n}\in\cB_j'$, $n, k,j\in\bN$, $j\leq k\leq n$, such that, for all $|\lambda|\leq \lambda_0\sigma_n$, $\int_{M_k} h_{k,j,\lambda,n}d\mu_k=\ell_{k,j,\lambda,n}(\bbe_j)=1$, $\|h_{k,j,\lambda,n}\|_{k}+ \|\ell_{k,j,\lambda,n}\|_{j}' \leq K$,\footnote{ Where, for each $\ell\in \cB'_j$, we have $\|\ell\|_{j}' =\sup_{\|h\|_j\leq 1}|\ell(h)|$.}  $|\ell_{k,j,\lambda,n}(h_{j,l,\lambda,n})|\geq K^{-1}$, for each $l\leq j$. Moreover, there exist $\alpha_{k,j,\lambda,n}\in\bC$, $|\alpha_{k,j,\lambda,n}|\leq K$ such that for all  $j,k,n\in\bN$, $j \le k \le n$,  and $h\in \cB_{j}$ we have 
\begin{equation}
\label{eq:O3 decay}
\| \cL_{k-1,\lambda, n}\cdots \cL_{j,\lambda, n} h -\alpha_{k,j,\lambda,n}h_{k,j,\lambda,n}\ell_{k,j,\lambda,n}(h)\|_{k}\leq K |\alpha_{k,j,\lambda,n}|e^{-c(k-j)}\|h\|_j.
\end{equation}
\end{enumerate}
\begin{rem}
Obviously, we can choose the same $K,\lambda_0$ in the assumptions (O-\ref{cond:multiplier})--(O-\ref{cond:lambda_multiplier2}).  Note that (O-\ref{cond:lambda_multiplier2}) implies (O-\ref{cond:lambda_multiplier}).
\end{rem}
\begin{notation}\label{not:constants} 
In order to alleviate notation, from now on we suppress the subscript $n$ in $h_{k,j,\lambda,n}, \ell_{k,j,\lambda,n},\alpha_{k,j,\lambda,n}$, since its value is always clear from the context.\\
In the following, we will use $\const, \Const$ for a generic constant depending only on the constants in Conditions (C-1)--(C-4) and (O-\ref{cond:multiplier})--(O-\ref{cond:lambda_multiplier2}). Also, given a Banach space $\cB$, we will use the notation $\cO(a)$, $a\in\bR_+$, to stand for an arbitrary element of $h\in\cB$ such that $\|h\|\leq \Const a$; we will not specify explicitly the Banach space (which could be $\bR, \bC, \cB_k,\bB_k$ etc ...) since it will always be clear from the context.
\end{notation}
Our main results are the following.
\begin{thm}\label{thm:main}
If the conditions (C-1), (C-2), (C-3), (C-4) and (O-\ref{cond:multiplier}), with $j_0=4$, together with (O-\ref{cond:lambda_multiplier}) are satisfied and $\ln\sigma_n>1$, then for each $\varpi>3$ there exists $C_\varpi, \lambda_1>0$ such that, for each $n\in\bN$ and $|\lambda|\leq \lambda_1\sigma_n$, we have
\[
\left|\bE\left(e^{i\lambda \sigma_n^{-1}  \hS_n}\right)-e^{-\frac{\lambda^2}2}\right|\leq C_\varpi\left(\frac{(\ln \sigma_n)^2\lambda^4}{\sigma_n^2}+ n^2 \lambda^2 \sigma_n^{-\varpi}+ n\frac{|\lambda|^3(\ln \sigma_n)^2}{\sigma_n^3} \right).
\]
\end{thm}
We will see that the above result yields helpful information only if for each $C_0\in\bR_+$ there exists $n_0\in\bN$ such that $\sigma_n\geq C_0n^{\frac 13}$ for all $n\geq n_0$.
\begin{thm}\label{thm:main_bis}
If the conditions (C-1), (C-2), (C-3), (C-4) and (O-\ref{cond:multiplier}),  with $j_0=4$, together with hypothesis (O-\ref{cond:lambda_multiplier2}) are satisfied,
 then for each $\varpi>3$ there exist constants $C_\varpi, \lambda_1>0$ and  a $C^1$ function $A_n$ such that, for each $n\in\bN$ and $
|\lambda|\leq \lambda_1\sigma_n(\ln \sigma_n)^{-1}$,
we have
\[
\begin{split}
&\bE\left(e^{i\lambda \sigma_n^{-1} \hS_n}\right)=e^{-\frac{\lambda^2}2+A_n(\lambda)}\\
&\left|A_n(\lambda)\right|\leq C_\varpi n\left[ |\lambda|^3\sigma_n^{-3}(\ln\sigma_n)^2 +\sigma_n^{-\varpi}|\lambda|\right]\\
&\left|A_n'(\lambda)\right|\leq  C_\varpi n\left[ \lambda^2\sigma_n^{-3}(\ln\sigma_n)^2 +\sigma_n^{-\varpi}\right].
\end{split}
\]
\end{thm}
The proof of Theorem~\ref{thm:main} is the content of Section \ref{sec:real_PF}, while 
Theorem~\ref{thm:main_bis} is proven in Section \ref{sec:complex_PF}.

To compare with other results, let us compute what the above results imply for the distribution function
\begin{equation}
\label{eq:Fn def}
F_n(x)=\bP \left( \left\{ \frac{\hS_n}{\sigma_n} \leq x \right\} \right) = \bE \left( \Id_{ \{ \hat{S}_n/\sigma_n \le x \} } \right) .
\end{equation}
Our first result using (O-\ref{cond:lambda_multiplier}) is nontrivial only if $\sigma_n\gg n^{\frac 13}$. This is the same restriction in \cite{H20}. But there, the error is smaller.

\begin{cor}
If the conditions (C-1),(C-2),(C-3), (C-4) and (O-\ref{cond:multiplier}),  with $ j_0=4$, together with hypothesis (O-\ref{cond:lambda_multiplier}) are satisfied,
then for all $n\in\bN$ such that $\ln\sigma_n>1$ we have
\begin{equation}\label{eq:final_est1}
\left|F_n(x)-\frac{1}{\sqrt{2\pi}}\int_{-\infty}^xe^{-\frac {y^2}2}dy\right|\leq \Const\sigma_n^{-\frac 34}n^{\frac 14}(\ln \sigma_n)^{\frac 12},
\end{equation}
\end{cor}
\begin{proof}
We can assume $\sigma_n^{-\frac 34}n^{\frac 14}(\ln\sigma_n)^{\frac 12}\leq 1$, otherwise the statement is trivial; hence, $n\leq \sigma_n^{3}$.
By \cite[equation (3.13) of Chapter XVI.3]{Fe2}, Theorem \ref{thm:main} implies
\[
\begin{split}
\left|F_n(x)-\frac{1}{\sqrt{2\pi}}\int_{-\infty}^xe^{-\frac {y^2}2}dy\right|&\leq \frac {C_\varpi}\pi\int_{-T_n}^{T_n}\left[\frac{(\ln \sigma_n)^2\zeta^3}{\sigma_n^2}+
 n\frac{\zeta^2(\ln \sigma_n)^2}{\sigma_n^3} + n^2 \sigma_n^{-\varpi}\zeta\right] d\zeta+ \frac{24 }{\pi T_n}\\
&\leq \Const C_\varpi \left[\frac{(\ln \sigma_n)^2T_n^4}{\sigma_n^2}+ n\frac{T_n^3(\ln \sigma_n)^2}{\sigma_n^3} + n^2\sigma_n^{-\varpi}T_n^2 \right]+\frac{24 }{\pi T_n},
\end{split}
\]
for each sequence $\{T_n\}$. Choosing $T_n= \sigma_n^{\frac 34}n^{-\frac 14}(\ln \sigma_n)^{-\frac 12}$, we have $T_n\geq 1$ and recalling that \eqref{eq:var_upper} implies $\sigma_n\leq \Const \sqrt n$, the result follows provided we choose $\varpi\geq 7$.
\end{proof}
To obtain a stronger result, we must assume  (O-\ref{cond:lambda_multiplier2}) and apply Theorem \ref{thm:main_bis}. 
\begin{cor}
\label{cor:better corollary}
If the conditions (C-1),(C-2),(C-3), (C-4) and (O-\ref{cond:multiplier}),  with $j_0=4$, together with hypothesis (O-\ref{cond:lambda_multiplier2}) are satisfied, then for each $c_\star\in (0,1)$ and for all $n\in\bN$ such that $\sigma_n\geq \max\{1,c_\star n^{\frac 13}\ln (n+1)\}$, there exists $\bar C>0$ such that
\begin{equation}\label{eq:final_est2}
\left|F_n(x)-\frac{1}{\sqrt{2\pi}}\int_{-\infty}^xe^{-\frac {y^2}2}dy\right|\leq \bar C \sigma_n^{-3}(\ln \sigma_n)^2 n.
\end{equation}
\end{cor}
\begin{proof}
We use again \cite[equation (3.13) of Chapter XVI.3]{Fe2} to write
\[
\left|F_n(x)-\frac{1}{\sqrt{2\pi}}\int_{-\infty}^xe^{-\frac {y^2}2}dy\right|\leq \frac {1}\pi\int_{-T_n}^{T_n}\left|\frac{e^{-\frac{\zeta^2}2}-e^{-\frac{\zeta^2}2+A_n(\zeta)}}\zeta \right| d\zeta+ \frac{24 }{\pi T_n}.
\]
We make the choice $T_n=\frac{\sigma_n^3}{C_\star n(\ln \sigma_n)^2}$, for some constant $C_\star>0$ large enough, and $\varpi=4$. By Theorem \ref{thm:main_bis} it follows that 
\[
\begin{split}
&\frac{\zeta^2}4\geq C_\varpi n\left| \zeta^3\sigma_n^{-3}(\ln\sigma_n)^2\right|\\
&\left|\frac{d}{d\zeta} e^{-\frac{\zeta^2}2+A_n(\zeta)}\right|\leq C_0,
\end{split}
\]
for all $|\zeta|\leq T_n$.
Also, setting $S_n=\frac{\sigma_n}{n^{\frac 13}(\ln\sigma_n)^{\frac 23}}$, we have that $|A_n(\zeta)|$ is uniformly bounded for all $|\zeta|\leq S_n$. Hence,  for all $|\zeta|\leq S_n$ we have, by Theorem \ref{thm:main_bis},
\[
\begin{split}
\left|\frac{1-e^{A_n(\zeta)}}{\zeta}\right|&\leq \Const \left|\frac{\int_0^\zeta A'(z) dz}{\zeta}\right|\leq \Const\sup_{z\in[0,\zeta]} |A'(z)|\\
&\leq \Const C_\varpi n\left| \zeta^2\sigma_n^{-3}(\ln\sigma_n)^2 +\sigma_n^{-\varpi}\right|.
\end{split}
\]
Accordingly,
\[
\begin{split}
\left|F_n(x)-\frac{1}{\sqrt{2\pi}}\int_{-\infty}^xe^{-\frac {y^2}2}dy\right| \leq&  \frac {1}\pi\int_{-S_n}^{S_n}e^{-\frac{\zeta^2}2}C_\varpi n\left[ \zeta^2\sigma_n^{-3}(\ln\sigma_n)^2 +\sigma_n^{-4}\right] d\zeta\\
&+ \frac {2e^{C_\varpi c_\star^{-\varpi}}}\pi\int_{S_n}^{T_n}e^{-\frac{\zeta^2}4}+\frac{24 }{\pi T_n}\\
&\leq  \bar C \sigma_n^{-3}(\ln \sigma_n)^2 n,
\end{split}
\]
for some $\bar C$ large enough.
\end{proof}
We conclude with a brief discussion on the relation with previous results.
\begin{rem}\label{rem:momenta}
The above results are only slightly weaker than \cite{H20}. However, they apply to a much wider range of systems. In particular, they are tailored for possible applications to the problem of a random Lorenz gas where the observable may not the bounded (e.g., see the relevant observables in the recent \cite{DK23}).  For bounded observables, one can take advantage of the fact that 
$\bE \left( e^{i \lambda \sigma_n^{-1} \hat{S}_n} \right)$ 
is an analytic function of $\lambda$. More generally, our results can be improved by following the same strategy used here and simply computing more terms in the various Taylor expansions,
at the price of requiring a larger value of $j_0$ in Condition (O-\ref{cond:multiplier}).
To obtain optimal results, it is however necessary to have bounds on the higher moments in terms of the variance.  The latter strategy is used in the recent paper \cite{DH24}, where much stronger results are obtained but limited to bounded observables and expanding maps. Unfortunately, it is not obvious how to extend the martingale approximation argument used to control the higher momenta in \cite{DH24} to the present general setting.
\end{rem}

\section{Proof of Theorem \ref{thm:main}}\label{sec:real_PF}
In order to compute
\begin{equation}\label{eq:upsilon}
\Upsilon_n(\lambda)=\bE\left(e^{i\lambda \sigma_n^{-1} \hS_n}\right)
\end{equation}
we show that it satisfies a differential inequality.
\begin{prop}\label{prop:diff1}
For each $\varpi>0$, there exist $C_\varpi,\lambda_1>0$ such that for all $n\in\bN$ and $\lambda\in\bR$, $|\lambda|\leq \lambda_1\sigma_n$, we have
\[
\begin{split}
&\frac{d}{d\lambda}\Upsilon_n(\lambda)=-\lambda \Upsilon_n(\lambda)+\cE_n(\lambda)\\
&|\cE_n(\lambda)|\leq C_\varpi\left(\frac{(\ln \sigma_n)^2  |\lambda|^3 }{\sigma_n^2}+ n^2 |\lambda|\sigma_n^{-\varpi}+ n\frac{\lambda^2(\ln \sigma_n)^2}{\sigma_n^3} \right).
\end{split}
\]
\end{prop}
Proposition \ref{prop:diff1} immediately implies Theorem \ref{thm:main}.
\begin{proof}[\bf \em  Proof of Theorem \ref{thm:main} using Proposition~\ref{prop:diff1}]\ \\
Setting $\theta_n(\lambda) = e^{\frac{\lambda^2}2} \Upsilon_n(\lambda)$, we compute 
\[
\theta_n'(\lambda)=e^{\frac{\lambda^2}2}\cE_n(\lambda).
\]
Since $\Upsilon_n(0)=1$, we have $\theta_n(0)=1$, hence
\[
\Upsilon_n(\lambda)=e^{-\frac{\lambda^2}2}+\int_0^\lambda e^{-\frac{\lambda^2-\xi^2}2}\cE_n(\xi)d\xi.
\]
Moreover, using that $e^{- (\lambda^2 - \xi^2)/2} \le 1$, we have
\[
\left|\int_0^\lambda e^{-\frac{\lambda^2-\xi^2}2}\cE_n(\xi)d\xi\right|\leq\sup_{|\xi|\leq|\lambda|}{|\cE_n(\xi)| }\,  |\lambda| \leq C_\varpi\left(\frac{(\ln \sigma_n)^2\lambda^4}{\sigma_n^2}+ n^2 \lambda^2 \sigma_n^{-\varpi}+  n\frac{|\lambda|^3(\ln \sigma_n)^2}{\sigma_n^3} \right).
\]
\end{proof}
For the remainder of this section, we fix $n \in \mathbb{N}$ and prove Proposition~\ref{prop:diff1}.  Since $n$ is fixed, to simplify
notation a little, we will drop the last index $n$ in the definition of $\cL_{k, \lambda, n}$ from \eqref{eq:twisted}.  That is,
\begin{equation}
\label{eq:simpler}
\cL_{k , \lambda} h := \cL_{k, \lambda, n} h = \cL_k(e^{i \sigma_n^{-1} \lambda {\hat g_k } } h) \, .
\end{equation}

\begin{proof}[\bf \em Proof of Proposition \ref{prop:diff1}]\ \\
For $\mu_0$-almost every $x$ we have
\[
\frac{d}{d\lambda}e^{i\lambda \sigma_n^{-1} \hS_n(x)}=i\sigma_n^{-1}\sum_{k=0}^{n-1}\hat g_k\circ f_{k,0}(x)e^{i\lambda \sigma_n^{-1} \hS_n(x)}.
\]
By the Schwartz inequality, recalling \eqref{eq:transfer_op}, \eqref{eq:int_bound},  \eqref{eq:ghat_bound} and (C-\ref{cond:powerbounded}), we have
\begin{equation}\label{eq:basic_bound}
\begin{split}
\bE(|\hat g_k\circ f_{k,0}(x)|)&\leq \left[\int_{M_{k+1}}\cL_k\left[(\hat g_k)^2\cL_{k-1}\cdots\cL_0\rho \right] d\mu_{k+1}\right]^{\frac 12}\\
&\leq \hat K C_*^\frac 12 \|\rho\|_0<\infty.
\end{split}
\end{equation}
Accordingly, by Leibniz's integral rule,
\begin{equation}\label{eq:step1}
\frac{d}{d\lambda}\Upsilon_n(\lambda)=i\sigma_n^{-1}\sum_{k=0}^{n-1}\bE\left(\hat g_k\circ f_{k,0}e^{i\lambda \sigma_n^{-1} \hS_n}\right).
\end{equation}
For some constant $C_L>0$, to be chosen later large enough, we define
\begin{equation}\label{eq:Ln}
L_n=C_L\ln\sigma_n.
\end{equation}
 Next, for each $k\leq n$ let
 \[
 \begin{split}
&\hS^i_{n,k}=\sum_{\substack{ j\in\{0,\dots, n-1\}\\ |j-k|\leq L_n}} \hat g_j\circ f_{j,0}\\
&\hS^e_{n,k}=\sum_{\substack{ j\in\{0,\dots, n-1\}\\ |j-k|>  L_n}} \hat g_j\circ f_{j,0}.
 \end{split}
 \]
Also, we define
 \[
 \begin{split}
& \Theta^0_k(\xi,\lambda)=\bE\left(e^{i\lambda \sigma_n^{-1}\hS^e_{n,k}+i\xi\sigma_n^{-1}\hS^i_{n,k} }\right)\\
& \Theta^1_{k}(\xi,\lambda)=\bE\left(\hat g_k\circ f_{k,0}e^{i\lambda \sigma_n^{-1}\hS^e_{n,k}+i\xi\sigma_n^{-1}\hS^i_{n,k} }\right)\\
& \Theta^2_{k,j}(\xi,\lambda)=\bE\left(\hat g_k\circ f_{k,0}\hat g_j\circ f_{j,0}e^{i\lambda \sigma_n^{-1}\hS^e_{n,k}+i\xi\sigma_n^{-1} \hS^i_{n,k} }\right)\\
& \Theta^3_{k,j,l}(\xi,\lambda)=\bE\left(\hat g_k\circ f_{k,0}\hat g_j\circ f_{j,0}\hat g_l\circ f_{l,0}e^{i\lambda \sigma_n^{-1}\hS^e_{n,k}+i\xi\sigma_n^{-1} \hS^i_{n,k} }\right).
\end{split}
\]
\begin{lem}\label{lem:taylor}
Let $j_*=j_0$ if $j_0$ is even, and $j_*=j_0-1$ if $j_0$ is odd.
 For each $\lambda \in\bR$ and $i\leq j_*$, the functions $\Theta^i(\cdot,\lambda)$ are everywhere $(j_*-i)$-times differentiable with bounded derivatives.
\end{lem}
\begin{proof}
For $\mu_0$-almost every $x$, the integrand is analytic. Its  $s-i$ derivative, for $s\leq j_*$, is bounded by a sum of terms of the type
\[
|\hat g_{a_1}\circ f_{a_1,0}\cdots \hat g_{a_s}\circ f_{a_s,0}|
\]
for a collection $\{a_i\}\subset \{0,\dots,n-1\}$. Suppose for simplicity that $a_i<a_{i+1}$, the case in which some $a_i$ are equal can be treated similarly. Let $ \varsigma=s$ if $s$ is even and $s+1$ if $s$ is odd.  Then, $\hat g_{a_i}\circ f_{a_i,0}\in L^\varsigma$, since, by  \eqref{eq:ghat_bound} and (C-\ref{cond:powerbounded}),
\[
\bE\left(\hat g_{a_i}^\varsigma\circ f_{a_i,0}\right)=\int_{M_{a_i+1}}\left( \cL_{a_i}\left[(\hat g_{a_{i}})^\varsigma\cL_{a_i-1}\cdots\cL_0\rho \right]\right)d\mu_{a_{i}+1}\leq C_*  \hat K^\varsigma \|\rho\|_0.
\]
Accordingly,
\[
\begin{split}
\bE&\left(|\hat g_{a_1}\circ f_{a_1,0}\cdots \hat g_{a_s}\circ f_{a_s,0}|\right)\leq
\bE\left(\hat g_{a_s}^\varsigma\circ f_{a_s,0}\right)^{\frac 1  \varsigma}
\bE\left(|\hat g_{a_1}\circ f_{a_1,0}\cdots \hat g_{a_{s-1}}\circ f_{a_{s-1},0}|^{\frac{\varsigma}{\varsigma-1}}\right)^{\frac{\varsigma-1}\varsigma}\\
\leq &[C_* \hat K^\varsigma\|\rho\|_0]^{\frac 1\varsigma}\bE\left(\hat g_{a_{s-1}}^\varsigma\circ f_{a_{s-1},0}\right)^{\frac 1  \varsigma}
\bE\left(|\hat g_{a_1}\circ f_{a_1,0}\cdots \hat g_{a_{s-2}}\circ f_{a_{s-2},0}|^{\frac{\varsigma}{\varsigma-2}}\right)^{\frac{\varsigma-2}{\varsigma}}\\
\leq &\left[C_* \|\rho\|_0\right]^{\frac s\varsigma} \hat K^s,
\end{split}
\]
where the last inequality is obtained by induction, using Hölder's inequality.
The lemma follows by repeated application of Leibniz's Integral Rule since all the derivatives are bounded by an $L^1$ function.
\end{proof}

By \eqref{eq:step1} and Lemma \ref{lem:taylor} we can Taylor expand with respect to $\xi$ and write
\begin{equation}\label{eq:step2}
\begin{split}
\frac{d}{d\lambda}\Upsilon_n(\lambda)=&i\sigma_n^{-1}\sum_{k=0}^{n-1}\Theta^1_{k}(\lambda,\lambda)\\
=&i\sigma_n^{-1}\sum_{k=0}^{n-1}\Theta^1_{k}(0,\lambda)-\sigma_n^{-2}\lambda\sum_{k=0}^{n-1}
\sum_{\substack{ j\in\{0,\dots, n-1\}\\ |j-k|\leq L_n}}\Theta^2_{k,j}(0,\lambda)\\
&-i\sigma_n^{-3}\int_0^\lambda d\xi\int_0^\xi dz \sum_{k=0}^{n-1}\sum_{\substack{ j\in\{0,\dots, n-1\}\\ |j-k|\leq L_n}}\sum_{\substack{ l\in\{0,\dots, n-1\}\\ |l-k|\leq L_n}}\Theta^3_{k,j,l}(z,\lambda).
\end{split}
\end{equation}
\begin{rem}
Note that nothing prevents us from expanding $\Theta^1_k$ to higher orders,
although this would require $j_0$ to be correspondingly larger in (O-\ref{cond:multiplier}). Doing so yields better estimates, but at the price of a much lengthier computation. We refrain from doing so as our goal is to present clearly the idea rather than to state the optimal results.
\end{rem}
Next, we must compute the terms in the above equation. To this end, we use the twisted transfer operators  $\cL_{k, \lambda}$ (recall \eqref{eq:simpler}).
Let us define, for all $k,j\in\bZ$,
\begin{equation}\label{eq:transfer_gen}
\cL_{k,j,\lambda}=\begin{cases}
\Id &\textrm{ if } k<0\\
\cL_{k,\lambda}\cdots\cL_{\max\{j,0\},\lambda} &\textrm{ if } k\geq j\\
\Id &\textrm{ if } k< j.
\end{cases}
\end{equation}
Note that, for $h\in L^1$,
\[
\int_{M_{k+1}}\left|\cL_{k,\lambda} h\right|\leq \int_{M_{k+1}}\cL_{k} \left|h\right|=\int_{M_k}|h|.
\]
Hence, by \eqref{eq:int_bound2}, (O-\ref{cond:multiplier})  (see also \eqref{eq:ghat_bound}), and (O-\ref{cond:lambda_multiplier}), we have 
\[
\left|\int_{M_{k+1}}\cL_{k,\lambda}\cdots \cL_{j, \lambda}\hat g_j\cL_{j-1,\lambda}\cdots \cL_{0,\lambda} \rho\right|\leq \sqrt 2  \hat K K\|\rho\|_0.
\]
Note that, for all $L_n \le k \le n-1- L_n$, 
 \begin{equation}\label{eq:Theta_PF}
 \begin{split}
& \Theta^0_k(\xi,\lambda)=\int_{M_n}\!\!\!\!\!\!\cL_{ n-1,k+L_n+1,\lambda}\cL_{k+L_n,k-L_n,\xi}\,\cL_{ k-L_n-1,0,\lambda} \rho\\
& \Theta^1_{k}(\xi,\lambda)=\int_{M_n}\!\!\!\!\!\!\cL_{ n-1,k+L_n+1,\lambda}\cL_{k+L_n, k,\xi}\,\hat g_k\cL_{k-1,k-L_n,\xi}\,\cL_{k-L_n-1,0,\lambda} \rho\\
& \Theta^2_{k,j}(\xi,\lambda)=\int_{M_n}\!\!\!\!\!\!\cL_{ n-1,k+L_n+1,\lambda}\cL_{k+L_n, j,\xi}\,\hat g_j\cL_{j-1,k,\xi}\hat g_k\cL_{k-1, k-L_n,\xi}\,\cL_{k-L_n-1,0,\lambda} \rho.
\end{split}
\end{equation}
The above formulae correspond to the case $k+L_n \geq j\geq k$, the definition in the other cases being obvious. 
To estimate the above terms, the following lemmata will be instrumental.
\begin{lem}\label{lem:one}
For each $k, m_1, m_2 \in \bN$ and $h\in\bB_k$, we have
\[
\begin{split}
&\|\cL_{k+m_1+m_2,k+m_2,0}\,\hat g_{k+m_2}\cL_{k+m_2-1,k,0}h\|_{k+m_1+m_2+1}\leq\Const e^{-c\min\{m_1,m_2\}}\|h\|_k\\
&\left|\int_{M_{k+m_1+m_2}}\hskip-12pt\cL_{k+m_1+m_2,k+m_2,0}\,\hat g_{k+m_2}\cL_{k+m_2-1,k,0}h\right|\leq \Const e^{-c m_2}\|h\|_{k}\\
&\|\cL_{k+m_1,k,0}\,\hat g_{k}\cL_{k-1,0,0}\rho\|_{k+m_1+1}\leq\Const e^{-c m_1}
\end{split}
\]
\end{lem}
\begin{proof}
Let us define 
\begin{equation}\label{eq:phi_def}
\phi_{k,j}:=\cL_{k-1,0}\cdots\cL_{j,0}\bbe_j.
\end{equation}
By (C-\ref{cond:powerbounded}), $\|\phi_{k,j}\|_k\leq C_*$ while \eqref{eq:transfer_op} implies $\int_{M_k}\phi_{k,j} d\mu_k=1$.
Equation \eqref{eq:decay_real} yields
\begin{equation}\label{eq:memory}
\begin{split}
\left\|\cL_{k-1,0}\cdots \cL_{j,0}h-\phi_{k,j}\int_{M_j} h d\mu_j\right\|_{k}
&=\left\|\cL_{k-1,0}\cdots \cL_{j,0}\left[h-\bbe_j\int_{M_j} h d\mu_j\right]\right\|_{k}\\
&\leq  \Const e^{-c(k-j)}\|h\|_j.
\end{split}
\end{equation}
In particular, choosing $h=\cL_{j-1,0}\cdots\cL_{0,0}\bbe_0$,
 \begin{equation} \label{eq:project close}
 \left\|\phi_{k,0}-\phi_{k,j}\right\|_{k}\leq  \Const e^{-c(k-j)}.
\end{equation}
 Then, using \eqref{eq:memory} multiple times, conditions (O-\ref{cond:multiplier}), (C-\ref{cond:powerbounded}) and equation \eqref{eq:project close} yields
  \[
 \begin{split}
& \cL_{k+m_1+m_2,k+m_2,0}\,\hat g_{k+m_2}\cL_{k+m_2-1,k,0}h=\phi_{k+m_1+m_2+1,k+m_2+1}\\
&\times\int \cL_{k+m_2,0}\hat g_{k+m_2}\cL_{k+m_2-1,k,0}h
+\cO(e^{-cm_1})\|\cL_{k+m_2,0}\hat g_{k+m_2}\phi_{k+m_2,k}\int h\|_{k+m_2+1}\\
&=\phi_{k+m_1+m_2+1,k+m_2+1}\int \cL_{k+m_2,0}\hat g_{k+m_2}\cL_{k+m_2-1,k,0}h
+\cO(e^{-cm_1}+e^{-c m_2})\|h\|_{k}\\
&=\phi_{k+m_1+m_2+1,0}\int \cL_{k+m_2,0}\hat g_{k+m_2}\cL_{k+m_2-1,k,0}h
+\cO(e^{-cm_1}+e^{-c m_2})\|h\|_{k}\\
&= \phi_{k+m_1+m_2+1,0} \int \cL_{k+m_2,0}\hat g_{k+m_2}  \cL_{k+m_2-1, 0,0} \rho  \int_{M_k} h+\cO(e^{-cm_1}+e^{-cm_2})\|h\|_k\\
&= \phi_{k+m_1+m_2+1,0}\int \hat g_{k+m_2}  \cL_{k+m_2-1, 0,0} \rho  \int_{M_k} h+\cO(e^{-cm_1}+e^{-cm_2})\|h\|_k\\
&=\cO(e^{-c\min\{m_1,m_2\}})\|h\|_k,
 \end{split}
 \]
 where, in the third equality, we have used again \eqref{eq:memory}  to estimate
 \[
 \begin{split}
&\left\| \cL_{k+m_2-1,k,0}h- \cL_{k+m_2-1, 0,0} \rho  \int_{M_k} h\right\|_{k+m_2}
=\left\| \cL_{k+m_2-1,k,0}\left(h- \cL_{k-1, 0,0} \rho  \int_{M_k} h\right)\right\|_{k+m_2}\\
&\leq\left\|\phi_{k+m_2,k}\int \left(h- \cL_{k-1, 0,0} \rho  \int_{M_k} h\right)\right\|_{k+m_2}+\Const e^{-cm_2}\left\|h- \cL_{k-1, 0,0} \rho  \int_{M_k} h\right\|_k\\
&\leq \Const e^{-cm_2}\|h\|_k.
\end{split}
 \]
This proves the first inequality. The other two are proven similarly.
\end{proof}
Let
 \begin{equation}\label{eq:Gamma_PF}
 \begin{split}
& \Gamma^0_k=\cL_{k+L_n,  k-L_n,0}\\
& \Gamma^1_{k}=\cL_{k+L_n, k,0}\,\hat g_k\cL_{k-1, k-L_n ,0}\\
& \Gamma^2_{k,j}=\cL_{k+L_n, j,0}\,\hat g_j\cL_{j-1,k,0}\hat g_k\cL_{k-1, k-L_n ,0}.\\
\end{split}
\end{equation}
\begin{lem}\label{lem:gamma_est}
For each $\varpi>0$ there exists $C_L=C_L(\varpi)>0$ such that for all $n-1-L_n \geq k\geq L_n$ and $h\in\cB_{k-L_n}$, we have
\[
\begin{split}
&\|\Gamma^1_{k}h\|_{ k+L_n+1}\leq \Const \sigma_{ n}^{-\varpi}\|h\|_{k-L_n}\\
&\|\Gamma^2_{k,j}h- \Sigma_{k, j}\Gamma_k^0h\|_{ k+L_n+1} \leq\Const \sigma_n^{-\varpi}\|h\|_{k-L_n},
\end{split}
\]
where
\begin{equation}\label{eq:variance}
\Sigma_{k,j}=\int \hat g_{k}\cL_{k-1,j,0}\hat g_j\cL_{j-1,0,0}\rho.
\end{equation}
\end{lem}
\begin{proof}
Equations \eqref{eq:decay_real}, \eqref{eq:int_bound} and conditions (C-\ref{cond:powerbounded}), (O-\ref{cond:multiplier}) imply
\[
\begin{split}
\left\|\Gamma^1_{k}h\right\|_{k+L_n+1} &\leq \left\|  \phi_{k+L_n+1,k+1}\int_{M_{k+1}}\cL_{k}\hat g_k\cL_{k-1, k-L_n,0}hd\mu_{k+1} \right\|_{ k+L_n+1} + \Const e^{-c L_n}\|h\|_{k-L_n}\\
&\leq \Const \left|\int_{M_k}\hat g_k\cL_{k-1, 0,0}\rho d\mu_{ k}\right|\cdot\left| \int_{M_{k-L_n}} h d\mu_{k-L_n}\right|+\Const K\sigma_n^{-\varpi}\|h\|_{k-L_n}\\
&\leq \Const \sigma_n^{-\varpi}\|h\|_{k-L_n}
\end{split}
\]
by \eqref{eq:centering} and provided we have chosen $C_L$ large enough. Next, we study $\Gamma^2_{k,j}$ for $j\geq k$, the case $j<k$ being identical. If $|k-j|\geq L_n/2$, then
\begin{equation}\label{eq:gamma1}
\begin{split}
&\Gamma^2_{k,j}h=\cL_{k+L_n, j,0}\,\hat g_j\cL_{j-1,k,0}\hat g_k\cL_{k-1, k-L_n,0}h\\
&=\cL_{k+L_n, j,0}\,\hat g_j\phi_{j,k+1}\int_{M_{k+1}}\cL_{k}\hat g_k\cL_{k-1,k-L_n,0}h+\cO(e^{-c|k-j|}\|h\|_{k-L_n})\\
&=\cO(\sigma_n^{-\varpi}\|h\|_{k-L_n})
\end{split}
\end{equation}
provided we choose $C_L$ large enough. On the other hand, if $|k-j|\leq L_n/2$, then $k+L_n-j\geq L_n/2 $ and $j-k+L_n\geq L_n/2$, hence, recalling \eqref{eq:memory},
\begin{equation}\label{eq:gamma2}
\begin{split}
\Gamma^2_{k,j}h=&\phi_{k+L_n,j+1}\int \cL_j\hat g_j\cL_{j-1,k,0}\hat g_k\cL_{k-1,k-L_n,0}h\\
&+\cO(e^{-c|k+L_n-j|}\|h\|_{k-L_n})\\
&=\phi_{k+L_n,j+1}\int \cL_j\hat g_j\cL_{j-1,k,0}\hat g_k\cL_{k-1,0,0}\rho \int h+\cO(\sigma_n^{-\varpi}\|h\|_{k-L_n})\\
&=\Sigma_{k, j}\Gamma_k^0h+\cO(\sigma_n^{-\varpi}\|h\|_{k-L_n}).
\end{split}
\end{equation}
The Lemma follows noting that, for $|k-j|\geq L_n/2$, $|\Sigma_{k,j}|\leq \Const \sigma_{n}^{-\varpi}$.
\end{proof}

We can now estimate the terms in \eqref{eq:step2} one at a time.

Note first that, by \eqref{eq:centering}, $\Theta^1_k(0,0)=\int_{M_k}\hat g_k\,\cL_{k-1,0,0}\rho\, d\mu_k=0$. It is thus convenient to estimate the $\lambda$-derivative of $\Theta^1_k$, rather than $\Theta^1_k$ itself. By Lemma \ref{lem:taylor} we can differentiate $\Theta^1_k$ with respect to $\lambda$ obtaining
\[
\partial_\lambda\Theta^1_k(0,\lambda)=i\sigma_n^{-1}\hskip-10pt\sum_{\substack{ j\in\{0,\dots, n-1\}\\ |j-k|>  L_n}}\hskip-10pt\Theta^2_{k,j}(0,\lambda).
\]
Accordingly,
\[
|\Theta^1_k(0,\lambda)|\leq \sigma_n^{-1}|\lambda| n\sup_{|j-k|>L_n}\sup_{|\lambda_1|\leq |\lambda|}|\Theta^2_{k,j}(0,\lambda_1)|.
\]
Note that, for $|j-k|>L_n$, the factor $\hat g_j$ occurs in one of the $\lambda$-twisted blocks, as in the estimate preceding \eqref{eq:Theta_PF}.
If $L_n\leq k\leq  n-1-L_n$, Lemma \ref{lem:gamma_est} implies
\[
\left| \Theta^2_{k,j}(0,\lambda)\right|\leq  \Const \sigma_n^{-\varpi}.
\]
If $k\geq n-L_n$, then only the indices $j<k-L_n$ occur and, by the second inequality of Lemma \ref{lem:one} and  conditions (O-\ref{cond:multiplier}), (O-\ref{cond:lambda_multiplier})
\[
\left| \Theta^2_{k,j}(0,\lambda)\right|\leq\Const \sigma_n^{-\varpi}\|\cL_{ k-L_n-1,j,\lambda}\hat g_j\cL_{j-1,0,\lambda} \rho\|_{k-L_n}\leq \Const \sigma_n^{-\varpi}.
\]
While, if $k\leq L_n$, then only the indices $j>k+L_n$ occur and, by the equation \eqref{eq:int_bound}, conditions (O-\ref{cond:multiplier}), (O-\ref{cond:lambda_multiplier}) and third inequality of Lemma \ref{lem:one}
\[
\left| \Theta^2_{k,j}(0,\lambda)\right|\leq \Const\|\cL_{k+L_n, k,0}\,\hat g_k\cL_{k-1,0,0} \rho\|_{k+L_n+1}\leq \Const \sigma_n^{-\varpi}.
\]
Collecting the above estimates, and provided $C_L$ has been chosen large enough, for all $k \in \bN$ we have
\begin{equation}\label{eq:theta1}
|\Theta^1_k(0,\lambda)|\leq \Const\,  |\lambda|\, n\sigma_n^{-\varpi}.
\end{equation}

Next, let us compute $\Theta^2_{k,j}(0,\lambda)$ for $|j-k|\leq L_n$. As before, we treat only the case $j\geq k$ since the case $j< k$ is identical and is left to the reader.

For $ n-1-L_n \geq k\geq L_n$  Lemma \ref{lem:gamma_est} implies
\begin{equation}\label{eq:theta2}
\Theta^2_{k,j}(0,\lambda)=\Sigma_{k, j}\Theta^0_k(0,\lambda)+\cO(\sigma_n^{-\varpi}).
\end{equation}
Note that, by equations \eqref{eq:variance0} and \eqref{eq:memory},
\begin{equation}\label{eq:variance_tot}
\sigma_n^2=\sum_{k=0}^{n-1}\sum_{|j-k|\leq L_n}\Sigma_{k,j}+\cO( n \sigma_n^{-\varpi}).
\end{equation}
The other possibilities can be treated as we did for $\Theta^1_k$ and yield a contribution of order $\sigma_n^{-\varpi}$.
Finally, we compute
\begin{equation}\label{eq:rem_tay3}
\begin{split}
& \left|\sum_{k=0}^{n-1}\sum_{\substack{ j,l\in\{0,\dots, n-1\}\\ |j-k|\leq L_n\\ |l-k|\leq L_n}}\Theta^3_{k,j,l}(z,\lambda)\right|\leq \sum_{k=0}^{n-1}\sum_{\substack{ j,l\in\{0,\dots, n-1\}\\ |j-k|\leq L_n\\ |l-k|\leq L_n}}\bE\left(\left|\hat g_k\circ f_{k,0}\,\hat g_j\circ f_{j,0}\,\hat g_l\circ f_{l,0}\right|\right) \\
& \phantom{\sum_{k=0}^{n-1}}\leq \Const\, n\,(2L_n+1)^2\leq \Const\, n L_n^2
\end{split}
\end{equation}
where we have used the H\"older estimate in the proof of Lemma \ref{lem:taylor} (case $s=3$, whence $\varsigma=4$), which yields $\bE\left(\left|\hat g_{a_1}\circ f_{a_1,0}\hat g_{a_2}\circ f_{a_2,0}\hat g_{a_3}\circ f_{a_3,0}\right|\right)\leq \left[C_*\|\rho\|_0\right]^{\frac 34} \hat K^3$.

To conclude, note that
\begin{equation}\label{eq:theta0}
\Upsilon_n(\lambda)=\Theta^0_k(\lambda,\lambda)=\Theta^0_k(0,\lambda)+\partial_\xi\Theta^0_k(0,\lambda)\lambda+\int_0^\lambda d\eta\int_0^\eta dz\partial_\xi^2\Theta^0_k(z,\lambda).
\end{equation}
Arguing as before, we have
\begin{equation}\label{eq:theta0_bis}
\begin{split}
&|\partial_\xi\Theta^0_k(0,\lambda)|\leq C\sigma_n^{-\varpi}\\
&|\partial_\xi^2\Theta^0_k(0,\lambda)|\leq \Const C_L^2\sigma_n^{-2}(\ln \sigma_n)^2.
\end{split}
\end{equation}
Hence, we can rewrite \eqref{eq:theta2} as
\begin{equation}\label{eq:theta2_bis}
\Theta^2_{k,j}(0,\lambda)=\Sigma_{k,j}\Upsilon_n(\lambda) +\cO\left(\Sigma_{k,j}\sigma_n^{-2}C_L^2\lambda^2(\ln \sigma_n)^2 + n \sigma_n^{-\varpi}\right).
\end{equation}
Collecting the above computations in \eqref{eq:step2} yields
\begin{equation}\label{eq:collect}
\begin{split}
\frac{d}{d\lambda}\Upsilon_n(\lambda)&=-\sigma_n^{-2}\sum_{\substack{k\in\{0,\cdots,n-1\}\\|k-j|\leq L_n}}\left[\Sigma_{k,j}\lambda \Upsilon_n(\lambda)+\cO\left(\Sigma_{k,j}\sigma_n^{-2}\lambda^3C_L^2(\ln \sigma_n)^2 + n \sigma_n^{-\varpi}\lambda\right)\right]\\
&\phantom{=}+  \cO\left(n\lambda^2\sigma_n^{-3}C_L^2(\ln \sigma_n)^2\right) \\
&=-\lambda \Upsilon_n(\lambda)+\cO\left(\frac{C_L^2(\ln \sigma_n)^2\lambda^3}{\sigma_n^2}+ n^2 \lambda \sigma_n^{-\varpi}+ n\frac{\lambda^2(\ln \sigma_n)^2}{\sigma_n^3} \right),
\end{split}
\end{equation}
and this concludes the proof recalling that $C_L$ is a constant depending on $\varpi$, see Lemma \ref{lem:gamma_est}.
\end{proof}


\section{Proof of Theorem \ref{thm:main_bis}}\label{sec:complex_PF}
We define $\Upsilon_n$ as in equation \eqref{eq:upsilon}.
\begin{prop}\label{prop:diff2}
For each $\varpi>0$, there exists $\lambda_1>0$ such that, for each $n\in\bN$ and $\lambda\in\bR$, $|\lambda|\leq \lambda_1\sigma_n (\ln\sigma_n)^{-1}$, we have
\[
\begin{split}
&\frac{d}{d\lambda}\Upsilon_n(\lambda)=- \left(\lambda+\cE_n(\lambda)\right)\Upsilon_n(\lambda) \\
&\Upsilon_n(0)=1\\
&\left|\cE_n(\lambda)\right|\leq C_\varpi \left(\lambda^2\sigma_n^{-3}(\ln\sigma_n)^2 n+n\sigma_n^{-\varpi}\right).
\end{split}
\]

\end{prop}
Proposition \ref{prop:diff2} immediately implies Theorem \ref{thm:main_bis}.
\begin{proof}[\bf \em  Proof of Theorem \ref{thm:main_bis}]\ \\
Let $A_n(\lambda)= - \int_0^\lambda\cE_n(\xi) d\xi$ then Proposition \ref{prop:diff2} implies $\Upsilon_n(\lambda)=e^{-\frac{\lambda^2}2 + A_n(\lambda)}$.
\end{proof}
To conclude we must thus prove Proposition \ref{prop:diff2}.
\begin{proof}[\bf \em  Proof of Proposition \ref{prop:diff2}]\ \\
We start again our proof by 
\begin{equation}\label{eq:step5.1}
\frac{d}{d\lambda}\Upsilon_n(\lambda)=i\sigma_n^{-1}\sum_{k=0}^{n-1}\int_{M_n}\cL_{n-1 ,k,\lambda}\hat g_k\cL_{k-1,0,\lambda}\rho d\mu_n,
\end{equation}
where we have used the same notation as in \eqref{eq:Theta_PF}.
If $2L_n\leq k\leq n-2L_n$, then, recalling conditions (O-\ref{cond:lambda_multiplier}) and (O-\ref{cond:lambda_multiplier2}) we can write
\begin{equation}\label{eq:complex_1}
\begin{split}
\cL_{ n-1,k,\lambda}&\hat g_k\cL_{k-1,0,\lambda}\rho=\alpha_{n,k+L_n+1,\lambda}\alpha_{ k-L_n,0,\lambda} h_{n,k+L_n+1,\lambda}\\
&\times \ell_{n,k+L_n+1,\lambda}(\cL_{k+L_n,k,\lambda}\hat g_k\cL_{k-1,k-L_n,\lambda}h_{ k-L_n,0,\lambda})\ell_{k-L_n,0,\lambda}(\rho)\\
&+\cO(\sigma_{ n}^{-\varpi}\alpha_{n,k+L_n+1,\lambda}\alpha_{k-L_n,0,\lambda})
\end{split}
\end{equation}
provided we have chosen $C_L$ large enough.

\begin{lem} \label{lem:coalesce} There exists $n_1\in\bN$ such that for each $n\geq k>j>l$, $k-l\geq n_1$, we have
\[
\alpha_{k,l,\lambda}=\alpha_{k,j,\lambda} \alpha_{j,l,\lambda} \ell_{k,j,\lambda}\left(h_{j,l,\lambda}\right)+\cO\left(  \alpha_{k,j,\lambda}\alpha_{j,l,\lambda} e^{-c\min\{|k-j |, |j-l|\}}\right).
\]
\end{lem}
\begin{proof}
By assumption (O-\ref{cond:lambda_multiplier2}) we have, on the one hand
\[
\left|\int_{M_k}\cL_{k-1,l,\lambda}(\bbe_l)d\mu_k-\alpha_{k,l,\lambda}\right|\leq K|\alpha_{k,l,\lambda}|e^{-c|k-l |}.
\]
On the other hand
\[
\begin{split}
\int_{M_k}&\cL_{k-1,l,\lambda}(\bbe_l)d\mu_k=\int_{M_k}\cL_{k-1,j,\lambda}\left(\cL_{j-1,l,\lambda}(\bbe_l)\right)d\mu_k\\
&=\alpha_{k,j,\lambda}\ell_{k,j,\lambda}\left(\cL_{j-1,l,\lambda}(\bbe_l)\right)+\cO\left(\alpha_{k,j,\lambda}e^{-c|k-j |}\left\|\cL_{j-1,l,\lambda}(\bbe_l)\right\|_j\right)\\
&=\alpha_{k,j,\lambda} \alpha_{j,l,\lambda} \ell_{k,j,\lambda}\left(h_{j,l,\lambda}\right)+\cO\left(\alpha_{k,j,\lambda} \alpha_{j,l,\lambda} e^{-c\min\{|k-j |, |j-l|\}}\right).
\end{split}
\]
Hence, if $Ke^{-cn_1}\leq 1/2$, the lemma follows.
\end{proof}

Taylor expanding with respect to $\lambda$, we have
\begin{equation}\label{eq:complex_2}
\begin{split}
\cL_{l,j,\lambda}h=&\cL_{l,j,0}h+i\sigma_n^{-1}\lambda\sum_{m=j}^l\cL_{l,m,0}\,\hat g_m\cL_{m-1,j,0} h\\
& -\sigma_n^{-2}\sum_{m=j}^{l} \int_0^{\lambda}(\lambda-\xi)\cL_{l,m,\xi}\,\hat g_m^2\cL_{m-1,j,\xi} h \,d\xi \\
&-2\sigma_n^{-2}\sum_{j\leq s< m\leq l} \int_0^{\lambda}(\lambda-\xi)\cL_{l,m,\xi}\,\hat g_m\cL_{m-1,s,\xi}\hat g_s\cL_{s-1,j,\xi} h \,d\xi.
\end{split}
\end{equation}
We can then compute
\begin{equation}\label{eq:complex_3}
\begin{split}
\cL_{k+L_n,k,\lambda}&\hat g_k\cL_{k-1,k-L_n,\lambda} h=\cL_{k+L_n,k,0}\hat g_k\cL_{k-1,k-L_n,0 }h\\
& +\, i\sigma_n^{-1}\lambda\sum_{j=k}^{k+L_n}\cL_{k+L_n,j,0}\,\hat g_j\cL_{j-1,k,0}\,\hat g_k\cL_{k-1,k-L_n,0} h \\
& +\, i\sigma_n^{-1}\lambda\sum_{j=k-L_n}^{k-1}\cL_{k+L_n,k,0}\,\hat g_k\cL_{k-1,j,0}\,\hat g_j\cL_{j-1,k-L_n,0} h \\
& +\,\cO(\lambda^2\sigma_n^{-2}L_n^2\|h\|_{k-L_n}),
\end{split}
\end{equation}
where, in the first order term, the two sums correspond to $j\geq k$ and $j<k$: in the
latter case $\hat g_j$ is inserted to the right of $\hat g_k$. We have written the second order
term directly as a remainder, which is all that is used below, and which does not require
$j_0\geq 4$ in (O-\ref{cond:multiplier}).

\begin{lem}\label{lem:apriori} 
For $\lambda\leq \lambda_0\sigma_n$, we have\footnote{ Recall that $\Sigma_{k,j}$ is defined in \eqref{eq:variance}, while $\lambda_0$ refers to condition (O-\ref{cond:lambda_multiplier2}).}
\[
\begin{split}
\Big\|\cL_{k+L_n,k,\lambda}\hat g_k\cL_{k-1,k-L_n,\lambda}h_{ k-L_n,0,\lambda}&-i\sigma_n^{-1}\lambda\!\!\!\sum_{j=k-L_n}^{k+L_n} \!\!\!\Sigma_{k,j}\cL_{k+L_n,k-L_n,\lambda}h_{ k-L_n,0,\lambda}\Big\|\leq\\
&\leq  \Const\left[\lambda^2\sigma_n^{-2}L_n^2+\sigma_n^{-\varpi}\right].
\end{split}
\]
\end{lem}
\begin{proof}
Equation \eqref{eq:complex_3}
\[
\begin{split}
\cL_{k+L_n,k,\lambda}\hat g_k&\cL_{k-1,k-L_n,\lambda}h_{ k-L_n,0,\lambda}
=\cL_{k+L_n,k,\lambda}\hat g_k\cL_{k-1,k-L_n,0}h_{ k-L_n,0,\lambda}\\
&+ i\sigma_n^{-1}\lambda\sum_{j=k}^{k+L_n} \cL_{k+L_n,j,0}\hat g_j\cL_{j-1,k,0}\hat g_k\cL_{k-1,k-L_n,0}h_{k-L_n,0,\lambda}\\
&+ i\sigma_n^{-1}\lambda\sum_{j=k-L_n}^{k-1} \cL_{k+L_n,k,0}\hat g_k\cL_{k-1,j,0}\hat g_j\cL_{j-1,k-L_n,0}h_{ k-L_n,0,\lambda}\\
&+\cO(\left[\lambda^2\sigma_n^{-2}L_n^2\right]\|h_{ k-L_n,0,\lambda}\|_{k-L_n}).
\end{split}
\]
By Lemma \ref{lem:gamma_est} and equation \eqref{eq:complex_2} it follows
\[
\begin{split}
\cL_{k+L_n,k,\lambda}\hat g_k\cL_{k-1,k-L_n,\lambda}h_{ k-L_n,0,\lambda}
=& i\sigma_n^{-1}\lambda\sum_{j=k-L_n}^{k+L_n} \Sigma_{k,j}\cL_{k+L_n,k-L_n,0}h_{ k-L_n,0,\lambda}\\
&+\cO(\left[\lambda^2\sigma_n^{-2}L_n^2+\sigma_n^{-\varpi}\right]\|h_{ k-L_n,0,\lambda}\|_{k-L_n})\\
=&i\sigma_n^{-1}\lambda\sum_{j=k-L_n}^{k+L_n} \Sigma_{k,j}\cL_{k+L_n,k-L_n,\lambda}h_{k-L_n,0,\lambda}\\
&+\cO(\left[\lambda^2\sigma_n^{-2}L_n^2+\sigma_n^{-\varpi}\right]\|h_{k-L_n,0,\lambda}\|_{k-L_n}).
\end{split}
\]
The Lemma follows recalling that, by condition (O-\ref{cond:lambda_multiplier2}), $\|h_{k-L_n,0,\lambda}\|_{k-L_n}\leq K$.
\end{proof}
 To continue, we need some rough estimates.
\begin{lem}\label{lem:rough}
There exists $\lambda_1\in (0,\lambda_0)$ such that, for all $|\lambda|\leq \lambda_1\sigma_n L_n^{-1}$ and $j\leq k\leq n$, $k-j\leq L_n$, we have, for each $h\in\cB_j$,

\begin{enumerate}
 \item[a)] $\displaystyle |\alpha_{k,j,\lambda}|\geq \frac{1}{2(1+K)}$
 \item[b)] $\displaystyle |\ell_{k,j,\lambda}(h)|\geq \const \left|\int_{M_j}h\right| -\Const \left[ |\lambda| \sigma_n^{-1}|k-j|+ e^{-c|k-j|}\right]\|h\|_j$
  \item[c)] $\displaystyle \|h_{k,j,\lambda}-\phi_{k,j}\|_k\leq \Const\left[ |\lambda| \sigma_n^{-1}|k-j|+ e^{-c|k-j|}\right].$
\end{enumerate}
\end{lem}
\begin{proof}
By condition  (O-\ref{cond:lambda_multiplier2}) we have
\[
 \left|\alpha_{k,j,\lambda}-\int_{M_k}\cL_{k-1,j,\lambda} \bbe_j d\mu_k\right|\leq K e^{-c|k-j|}|\alpha_{k,j,\lambda}|.
\]
While \eqref{eq:complex_2} and \eqref{eq:normalize_e} and Condition (O-\ref{cond:multiplier}) imply
\[
\begin{split}
\int_{M_k}\cL_{k-1,j,\lambda}\bbe_j=&\int_{M_k}\cL_{k-1,j,0}\bbe_j+i\sigma_n^{-1}\sum_{m=j}^k\int_0^\lambda d\xi \int_{M_k}\cL_{k-1,m,\xi}\,\hat g_m\cL_{m-1,j,\xi} \bbe_j\\
=& 1+\cO(\sigma_n^{-1}\lambda |k-j|).
\end{split}
\]
Which yields
\begin{equation}\label{eq:alpha_lower}
|\alpha_{k,j,\lambda}|\geq \frac{1-\Const \sigma_n^{-1}\lambda |k-j|}{1+Ke^{-c|k-j|}}.
\end{equation}
The above implies $|\alpha_{k,j,\lambda}|\geq(2+2K)^{-1}$, provided $|k-j|\leq L_n$ and $\lambda_1$ is chosen small enough.
Next, for each $h\in\cB_j$,
\[
\int_{M_k}\cL_{k-1,j,\lambda} h=\alpha_{k,j,\lambda}\ell_{k,j,\lambda}(h)+\cO(e^{-c|k-j|} \alpha_{k,j,\lambda}\|h\|_j).
\]
Thus, for $|k-j|\leq L_n$, by the first statement of the Lemma and equation \eqref{eq:complex_2} 
\begin{equation}
\label{eq:ell example}
\begin{split}
|\ell_{k,j,\lambda}(h)|&\geq \const\left|\int_{M_k}\cL_{k-1,j,\lambda} h\right|-\Const e^{-c|k-j|}\|h\|_j\\
&\geq \const \left|\int_{M_j}h\right| -\Const \left[ |\lambda| \sigma_n^{-1}|k-j|+ e^{-c|k-j|}\right]\|h\|_j, 
\end{split}
\end{equation}
proving item (b).

Finally, by Conditions (O-\ref{cond:lambda_multiplier2}), (O-\ref{cond:lambda_multiplier}) and equation \eqref{eq:complex_2} we have
\[
\begin{split}
&\alpha_{k,j,\lambda} \left(h_{k,j,\lambda}-\phi_{k,j}\right)=\cL_{k-1,j,\lambda}\bbe_j-\cL_{k-1,j,0}\bbe_j\int\cL_{k-1,j,\lambda}\bbe_j+\cO\left(\alpha_{k,j,\lambda} e^{-c|k-j|}\right)\\
&=\sum_{l=j}^{ k-1}\cL_{k-1,l+1,0}\cL_{l,j,\lambda}\bbe_j\int\cL_{k-1,l+1,\lambda}\bbe_{l+1}-\cL_{k-1,l,0}\cL_{l-1,j,\lambda}\bbe_j\int\cL_{k-1,l,\lambda}\bbe_l\\
&\phantom{=}+\cO\left(\alpha_{k,j,\lambda} e^{-c|k-j|}\right)=\cO\left(\lambda \sigma_n^{-1}|k-j|+\alpha_{k,j,\lambda} e^{-c|k-j|}\right).
\end{split}
\]
The Lemma follows remembering \eqref{eq:alpha_lower}.
\end{proof}
Using Lemma~\ref{lem:apriori} in equation \eqref{eq:complex_1} and recalling Lemma~\ref{lem:coalesce} we obtain, for $2L_n\leq k\leq n-2L_n$, choosing $C_L$ large enough,
\begin{equation}\label{eq:final_complex prelim}
\begin{split}
&\cL_{n-1,k,\lambda}\hat g_k\cL_{k-1,0,\lambda}\rho=i\sigma_n^{-1}\lambda\!\!\! \sum_{j=k-L_n}^{k+L_n} \!\!\!\Sigma_{k,j}\alpha_{n,k+L_n+1, \lambda}\alpha_{ k-L_n,0,\lambda} h_{n,k+L_n+1,\lambda}\\
&\times \ell_{n,k+L_n+1, \lambda}\left(\cL_{k+L_n,k-L_n,\lambda}h_{ k-L_n,0,\lambda}\right)\left[ 1+\cO(\left[\lambda \sigma_n^{-1}L_n+\sigma_n^{-\varpi}\right])\right]\\
&\times\ell_{ k-L_n,0,\lambda}(\rho)+\cO(\sigma_{n}^{-\varpi}\alpha_{n,0,\lambda})\\
&=\left(i\lambda\sigma_n^{-1}\!\!\!\sum_{j=k-L_n}^{k+L_n}\Sigma_{k,j}+\cO\left(\lambda^2\sigma_n^{-2}L_n^2+\sigma_n^{-\varpi}\right)\right)\cL_{ n-1,0,\lambda}\rho\\
&\phantom{=}
+\cO(\sigma_n^{-\varpi}\alpha_{n,0,\lambda}),
\end{split}
\end{equation}
where, in the second line, we have used the fact that Lemma \ref{lem:rough} implies, recalling Condition (O-\ref{cond:lambda_multiplier2}),
\[
\begin{split}
\left|\ell_{n,k+L_n+1,\lambda}(\cL_{k+L_n,k-L_n,\lambda}h_{ k-L_n,0,\lambda})\right|\geq&\left|\alpha_{k+L_n,k-L_n,\lambda}\ell_{n,k+L_n,\lambda}(h_{k+L_n,0,\lambda})\right|\\
&\times\left| \ell_{k-L_n,k-L_n,\lambda}(h_{ k-L_n,0,\lambda})\right|-\Const(\sigma_n^{-\varpi})\\
\geq& \const -\Const \sigma_n^{-\varpi}\geq \const.
\end{split}
\] 
To complete the estimate needed for Proposition~\ref{prop:diff2}, we need the following lemma.
\begin{lem}
\label{lem:lower rho}
If $C_L$ is chosen large enough, there  exists $\const >0$ such that for all $n \in \mathbb{N}$ and each $\lambda \in \mathbb{R}$, $|\lambda| \le \lambda_1 \sigma_n L_n^{-1}$, we have
\[
\left| \int_{M_n} \cL_{n-1,0,\lambda} \rho \right| \ge \const |\alpha_{n,0,\lambda}| \, .
\]
\end{lem}

We postpone the proof of the lemma and use it to complete the proof of Proposition~\ref{prop:diff2}.

Applying Lemma~\ref{lem:lower rho} to \eqref{eq:final_complex prelim}, we obtain
\begin{equation}
\label{eq:final_complex_1}
\begin{split}
\int_{M_n} \cL_{n-1,k,\lambda}\hat g_k\cL_{k-1,0,\lambda}\rho=&\left(i\lambda\sigma_n^{-1}\!\!\!\sum_{j=k-L_n}^{k+L_n}\Sigma_{k,j}+\cO\left(\lambda^2\sigma_n^{-2}L_n^2+\sigma_n^{-\varpi}\right)\right)\\
&\times \int_{M_n}\cL_{n-1,0,\lambda}\rho \, ,
\end{split}
\end{equation}
which is our main estimate in the case $2L_n \le k \le n-2L_n$.

It remains to consider the cases $k\leq 2L_n$ and $k\geq n-2L_n$. If $k> n-2L_n$, then, using \eqref{eq:complex_2},
\[
\begin{split}
&\int_{M_n} \cL_{n-1,k,\lambda}\hat g_k\cL_{k-1,0,\lambda}\rho d\mu_n=\int_{M_n} \cL_{n-1,k,0}\hat g_k\cL_{k-1,k-L_n,0}\cL_{k-L_n-1,0,\lambda}\rho d\mu_n\\
&+i\sigma_n^{-1}\lambda\sum_{m=k-L_n}^n\int_{M_n}\cL_{n-1,m,0}\,\hat g_m\cL_{m-1,k,0} \hat g_k\cL_{k-1,k-L_n,0}\cL_{k-L_n-1,0,\lambda}\rho d\mu_n\\
&+\cO(\sigma_n^{-2}\lambda^2L_n^2\alpha_{n,0,\lambda}).
\end{split}
\]
Arguing as in Lemma \ref{lem:apriori} yields
\begin{equation}\label{eq:final_complex-2}
\begin{split}
\int_{M_n} \cL_{n-1,k,\lambda}&\hat g_k\cL_{k-1,0,\lambda}\rho  = \sum_{j=k-L_n}^{n}\left(i\lambda\sigma_n^{-1}\Sigma_{k,j}+\cO\left(\lambda^2\sigma_n^{-2}L_n^2\right)\right)\\
&\times \int_{M_n}\cL_{n-1,0,\lambda}\rho+\cO(\sigma_n^{-\varpi}\alpha_{n,0,\lambda}) \\
&=  \sum_{j=k-L_n}^{n}\left(i\lambda\sigma_n^{-1}\Sigma_{k,j}+\cO\left(\lambda^2\sigma_n^{-2}L_n^2 + \sigma_n^{-\varpi} \right)\right)\int_{M_n}\cL_{n-1,0,\lambda}\rho , 
\end{split}
\end{equation}
where we have applied Lemma~\ref{lem:lower rho} in the last line.

Analogously, if $k\leq 2 L_n$, then
\begin{equation}\label{eq:final_complex-3}
\begin{split}
\int_{M_n}\hskip-10pt \cL_{n-1,k,\lambda}\hat g_k\cL_{k-1,0,\lambda}\rho =&\sum_{j=0}^{k+L_n}\left(i\lambda\sigma_n^{-1}\Sigma_{k,j}+\cO\left(\lambda^2\sigma_n^{-2}L_n^2+\sigma_n^{-\varpi}\right)\right)\\
&\times \int_{M_n}\cL_{n-1,0,\lambda}\rho.
\end{split}
\end{equation}

Substituting respectively \eqref{eq:final_complex_1}, \eqref{eq:final_complex-2} or \eqref{eq:final_complex-3}
in each case in equation \eqref{eq:step5.1}, summing over $k$, and recalling \eqref{eq:variance_tot} we finally have
\[
\begin{split}
&\frac{d}{d\lambda}\Upsilon_n(\lambda)=- \left(\lambda+\cE_n(\lambda)\right)\Upsilon_n(\lambda) \\
&\left|\cE_n(\lambda)\right|\leq \Const \left(\lambda^2\sigma_n^{-3}L_n^2 n+ n  \sigma_n^{-\varpi}\right),
\end{split}
\]
concluding the proof of Proposition~\ref{prop:diff2}.
\end{proof}

\begin{proof}[Proof of Lemma~\ref{lem:lower rho}]
We may assume that $n \ge \max \{2L_n, n_1 \}$, where $n_1$ is from Lemma~\ref{lem:coalesce} since otherwise the lemma
is trivial using \eqref{eq:complex_2}.

For such $n$, we write, using (O-\ref{cond:lambda_multiplier2})
\[
\begin{split}
\cL_{n-1,0, \lambda} \rho =& \cL_{n-1,L_n,\lambda} \cL_{L_n-1, 0, \lambda} \rho\\
&= \cL_{n-1, L_n,\lambda} \left( \alpha_{L_n,0,\lambda} \ell_{L_n,0,\lambda}(\rho) h_{L_n,0,\lambda}
 +\mathcal{O} (\alpha_{L_n, 0 ,\lambda} e^{-cL_n} ) \right) \, .
\end{split}
\]
Integrating and estimating as in \eqref{eq:ell example} for each $\varpi>0$ we can choose $C_L$ such that
\[
\alpha_{L_n,0,\lambda} \ell_{L_n,0,\lambda}(\rho)
= \int \cL_{L_n, 0, \lambda} \rho + \mathcal{O} (\alpha_{L_n, 0 ,\lambda} e^{-cL_n} )
= 1 + \mathcal{O} (\alpha_{L_n, 0 ,\lambda} (\lambda \sigma_n^{-1} L_n + \sigma_n^{-\varpi})  ) .
\]
Putting these estimates together, we have 
\[
\begin{split}
\int_{M_n} \cL_{n-1,0,\lambda} \rho & = \int_{M_n} \cL_{n-1, L_n, \lambda} \left( h_{L_n,0,\lambda} + \mathcal{O} (\alpha_{L_n, 0 ,\lambda} (\lambda \sigma_n^{-1} L_n + \sigma_n^{-\varpi})  ) \right) \\
& = \alpha_{n ,L_n+1, \lambda} \left( \ell_{n,L_n,\lambda}(h_{L_n,0,\lambda}) + \mathcal{O} (\alpha_{L_n, 0 ,\lambda} (\lambda \sigma_n^{-1} L_n + \sigma_n^{-\varpi}) + e^{-c |n-L_n|} ) \right)
\end{split}
\]
Next, applying first (O-\ref{cond:lambda_multiplier}) and Lemma~\ref{lem:rough}(a) and then Lemma~\ref{lem:coalesce} yields,
\[
\begin{split}
\left| \int_{M_n} \cL_{n-1,0,\lambda} \rho \right| 
& \ge  K^{-1} \left| \alpha_{n ,L_n+1, \lambda} \alpha_{L_n, 0, \lambda} \left( \ell_{n,L_n,\lambda}(h_{L_n,0,\lambda}) + \mathcal{O} \big( \lambda \sigma_n^{-1} L_n + \sigma_n^{-\varpi} \big) \right) \right| \\
& = K^{-1} | \alpha_{n, 0, \lambda} | \frac{  \ell_{n,L_n,\lambda}(h_{L_n,0,\lambda}) + \mathcal{O} \big( \lambda \sigma_n^{-1} L_n + \sigma_n^{-\varpi} \big)} {  \ell_{n,L_n,\lambda}(h_{L_n,0,\lambda}) + \mathcal{O} \big( e^{-c \min \{ L_n, n-L_n \} } \big) } \\
& \ge \const | \alpha_{n, 0, \lambda} | \, ,
\end{split}
\]
for $\varpi$ large enough and $\lambda$ sufficiently small, where in the last line we have used the lower bound on 
$|\ell_{n,L_n,\lambda}(h_{L_n,0,\lambda})|$ from (O-\ref{cond:lambda_multiplier2}).
\end{proof}


\section{Cones}
\label{sec:complex}

In this section, we describe the machinery of complex cones necessary to verify the contraction
required by (O-\ref{cond:lambda_multiplier2}), and which can also imply (O-\ref{cond:lambda_multiplier}).  For the reader's convenience,
we provide a self contained theory of complex cones more than sufficient for our needs. 
Since we strive for simplicity, our results are not optimal; see \cite{Ru10, Du09} for a more complete, general, but, we believe, less readable exposition.
We will present the theory in a general context, with some tools to aid in the application of
this method to the examples: expanding maps in Section~\ref{sec:expanding} and billiards in Section~\ref{sec:billiard}.
We start by recalling few needed facts about real cones.


\subsection{Real Cones}\ \\
\label{sec:realcones}

Let $V$ be a topological real vector space and $\cS\in V'$ such that $\ell(x)=0$ for all $\ell\in\cS$ implies $x=0$. Define
\begin{equation}\label{eq:com_cone_def0}
C_\bR=\{h\in V\setminus \{0\}\;:\; \ell(h)\geq 0, \forall \ell\in\cS\}.
\end{equation}  
Note that
\begin{equation}\label{eq:no_seminorm}
C_{\bR}\cap-C_{\bR}=\emptyset,
\end{equation}
since if $h\in C_{\bR}\cap-C_{\bR}$, then $\ell(h)=0$ for all $\ell\in\cS$; hence $h=0$, contrary to definition \eqref{eq:com_cone_def0}.

Moreover, assume that there exists $\bbe\in C_\bR$ such that 
\begin{equation}
\label{eq:archimedean}
\mbox{for all $h\in  V$ there exists $\cop\in\bR_+$ such that $\cop\bbe-h\in C_\bR$.}
\end{equation}  
Then one can easily check that
\begin{equation}\label{eq:b_norms}
\| h \| :=\inf\{\cop\in\bR_+\;:\; \ell(\cop\bbe\pm h)\geq 0, \;\forall \ell \in \cS\}
\end{equation}
is a norm.\footnote{Alternatively, it is equivalent to define this norm by $\|h\|=\inf\{\cop\in\bR\;:\; -\cop \bbe\preceq h\preceq \cop \bbe\}$, where $h\preceq g$ iff $g-h\in C_\bR $.} 
Note that $\| \bbe \| = 1$.
 Let $\cB_\bR$ be the completion of $V$ in such a norm.\footnote{ Note that the expression \eqref{eq:b_norms} defines the norm also on $\cB_\bR$.} Since, for all $\ell\in \cS$, $\ell(\|h\|\bbe\pm h)\geq 0$, it follows that $\cS\subset \cB_\bR'$.

Finally, let $\cS_*$ be the weak-$*$ closure of the convex hull of $\{\cop\ell\;:\;\cop\in\bR_+,\ell\in\cS\}$.
To simplify matters further, we restrict to the special case in which there exists $\bbm\in\cS_*$ and $\kappa\in (0,1)$ such that 
\begin{equation}\label{eq:m-normalization}
\bbm(\bbe)=1
\end{equation}
and, for all $h\in\cC_\bR$,
\begin{equation}\label{eq:maybetoomuch}
\bbm(h)\geq \kappa \|h\|.
\end{equation}
Having described the abstract setting, we can now proceed to explore its properties. 
We start defining a larger cone $\cC_\bR \supset C_\bR$ by
\begin{equation}\label{eq:com_cone_def}
\cC_\bR=\{h\in\cB_\bR\setminus \{0\}\;:\; \ell(h)\geq 0, \forall \ell\in\cS\}.
\end{equation}
Next, we define the dual cone as
\begin{equation}\label{eq:dual_cone_def}
\cC'_\bR=\{\ell\in \cB'_\bR\;:\; \ell(h)\geq 0\quad \forall h\in\cC_\bR\}\supset \cS
\end{equation}
Note that we have
\begin{equation}\label{eq:com_cone_def2}
\cC_\bR=\{h\in\cB_\bR\setminus \{0\}\;:\; \ell(h)\geq 0, \forall \ell\in\cC'_\bR\}.
\end{equation}
Hence, again, $\cC_\bR\cap -\cC_\bR=\emptyset$.\footnote{ Indeed, if $h\in \cC_\bR\cap -\cC_\bR$, then $\ell(h)=0$ for all $\ell\in \cS$, and $\|h\|=0$, by \eqref{eq:b_norms}.} 
The cone $\cC_{\bR}$ is associated with the Hilbert metric (see \cite[Equation (D.1.2)]{DKL}): for all $g,h\in \cC_{\bR}$,
\[
\aligned
\alpha (h,\,g)=&\sup\{\cop\in\mathbb{R}^+\;|\; g-\cop h\in\cC_{\bR}\}\\
\beta(h,\,g)=&\inf\{\cop\in\mathbb{R}^+\;|\; \cop h-g\in\cC_{\bR}\}\\
d_{H,\cC_\bR}(h,\,g)=&\ln \left[\frac{\beta(h,\,g)}{\alpha(h,\,g)}\right].
\endaligned
\]
From now on we will write simply $d_H$, if the cone is clear from the context.\\
Since by \eqref{eq:com_cone_def} $g - \cop h \in \cC_{\bR}$ iff $\ell(g-\cop h) \ge 0$ for all $\ell \in \cS$, taking the
limit $\cop \uparrow \alpha$ yields $\alpha = \inf_{\ell \in \cS} \frac{\ell(g)}{\ell(h)}$.
Alternatively, using \eqref{eq:com_cone_def2} in the same way, we have
$\alpha = \inf_{\ell \in \cC'_{\bR}} \frac{\ell(g)}{\ell(h)} = \inf_{\ell \in \cS} \frac{\ell(g)}{\ell(h)}$.
Arguing similarly for $\beta$ yields
\begin{equation}\label{eq:hilbert_real}
d_H(h,g)=\sup_{\ell, p\in\cS}\ln\frac{\ell(h)p(g)}{\ell(g)p(h)}=\sup_{\ell, p\in\cC_{\bR}'}\ln\frac{\ell(h)p(g)}{\ell(g)p(h)}.
\end{equation}

Let us recall the relevance of real cones for our context. 
\begin{thm}[{\cite[Theorem 1.1]{Li95}}]\label{thm:real_contraction}
Let $\cB_1,\cB_2$ be real Banach spaces with real cones $\cC_1,\cC_2$ as above.  Let $L\in L(\cB_1,\cB_2)$, such that $L\cC_1\subset \cC_2$ and 
\[
\Delta=\sup_{x,y\in\cC_1}d_{H,\cC_2}(Lx,Ly). 
\]
Then, for all $h,g\in\cC_1$ we have
\[
d_{H,\cC_2}(Lh,Lg)\leq \tanh(\Delta/4)d_{H,\cC_1}(h,g).
\]
\end{thm}
\begin{lem}[{\cite[Lemma 1.3]{Li95}}] 
Given $h,g\in\cC$ for which $\| h\|=\|g\|$,
\[
\|h-g\|\leq \left(e^{d_H(h,g)}-1\right)\| h \|.
\]
\end{lem}
In the next section, we describe how the above results can be extended to the case of operators acting on complex Banach spaces, see, in particular, Theorem \ref{thm:complex_contraction} and Lemma \ref{lem:compare}.


\subsection{Complex Cones}\ \\
\label{sec:complexcones}

Let $\cB_\bC$ be the complexification of $\cB_\bR$ as detailed in Lemma \ref{lem:complexification}.

\begin{defin}[Complex cone and its dual]
\label{def:complex cone}
The associated complex cone is defined as $\cC_\bC:=\bC_* \cdot (\cC_\bR+i\cC_{\bR})$, where $\bC_*=\bC\setminus\{0\}$.
We may also write $\cC_\bC = \bC_1 \cdot (\cC_\bR + i \cC_\bR)$, where $\bC_1 = \{ z \in \bC : |z| =1 \}$.
We define the dual cone by 
\[
\cC'_\bC=\{\ell\in \cB'_\bC\;:\; \ell(h)\neq 0\quad \forall h\in\cC_\bC\} .
\]
\end{defin}
For any element $\ell\in\cB'_\bR$, its action on $h+ig\in\cB_\bC$, $h,g \in \cB_\bR$, is naturally defined by $\ell(h+ig)=\ell(h)+i\ell(g)$.

Then, for all $h\in\cC_{\bC}$, letting $z\in\bC$  such that $|z|=1$ and  $zh=x+iy$ with $x,y\in \cC_{\bR}$,
\begin{equation}\label{eq:m_lower}
\begin{split}
| \bbm(h) |&=|\bbm(x+iy)|=\sqrt{\bbm(x)^2+\bbm(y)^2}\geq \kappa\sqrt{\|x\|^2+\|y\|^2}\\
&\geq \frac{\kappa}{\sqrt 2}\|x+iy\|= \frac{\kappa}{\sqrt 2}{ \| h \| } ,
\end{split}
\end{equation}
where in the second line we have used \eqref{eq:norms_equiv}.
In addition, if $\|u\|\leq 1$, then
\begin{equation}
\label{eq:real int}
\bbe+\frac 12u\succeq \bbe-\frac 12\|u\|\bbe\succeq \frac 12\bbe\succeq 0
\end{equation}
that is the interior of $\cC_\bR$, and hence of $\cC_\bC$, is not empty. That is $\cC_\bC$ is a {\em regular cone} in the language of \cite[ Definition 3.2-(3)]{Ru10}. Moreover, for  each $\ell\in\cB_\bR'$ such that $\|\ell\|<\kappa$ we have $\bbm+\ell\in\cC_\bR'$, so also $\cC_{\bR}'$ has non empty interior.

The first statement of the following Lemma is contained in \cite[Proposition 5.2]{Ru10}, while the second is the first part of the proof of \cite[Lemma 4.1]{Du09}.\footnote{ Note however that in \cite[Lemma 4.1]{Du09} it is assumed that
$x,y$ are linearly independent, which is not necessary.}
\begin{lem}\label{lem:dual_comp}
We have the following characterizations.
\begin{itemize}
  \item[a)] $\displaystyle \cC_\bC=\{h\in\cB_\bC\setminus\{0\}\;:\; \forall \ell, m\in \cC'_\bR,\; \Re\left(\ell(h)\overline{m(h)}\right)\geq 0\}$.
  \item[b)] $\displaystyle \cC_\bC'=\{ \ell\in \cB_{\bC}'\;:\; \forall x, y\in \cC_\bR,\; \Re(\ell(x)\overline{\ell(y)})>0 \}$. \smallskip
  \item[c)] $\displaystyle \cC'_\bC\supset\{\pm\ell\pm i p\;:\; \ell, p\in \mathring{\cC}'_\bR\}=:\hat\cC'_{\bC}$,
where $\mathring{\cC}'_\bR=\{\ell\in\cC'_{\bR}\;:\; \ell(x)>0 \, \forall x\in\cC_{\bR}\}$.
\end{itemize}
\end{lem}
\begin{proof}
(a) To prove that $\cC_\bC\subseteq\{h\in\cB_\bC\setminus\{0\}\;:\; \forall \ell, m\in \cC'_\bR,\; \Re\left(\ell(h)\overline{m(h)}\right)\geq 0\}$,  check that $\Re\left(\ell(z h)\overline{m(zh)}\right)=|z|^2\Re\left(\ell(h)\overline{m(h)}\right)$, for all $z\in\bC_*$. Also, by definition, for each $h\in\cC_\bC$ there exists $z\in\bC$ such that $zh=x+iy$ with $x,y\in\cC_\bR$. Finally, for $x,y\in\cC_\bR$,  $\Re\left(\ell(x+iy)\overline{m(x+iy)}\right)=\ell(x)m(x)+\ell(y)m(y)\geq 0$.

To prove the opposite inclusion, let $h\in \{h\in\cB_\bC\;:\; \forall \ell, m\in \cC'_\bR,\; \Re\left(\ell(h)\overline{m(h)}\right)\geq 0\}$. We claim that there exists $z\in\bC$ such that $\tilde h:=zh=x+iy$ with $x\in\cB_\bR$ and $y\in\cC_\bR$. Indeed, let $h=u+iv$, $u,v\in\cB_\bR$ and assume that $u,v\not\in\cC_\bR\cup-\cC_\bR$. If $u= \chi_* v$, $\chi_*\in\bR$, then we can choose $\ell, m\in\cC_\bR'$ such that $\ell(v)>0$ and $m(v)<0$, so  $\Re\left(\ell(h)\overline{m(h)}\right)= (1+\chi_*^2)\ell(v)m(v)<0$ contrary to assumptions. Hence $u,v$ must be linearly independent in $\cB_\bR$. Note that there must exist $\ell \in\cC_\bR'$ such that 
$\ell(u)=0$ and $\ell(v)\neq 0$.\footnote{ Indeed, by definition $\cS\subset \cC_\bR'$ and if $\pm u\not\in\cC_\bR$, then there must exist $\ell_1, p\in\cS$ such that $\ell_1(u), p(-u)<0$. Hence, $\ell=\ell_1-\frac{\ell_1(u)}{p(u)}p\in\cC_{\bR}'$ and $\ell(u)=0$. Next, if for all $\ell\in\cC_\bR'$ such that $\ell(u)=0$ we have also $\ell(v)=0$, then the same must happen for all $\ell\in\cB_\bR'$ (since $\cC_\bR'$ has a non empty interior) and this is impossible by the Hahn-Banch theorem.} Thus $\Re\left(\ell(h)\overline{m(h)}\right)=\ell(v)m(v)$, and choosing $m\in\cC_\bR'$ such that $\ell(v)m(v)<0$ leads to a contradiction. Hence, if $v\in\cC_\bR$ we are done, if $v\in -\cC_\bR$ then we choose $z=-1$, if $u\in\cC_\bR$ we choose $z=i$ and finally if $u\in-\cC_\bR$ we choose $z=-i$.

Then for each $\ell\in \cC'_\bR$ such that $\ell(\tilde h)\neq 0$, we can write  $\ell(\tilde h)|\ell(\tilde h)|^{-1}=e^{i\theta_\ell}$, $\theta_\ell\in [0,\pi]$. Suppose that there exists $\ell_1,\ell_2\in \cC'_\bR$ such that $\pi\geq |\theta_{\ell_2}-\theta_{\ell_1}|> \pi/2$, then 
\[
\Re(\ell_2 (h)\overline{\ell_1(h)})= |z|^{-2} |\ell_2(\tilde h)\ell_1(\tilde h)|\cos(\theta_{\ell_2}-\theta_{\ell_1})<0
\]
contrary to the assumption. Thus there exists $\vf\in [0,\pi/2]$ such that for each $\ell\in \cC'_\bR$, $\ell(\tilde h)\neq 0$, $\theta_\ell\in [\vf-\pi/4,\vf+\pi/4]$. Let $e^{-i\vf}\tilde h=u+iv$, then $\ell(u+iv)=|\ell(u+iv)|e^{i(\theta_\ell-\vf)}$ which implies $\ell(u)\geq |\ell(v)|$, that is $\ell(u\pm v)\geq 0$. This means that $u\pm v\in\cC_\bR$. Accordingly, setting 
$x+iy=(1+i)e^{-i\vf} \tilde h = (1+i)(u+iv)=u-v+i(u+v)$ we have that $x,y\in\cC_{\bR}$, hence $h\in\cC_\bC$.

Next, we prove statement (b).
Let $\ell\in \{ \ell\in \cB_{\bC}'\;:\; \forall x, y\in \cC_\bR,\; \Re(\ell(x)\overline{\ell(y)})>0 \}$. 
If there exist $x,y\in\cC_\bR$ such that $\ell(x+iy)=0$, then $\ell(x)=-i\ell(y)$ and $\Re(\ell(x)\overline{\ell(y)})=\Re(-i|\ell(y)|^2)=0$ contrary to the assumption.
Hence $\{ \ell\in \cB_{\bC}'\;:\; \forall x, y\in \cC_\bR,\; \Re(\ell(x)\overline{\ell(y)})>0 \}\subset \cC_\bC'$.

To prove the opposite inclusion, let $\ell \in \cC_\bC'$.  Suppose that there exist $x,y\in\cC_\bR$ such that $\Re(\ell(x)\overline{\ell(y)})\leq 0$ and write $\ell(x)=re^{i\theta}\neq 0$, $\ell(y)=s e^{i\vf}\neq 0$, with $\theta,\vf\in [-\pi,\pi)$. Then $0\geq \Re(\ell(x)\overline{\ell(y)})=rs \cos(\theta-\vf)$, so $\frac 32\pi\geq |\theta-\vf|\geq \pi/2$.
Define 
\[
z=-se^{-i(\theta-\vf)}x+r y=-s\cos(\theta-\vf) x+r y+is \sin (\theta-\vf)x=: u+iv 
\]
then $\ell(z)=0$. Note that $z\neq 0$ since $-s\cos(\theta-\vf)\geq 0$ hence $0\neq u\in\cC_\bR$. If $\sin(\theta-\vf)\geq 0$, then also $v\in\cC_\bR$.
Otherwise $i z=u_1+iv_1$ with $u_1,v_1\in\cC_\bR$. Thus there exists $w\in\cC_{\bC}$ such that $\ell(w)=0$, contradicting the hypothesis.

Finally, we prove statement (c).
Let $\pm \ell,\pm p\in\mathring{\cC}'_\bR$, then for each $x,y\in\cC_\bR$ we have
\[
\Re((\ell+ip)(x)\cdot (\ell-ip)(y))=\ell(x)\ell(y)+p(x)p(y)> 0,
\]
which implies $\ell + i p \in \cC'_\bC$ by part (b) of the lemma.
\end{proof}

\begin{lem}\label{lem:upper_dual} For each $\ell\in\cC_{\bC}'$ we have
\[
\|\ell \|' \leq \sqrt 2 |\ell(\bbe)|.
\]
\end{lem}
\begin{proof}
If $\ell\in\cC_{\bC}'$, then, when acting on $\cB_{\bR}$ it can be written as $\ell(x)=\ell_1(x)+i\ell_2(x)$.
By Lemma \ref{lem:dual_comp}(b) we have, for each $x\in\cB_{\bR}$,
\[
0\leq \Re(\ell(\|x\|\bbe-x)\overline{\ell(\|x\|\bbe+x)})=\ell_1(\|x\|\bbe-x)\ell_1(\|x\|\bbe+x)+\ell_2(\|x\|\bbe-x)\ell_2(\|x\|\bbe+x)
\]
which implies
\[
\begin{split}
|\ell(x)|&\leq \|x\| |\ell(\bbe)|.
\end{split}
\]
If $x+iy\in\cB_{\bC}$, then
\[
|\ell(x+iy)|\leq |\ell(\bbe)|(\|x\|+\|y\|)\leq \sqrt 2 |\ell(\bbe)|\|x+iy\|.
\]
\end{proof}

\begin{defin}\label{def:convex}
We call the complex cone {\em linearly convex} if for each $g\in\cB_\bC\setminus \cC_\bC$
we can find $\ell\in \hat\cC'_\bC$ such that $\ell(g)=0$.\footnote{ Note that this definition is a bit different from the one in \cite[Section 3]{Du11}. But it has the advantage of simplifying the proof of Lemma \ref{lem:lincox}.}
\end{defin}

\begin{rem} \label{rem:cstar=c}
Note that if $\cC_\bC$ is linearly convex, then 
\[
\cC_*:=\{h\in\cB_\bC\;:\; \ell(h)\neq 0\;\forall\; \ell \in \hat\cC'_\bC\}=\cC_{\bC}.
\]
Indeed, if $h\in\cC_*$, then $h\in \cC_{\bC}$, otherwise, by linear convexity, there exists $\ell\in\hat\cC'_\bC$ such that $\ell(h)=0$, contrary to the assumption. On the other hand, if $h\in\cC_\bC$, then for each $\ell \in \hat\cC'_\bC$, by Lemma \ref{lem:dual_comp}(c) and Definition \ref{def:complex cone}, we have that $\ell(h)\neq 0$; hence, $h\in\cC_*$.
\end{rem}

The next Lemma corresponds to \cite[Lemma 4.1]{Du11}.

\begin{lem}\label{lem:lincox} The cone $\cC_\bC$ is linearly convex.
\end{lem}
\begin{proof}
Let $x+iy\not \in\cC_\bC$, $x,y\in\cB_\bR$, then by the first statement of Lemma \ref{lem:dual_comp}, there exists $\tilde \ell, \tilde p\in\cC_\bR'$ such that  
\[
\tilde \ell(x)\tilde p(x)+\tilde \ell(y)\tilde p(y)=\Re(\tilde \ell(x+iy)\tilde p(x-iy))<0.
\]
Note that, setting $p=\tilde p+\ve \bbm$ and $\ell=\tilde \ell+\ve\bbm$, with $\ve>0$, then $\ell, p\in\mathring{\cC}_\bR'$ and for $\ve$ sufficiently small, then we still have
\begin{equation}\label{eq:basic_lincovex}
\ell(x) p(x)+\ell(y) p(y)<0.
\end{equation}
If $y=0$ and $p(x)>0$, then $q=(1+i)[p(x)\ell-\ell(x) p ]\in \hat \cC_\bC'$ and $q(x)=0$, while if $p(x)<0$ the same holds for $q=(1+i)[-p(x)\ell+\ell(x) p]$. If $y\neq 0$, then we can perturb so that $\ell(y)\neq 0$ and still $\ell, p\in\mathring{\cC}_\bR'$ and \eqref{eq:basic_lincovex} holds.
Define
\[
\begin{split}
&q=[\alpha_p \ell(x)+p(x)]\ell-i\ell(y)[\alpha_p\ell+p]\\
&\alpha_p=-\frac{p(x)\ell(x)+\ell(y)p(y)}{\ell(x)^2+\ell(y)^2}>0.
\end{split}
\]
One can check that $q(x+iy)=0$.
Note that $\alpha_p\ell+p\in\mathring{\cC}_\bR'$, thus $q\in\hat \cC_\bC'$ if $\alpha_p \ell(x)+p(x)\neq 0$. If instead $\alpha_p \ell(x)+p(x)= 0$, then $(1+i)q\in\hat \cC_\bC'$, concluding the proof of the lemma.
\end{proof}

In order to construct a metric on the cone, for each $h,g\in\cC_{\bC}$ define
\begin{equation}\label{eq:gauge}
E_\cC(h,g)=\left\{\frac{\ell(h)}{\ell(g)}\;:\; \ell\in\cC'_{\bC}\right\} \, .
\end{equation}

For future use, it is convenient to record the following facts.
\begin{lem}\label{lem:are_the_same}
For each $h,g\in\cC_{\bC}$ we have
\[
E_\cC(h,g)=\left\{\frac{\ell(h)}{\ell(g)}\;:\; \ell\in\hat \cC'_{\bC}\right\}=:\widehat E_\cC(h,g).
\]
\end{lem}
\begin{proof}
Clearly $\widehat E_\cC(h,g)\subset  E_\cC(h,g)$.
On the other hand, note that if $\zeta\in E_\cC(h,g)$, then there exists $\ell\in\cC_{\bC}'$ such that $\ell(h-\zeta g)=0$. But this implies $h-\zeta g\not\in\cC_{\bC}$, otherwise, by definition, there would exist $z$ such that $z(h-\zeta g)=x+iy$ with $x,y\in\cC_{\bR}$ and this would imply $\Re(\ell(x)\overline{\ell(y)})=0$ contrary to Lemma \ref{lem:dual_comp}(b). Accordingly, by Lemma~\ref{lem:lincox} there exists
 $\ell\in\hat\cC_{\bC}'$ such that $\ell(h-\zeta g)=0$, hence $\zeta\in \widehat E_\cC(h,g)$. That is $E_\cC(h,g)= \widehat E_\cC(h,g)$ and the lemma follows.
\end{proof}

\begin{lem}\label{lem:E-to-cone}
For each $h,g\in\cC_{\bC}$, $z\not \in E_{\cC}(h,g)$ iff $zg-h\in\cC_\bC$.
\end{lem}
\begin{proof}
If $z\not \in E_{\cC}(h,g)$, then, by Lemma \ref{lem:are_the_same}, we have that for all $\ell\in \hat \cC'_{\bC}$ $z\neq \frac {\ell(h)}{\ell(g)}$, that is $\ell(h-z g)\neq 0$. Hence, by Remark \ref{rem:cstar=c}, $zg-h\in\cC_\bC$. On the other hand, if $zg-h\in\cC_\bC$, then again by  by Remark \ref{rem:cstar=c}, we have that, for each $\ell\in \hat \cC'_{\bC}$, $\ell(zg-h)\neq 0$. Hence, $z\neq \frac {\ell(h)}{\ell(g)}$, that is $z\not \in E_{\cC}(h,g)$.
\end{proof}

Next, we define the key object in the theory,
\begin{equation}\label{eq:gauge_dist}
\delta_\cC(h,g)=\ln\frac{\sup_{z\in E_\cC(h,g)}|z|}{\inf_{z\in E_\cC(h,g)}|z|}=\sup_{z,w\in E_\cC(h,g)}\ln\left|\frac z w\right| \, .
\end{equation}
As a first result, we have a version of \cite[Lemma 2.6]{Du09} showing the relevance of $\delta_\cC$.
\begin{lem}\label{lem:compare}
Let $h,g\in\cC_{\bC}$, such that $\bbm(h)=\bbm(g)$, $|\bbm(g)|=1$.
Then
\[
\|h-g\|\leq \frac{\sqrt2}{\kappa}\delta_{\cC}(h,g).
\]
\end{lem}
\begin{proof}
There exists $z\in\bC$, $|z|=1$ such that $\bbm(zh)=\bbm(zg)=1$. Let $\tilde h=z h$ and $\tilde g=z g$. Also, by definition $\delta_{\cC}(h,g)=\delta_{\cC}(\tilde h,\tilde g)$. Then, since $\bbm\in\cS_*$,
\[
a:=\inf_{z\in E_\cC(\tilde h,\tilde g)}|z|\leq \left|\frac{\bbm(\tilde h)}{\bbm(\tilde g)}\right|=1\leq \sup_{z\in E_\cC(\tilde h,\tilde g)}|z|=:b
\]
Note that if $a=0$ or $b=\infty$, then the Lemma is trivially satisfied. For each $0<\alpha<a\leq 1\leq b<\beta$  we have $\alpha\tilde g-\tilde h, \beta\tilde g-\tilde h\in\cC_\bC$. Indeed, if $\alpha\tilde g-\tilde h\not\in\cC_\bC$, then, since the cone is linearly convex (see Definition \ref{def:convex}), there would exist $\ell\in\hat\cC_{\bC}'$ such that $\ell(\alpha\tilde g-\tilde h)=0$, but then $\alpha=\frac{\ell(\tilde h)}{\ell(\tilde g)}\geq a$ contrary to the assumption; the argument for $\beta\tilde g-\tilde h\in\cC_\bC$ is the same.

We can then write $(\beta-\alpha)(\tilde g-\tilde h)=(\beta-1)(\alpha\tilde g-\tilde h)+(1-\alpha)(\beta\tilde g-\tilde h)$, and, recalling \eqref{eq:m_lower}, we have 
\[
\begin{split}
\|\tilde g-\tilde h\|&\leq \frac {\beta-1}{\beta-\alpha}\|\alpha\tilde g-\tilde h\|+\frac {1-\alpha}{\beta-\alpha}\|\beta\tilde g-\tilde h\|\\
&\leq \frac{\sqrt2}{\kappa}\left[ \frac {\beta-1}{\beta-\alpha}|\bbm(\alpha\tilde g-\tilde h)|+\frac {1-\alpha}{\beta-\alpha}|\bbm(\beta\tilde g-\tilde h)|\right]\\
&\leq \frac{2\sqrt2}{\kappa}\frac {(\beta-1)(1-\alpha)}{\beta-\alpha}=\frac{2\sqrt2}{\kappa}\frac {(\sqrt \beta-\sqrt \alpha)^2-(1-\sqrt{\alpha \beta})^2}{\beta-\alpha}\\
&\leq \frac{2\sqrt2}{\kappa}\frac {\sqrt \beta-\sqrt \alpha}{\sqrt \beta+\sqrt \alpha}\leq  \frac{2\sqrt2}{\kappa}\ln \sqrt{\frac \beta \alpha}
= \frac{\sqrt2}{\kappa}\ln \frac \beta \alpha.
\end{split}
\]
where, in the last line, we have used that, for $x\geq 1$, $\ln x\geq \frac{x-1}{x+1}$. Taking the inf on $\beta$ and the sup on $\alpha$ yields
\[
\|g-h\|=\|\tilde g-\tilde h\|\leq \frac{\sqrt2}{\kappa}\delta_{\cC}(\tilde h,\tilde g)=\frac{\sqrt2}{\kappa}\delta_{\cC}(h,g).
\]
\end{proof}
For $z\in\bC$ and all $h\in\cC_\bC$ we have $\delta_\cC(h,z h)=0$. It is then natural to define the equivalence relation $h\sim g$ iff $h=zg$ for some $z\in\bC$. Let $\widetilde \cC_{\bC}$ be the quotient of $\{h\in\cC_\bC\;:\; \delta_\cC(h,\bbe)<\infty\}$.
The following is similar to \cite[Theorem 3.1]{Du11}.\footnote{This result is not used in the present paper; we include it for completeness.}
\begin{thm}\label{thm:complete}
The space $(\widetilde \cC_{\bC}, \delta_\cC)$ is a complete metric space.
\end{thm}
\begin{proof}
Clearly, $\delta_\cC(g,h)=\delta_\cC(h,g)$. For each $h,g,f\in \widetilde \cC_{\bC}$ and $\ell\in\hat \cC'_{\bC}$ we have
\[
\inf_{z\in E_\cC(h,g)}|z|\inf_{z\in E_\cC(g,f)}|z|\leq\left|\frac{\ell(h)}{\ell(g)}\right|\,\left|\frac{\ell(g)}{\ell(f)}\right|=\left|\frac{\ell(h)}{\ell(f)}\right|\leq \sup_{z\in E_\cC(h,g)}|z|\sup_{z\in E_\cC(g,f)}|z|,
\]
hence, by \eqref{eq:gauge_dist} and taking the sup and inf on $\ell$, the triangle inequality follows. 

Next, if $\delta_\cC(h,g)=0$, then by Lemma \ref{lem:compare} we have $\frac{h}{\bbm(h)}=\frac{g}{\bbm(g)}$,  that is $h\sim g$.

As for the completeness, it follows from Lemma \ref{lem:compare} that if $\{h_n\}\subset\cC_{\bC}$ is a Cauchy sequence with respect to $\delta_\cC$, then
\[
\left\|\frac{1}{\bbm(h_n)}h_n- \frac{1}{\bbm(h_m)}h_m\right\|\leq \frac{\sqrt2}{\kappa}\delta_{\cC}(h_n,h_m).
\]
So $\frac{1}{\bbm(h_n)}h_n$ is a Cauchy sequence in $\cB_{\bC}$, hence it has a limit $h_*$, and by Lemma \ref{lem:dual_comp}(a), $h_*\in\cC_{\bC}$. 
Note that, for each $g\in\cC_\bC$ and $q\in\hat \cC_\bC'$ with $|q(\bbe)|=1$, recalling equation \eqref{eq:m_lower},\footnote{Note that, since $\bbm\in\cS_*$, $\left|\frac{\bbm(h)}{\bbm(g)}\right|\in \overline{E_{\cC}(h,g)}$.}
\[
e^{\delta_{\cC}(\bbe,g)}\geq \left|\frac{\bbm(g) q(\bbe)}{\bbm(\bbe)q(g)}\right|\geq \const \frac{\|g\| }{|q(g)|} \, .
\]
Hence, $|q(h_*)|\geq \const e^{-\delta_{\cC}(\bbe,h_*)}\|h_*\|\geq \const$, since $\delta_{\cC}(\bbe,h_*)\leq \delta_{\cC}(\bbe,h_m)+1$, for $m$ large enough.
For all $\ell, p\in\hat\cC_{\bC}'$ let $\tilde \ell=\ell(\bbe)^{-1}\ell$ and $\tilde p=p(\bbe)^{-1}p$, then, 
setting $\tilde h_n = h_n/\bbm(h_n)$ and recalling Lemma \ref{lem:upper_dual},
\[
\left|\frac{\ell(h_*)p(h_n)}{\ell(h_n)p(h_*)}\right|=\left|\frac{\tilde \ell(h_*)\tilde p(\tilde h_n )}{\tilde \ell(\tilde h_n)\tilde p(h_*)}\right|\leq \frac{1+\frac{|\tilde p(h_*- \tilde h_n)|}{\tilde p(h_*)}}{1-\frac{|\tilde \ell(h_*- \tilde h_n)|}{\tilde \ell(h_*)}}\leq \frac{1+\Const \|h_*- \tilde h_n \|}{1-\Const \|h_*- \tilde h_n \|}
\]
which, taking the sup on $\ell,p$ implies $\lim_{n\to\infty}\delta_\cC(h_n,h_*)=0$.
\end{proof}

The previous result is nontrivial provided $\{ h \in \cC_{\bC} : \delta_{\cC}(h , \bbe) < \infty \} \neq \emptyset$.  This
follows from the fact that $\bbe$ is in the interior of $\cC_{\bC}$, which is verified by the following lemma. 

\begin{lem}
\label{lem:interior}
For any $h \in \cB_{\bC}$ and for all $z \in \mathbb{\bC}$ with $|z| < (\sqrt{2} \|h\|)^{-1}$ we have
$\bbe + z h \in \cC_{\bC}$ and $\delta_{\cC}(\bbe, \bbe + z h) \le \ln \left( \frac{1 + |z| \sqrt{2} \| h \| }{ 1 - |z| \sqrt{2} \| h \| } \right) < \infty$.
\end{lem}

\begin{proof}
We begin by claiming that if $h \in \cB_{\bC}$, then $h + z \bbe \in \cC_{\bC}$ for all $z \in \bC$ with $|z| > \sqrt{2} \| h\|$.  
Indeed, by Remark~\ref{rem:cstar=c} and Lemma~\ref{lem:lincox}, if $h + z \bbe \notin \cC_{\bC}$, then
there exists $\ell \in \hat{\cC}'_{\bC}$ such that $\ell(h + z \bbe) = 0$.  But then applying Lemma~\ref{lem:upper_dual}, one has
\[
| z \ell(\bbe) | = | \ell(h)| \le  \sqrt{2} \| h \| | \ell(\bbe) | ,
\]
which is a contradiction.

With the claim proved, it follows that $\bbe + z h \in \cC_{\bC}$ for all $|z| < (\sqrt{2} \|h \|)^{-1}$.  Moreover, for any
$\ell \in \hat{\cC}'_{\bC}$,
\[
\frac{\ell(\bbe + z h)}{\ell(\bbe)} = 1 + z \frac{\ell(h)}{\ell(\bbe)} \, ,
\]
so that according to the definition \eqref{eq:gauge_dist} and using again Lemma~\ref{lem:upper_dual}, 
\[
\delta_{\cC}(\bbe + z h, \bbe) \le \ln \left( \frac{1 + |z| \sqrt{2} \| h \| }{ 1 - |z| \sqrt{2} \| h \| } \right) \, ,
\]
as required.
\end{proof}

The metric $\delta_{\cC}$ also provides a partial converse to Lemma~\ref{lem:upper_dual}, as follows.
\begin{lem}
\label{lem:partial converse}
For any $h \in \cC_{\bC}$ and all $\ell \in \cC_{\bC}'$,
\[
|\ell(h)| \ge e^{- \delta_{\cC}(h, \bbe)} \tfrac{\kappa}{\sqrt{2}} \| h \| |\ell(\bbe)| \, .
\]
\end{lem}

\begin{proof}
Let $h \in \cC_{\bC}$ and suppose $D = \delta_{\cC}(h , \bbe)$.  Then by definition \eqref{eq:gauge_dist},
for any $\ell \in \cC_{\bC}'$, 
\[
e^D \ge \frac{|\ell(\bbe) \bbm(h)|}{|\ell(h) \bbm(\bbe)|} \ge \frac{|\ell(\bbe)|}{|\ell(h)|} \frac{\kappa}{\sqrt{2}} \| h \| \, ,
\]
where we have used \eqref{eq:m_lower}.  The lemma follows.
\end{proof}

The next is a basic result in \cite{Du09}. We provide a detailed proof, although it is essentially as in \cite{Du09}, because our definitions differ slightly from the ones in \cite{Du09}.
\begin{thm}[{\cite[Theorem 3.1 (iii)]{Du11}}]\label{thm:complex_contraction}
Let $\cB_1,\cB_2$ be complex Banach spaces with complex cones $\cC_1,\cC_2$, satisfying hypotheses \eqref{eq:com_cone_def0}, \eqref{eq:archimedean} and \eqref{eq:maybetoomuch}. Let $L\in L(\cB_1,\cB_2)$, such that $L\cC_1\subset \cC_2$ and 
$\Delta=\sup_{x,y\in\cC_1}\delta_{\cC_2}(Lx,Ly)$. Then, for all $h,g\in\cC_1$ we have
\[
\delta_{\cC_2}(Lh,Lg)\leq \tanh(\Delta/4)\delta_{\cC_1}(h,g).
\]
\end{thm}

\begin{proof}
Let $g,h\in \cC_1$. If $g=z h$, $z\in\bC$, then we have seen that $\delta_{\cC_2}(Lh,Lg)=\delta_{\cC_1}(h,g)=0$, hence the theorem is trivially true. We can then assume $g,h$ are 
linearly independent. Let $\zp,\zm\in \bC$ be such that 
\begin{equation}\label{eq:Mm_def}
\begin{split}
&|\zp|=M := \sup\{|z|\;:\; z\in E_{\cC_1}(h,g)\}\\
&|\zm|=m := \inf\{|z|\;:\; z\in E_{\cC_1}(h,g)\}.
\end{split}
\end{equation}
 Since 
if $M =\infty$ or $m= 0$, the statement is trivial, we can assume $\infty >M > m>0$. 
Indeed it must be that $M>m$, otherwise $M=m$ and $E_{\cC_1}(h,g)\subset \{z\in \bC\;:\; |z|=m\}$. But Lemma \ref{lem:dual_comp}(a) implies that $\cC_\bC\cup\{0\}$ is closed, hence its complement is open, and so is $E_{\cC_1}(h,g)$ by Lemma \ref{lem:E-to-cone},
contradicting the hypothesis that it is contained in a circle.
The openness of $E_{\cC_1}(h,g)$ implies that
 $\zm,\zp \not\in E_{\cC_1}(h,g)$. Hence, by Lemma \ref{lem:E-to-cone},
\[
\zp g-h\in \cC_1\, ,\quad h-\zm g\in \cC_1.
\]
Accordingly,  by assumption,
\[
L(\zp g-h)\in \cC_2\,,\quad
L(h-\zm g)\in \cC_2
\]
and,
\begin{equation}\label{eq:diam_hyp}
\delta_{\cC_2}(L(\zp g-h),L(h-\zm g))\leq \Delta.
\end{equation}
Let $z\not\in E_{\cC_2}(L(h-\zm g), L(\zp g-h))$. By Lemma \ref{lem:E-to-cone} this is equivalent to
\[
z L(\zp g-h)-L(h-\zm g)\in \cC_2
\]
or
\[
\frac{z \zp+\zm}{z +1}Lg-Lh\in \cC_2.
\]
Accordingly, if we define the Möbius transformation $\psi(z)=\frac{z \zp+\zm}{z +1}$, then 
\[
\psi(E_{\cC_2}(L(h -\zm g), L(\zp g-h)))=E_{\cC_2}(Lh,Lg).
\]
Recalling \eqref{eq:diam_hyp} and \eqref{eq:gauge_dist}
\begin{equation}\label{eq:direction_un}
e^\Delta\geq \sup_{z,w\in\psi^{-1}(E_{\cC_2}(Lh,Lg))}\left|\frac z w\right|=\sup_{\alpha,\beta\in E_{\cC_2}(Lh,Lg)}\left|\frac  {\zp-\alpha}{\zp-\beta}\right|\,\left| \frac {\zm-\beta}{\zm-\alpha }\right| \, .
\end{equation}
Note that, by Lemma  \ref{lem:E-to-cone} and the cone invariance, if $z\not\in E_{\cC_1}(h,g)$, then $zLg-Lh\in\cC_2$, hence $z\not\in  E_{\cC_2}(Lh,Lg)$. It follows that $E_{\cC_2}(Lh,Lg)\subset E_{\cC_1}(h,g)$, 
hence in equation \eqref{eq:direction_un} we have $m <  |\alpha|,|\beta| < M$.

Next, we choose $\zm,\zp$ to maximize the right-hand side of the above expression. Note that $\phi(z)=\frac {z-\beta}{z-\alpha }$ is a M\"obius transformation, hence it maps the circle $|z|=r$ into the circle of center
$c(r)=\frac{\beta\overline\alpha-r^2}{|\alpha|^2-r^2}$ and radius $R(r)=\frac{ r|\beta-\alpha|}{|\,|\alpha|^2-r^2\,|}$.\footnote{ Indeed, the reader can easily check that $|\phi^{-1}(x+iy)|=r$ implies 
\[
x^2+y^2-2\frac{\Re(\beta\overline\alpha-r^2)x+\Im(\beta\overline\alpha-r^2)y}{|\alpha|^2-r^2}+\frac{|\beta|^2-r^2}{|\alpha|^2-r^2}=0 .
\]
}
As observed above, the entire circles of radius $m$ and $M$ lie in the
complement of $E_{\cC_1}(h,g)$.
Thus,
\[
\begin{split}
\sup_{|\zm|=m}\left|\frac  {\zm-\beta}{\zm-\alpha}\right|&
=|c(m)|+R(m)=\frac{|\beta \overline\alpha-m^2|+m|\alpha-\beta|}{|\alpha|^2-m^2} \\
& \geq \frac{|\beta||\alpha| - m^2 + m|\beta| - m|\alpha|}{|\alpha|^2 - m^2} 
= \frac{|\beta|-m}{|\alpha|-m}.
\end{split}
\]
Reversing the roles of $\alpha$ and $\beta$ yields the analogous estimate on $C(M)$,
\[
\sup_{|\zp|=M}\left|\frac  {\zp-\alpha}{\zp-\beta}\right| = \frac{|M^2- \alpha \overline\beta|+M|\beta-\alpha|}{M^2-|\beta|^2} \geq \frac{M-|\alpha|}{M-|\beta|} .
\]
Substituting in \eqref{eq:direction_un} yields
\[
e^\Delta\geq \sup_{\alpha,\beta\in E(Lh,Lg)}\frac{(M-|\alpha|)(|\frac{\beta}\alpha||\alpha|-m)}{(M-|\frac{\beta}\alpha||\alpha|)(|\alpha|-m)}.
\]
Remark that the above supremum occurs when $|\alpha| < |\beta|$, which we will assume going forward.  With this condition, the ratio 
$\frac{ (M- |\alpha|)(|\beta| - m)}{(M - |\beta|)(|\alpha| - m)} \ge 1$.
We can then study the function $\phi(t)=\frac{(M-t)(|\frac{\beta}\alpha| t-m)}{(M-|\frac{\beta}\alpha|t)(t-m)}$  for $m<t<M \frac{|\alpha|}{|\beta|}$.
This function is positive, has vertical asymptotes at the endpoints of its domain\footnote{Note that
$t = |\alpha|$ is indeed in this domain since, recalling also \eqref{eq:Mm_def}, $m < |\alpha|$ by definition and $M|\alpha|/|\beta| > |\alpha|$ since $M > |\beta|$.}  and has a minimum at $t_0=\sqrt{\frac{m M|\alpha|}{|\beta|}}\in (m,M \frac{|\alpha|}{|\beta|})$. 
Hence,
\[
e^\Delta\geq  \sup_{\alpha,\beta\in E_{\cC_2}(Lh,Lg)} \phi(t_0)=\sup_{\alpha,\beta\in E_{\cC_2}(Lh,Lg)}\frac{\left(\sqrt{\frac{M|\beta|}{m|\alpha|}}-1\right)^2}{\left(\sqrt{\frac M m}-\sqrt{\frac{|\beta|}{|\alpha|}}\right)^2}.
\]
Since the above right-hand side is increasing with respect to $\frac{|\beta|}{|\alpha|}$ for $0<\frac{|\beta|}{|\alpha|} <\frac M m$, and, by definition, $ \sup_{\alpha,\beta\in E_{\cC_2}(Lh,Lg)}\frac{|\beta|}{|\alpha|}= e^{\delta_{\cC_2(Lg,Lh)}}=:e^{\theta_2}$, we have, setting 
$e^{\theta_1}:= \frac{M}{m} = e^{\delta_{\cC_1}(h,g)}$, and recalling that by construction, $\theta_1 > \theta_2$,
\[
\begin{split}
e^{\frac{\Delta}2}&\geq \frac{e^{\frac{\theta_2+\theta_1}2}-1}{e^{\frac{\theta_1}2}-e^{\frac{\theta_2}2}}
=\frac{e^{\frac{\theta_2+\theta_1}4}-e^{-\frac{\theta_2+\theta_1}4}}{e^{\frac{\theta_1-\theta_2}4}-e^{-\frac{\theta_1-\theta_2}4}}=\frac{\sinh\left(\frac{\theta_2+\theta_1}4\right)}{\sinh\left(\frac{\theta_1-\theta_2}4\right)}\\
&=\frac{\sinh\left(\frac{\theta_1}4\right)\cosh\left(\frac{\theta_2}4\right)+\cosh\left(\frac{\theta_1}4\right)\sinh\left(\frac{\theta_2}4\right)}{\sinh\left(\frac{\theta_1}4\right)\cosh\left(\frac{\theta_2}4\right)-\cosh\left(\frac{\theta_1}4\right)\sinh\left(\frac{\theta_2}4\right)}=\frac{\tanh\left(\frac{\theta_1}4\right)+\tanh\left(\frac{\theta_2}4\right)}{\tanh\left(\frac{\theta_1}4\right)-\tanh\left(\frac{\theta_2}4\right)}.
\end{split}
\]
The above implies 
\[
\tanh\left(\frac{\theta_2}4\right)\leq \frac{e^{\frac{\Delta}2}-1}{1+e^{\frac{\Delta}2}}\tanh\left(\frac{\theta_1}4\right)
=\tanh\left(\frac{\Delta}4\right)\tanh\left(\frac{\theta_1}4\right)
\leq \tanh\left(\frac{ \theta_1  \tanh\left(\frac{\Delta}4\right)}4\right)
\]
since, for $\gamma\in (0,1)$ and $x\geq 0$, $\tanh(\gamma x)\geq \gamma\tanh(x)$. Because the hyperbolic tangent is increasing, we finally have $\theta_2\leq\tanh\left(\frac{\Delta}4\right)\theta_1$,
and the theorem follows.
\end{proof}

Next, we provide a simplified but less optimal version of \cite[Theorem 4.5]{Du11} that suffices for our purposes.\\
For $i\in\{1,2\}$, let  $\cB_{i,\bR}$, $\cB_{i,\bC}$, $\cS_i$, $\bbe_i$ and $\bbm_i$ be as in Section \ref{sec:realcones} and Appendix \ref{app:complexification}.
We assume \eqref{eq:b_norms} and \eqref{eq:maybetoomuch} hold for both $i=1$ and $i=2$.
For $A \in L(\cB_{1,\bR}, \cB_{2,\bR})$, set
\[
B := d_H(A \bbe_1, \bbe_2).
\]
We start by relating the real diameter to the complex one.

\begin{lem}\label{lem:diam_rel}
Let $A\in L(\cB_{1,\bR},\cB_{2,\bR})$ such that $A(\cC_{1,\bR})\subset \cC_{2,\bR}$ with 
\[
\diam_{d_H}(A(\cC_{1,\bR}))=\Delta_\bR<\infty
\]
then\footnote{A more sophisticated argument would yield a better constant, independent of $\kappa$, see \cite{Du11}, but this estimate is more than sufficient for our needs and in the spirit to present the argument in its simplest form.   We remark also that if $\bbe_2 \in A( \cC_{1,\bR})$, as is often the 
case, then one can take $B=0$ in the statement of the lemma. \label{foot:B} }
\[
\Delta_{\bC} :=\diam_{\delta_{\cC_{2,\bC}}}(A(\cC_{1,\bC}))\leq 8\Delta_\bR + 8B  + 2 \ln[  3\sqrt{2} \kappa^{-2} ].
\]
\end{lem}
\begin{proof}
By \eqref{eq:hilbert_real} we have using the triangle inequality, for all $h\in\cC_{1,\bR}$ and $\ell\in\cC_{2,\bR}'$
\[
e^{ \Delta_\bR + B}\geq \frac{\ell(Ah)\bbm_2(\bbe_2)}{\bbm_2(Ah)\ell(\bbe_2)}\geq e^{ -\Delta_\bR -B}
\]
which, recalling \eqref{eq:maybetoomuch} and the normalization $\bbm_i(\bbe_i)  =1$, implies
\begin{equation}\label{eq:up_lower}
\ell(\bbe_2)\|Ah\|e^{ \Delta_{\bR} + B}\geq \ell(Ah)\geq \ell(\bbe_2)\kappa\|Ah\|e^{-\Delta_{\bR} -B}.
\end{equation}
Accordingly, for all $h=z(x+iy)$, $z\in \bC\setminus\{0\}$, $x,y \in \cC_{1,\bR}$, and $\ell,p\in \cC_{2,\bR}'$
\[
\Re(\ell(A h)\overline{p(Ah)})=|z|^2\left[\ell(Ax)p(Ax)+\ell(Ay)p(Ay)\right]\geq 0.
\]
It follows from Lemma~\ref{lem:dual_comp}-(a) that $A(\cC_{1,\bC})\subset \cC_{2,\bC}$, and, for all $\tilde\ell, \tilde p\in \hat \cC_{2,\bC}'$
and all $z\in\bC$ and $x,y \in \cC_{1,\bR}$, setting $\ell=|\tilde\ell(\bbe_2)|^{-1}\tilde\ell=\ell_1+i\ell_2$ and $p=|\tilde p(\bbe_2)|^{-1}\tilde p=p_1+ip_2$
\begin{equation}
\label{eq:top_bottom}
\left|\frac{z\tilde\ell(A(x+iy))\tilde p(\bbe_2)}{\tilde\ell(\bbe_2)z\tilde p(A(x+iy))}\right|^2\leq
\frac{\left[\ell_1(Ax)-\ell_2(Ay)\right]^2+\left[\ell_2(Ax)+\ell_1(Ay)\right]^2}
{\left[p_1(Ax)-p_2(Ay)\right]^2+\left[p_2(Ax)+p_1(Ay)\right]^2}.
\end{equation}

We use \eqref{eq:up_lower} to bound the numerator,
\[
\begin{split}
\left[\ell_1(Ax)-\ell_2(Ay)\right]^2+\left[\ell_2(Ax)+\ell_1(Ay)\right]^2
& \le 2 \big( \ell_1(Ax)^2 + \ell_2(Ay)^2  \big)+ 2 \big(\ell_2(Ax)^2+\ell_1(Ay)^2 \big) \\
& \le 2 e^{2 (\Delta_\bR + B)} \big( \|Ax\|^2 + \| Ay \|^2 \big) \, ,
\end{split}
\]
where we have used the normalization $\ell_1(\bbe_2)^2 + \ell_2(\bbe_2)^2 = 1$.

To derive a lower bound on the denominator in \eqref{eq:top_bottom}, we argue by cases.
Fix $b \ge 1$ and $c>0$, $c < \frac{e^{-2\Delta_{\bR}-2B} \kappa}{\sqrt{e^{-4 \Delta_{\bR}-4B} \kappa^2 + 4b^2} } < 1/\sqrt{5}$.

\medskip
\noindent
{\em Case 1: $p_1(\bbe_2)^2, p_2(\bbe_2)^2 \ge c^2$.}
Note that, by definition, there exists $\sigma_i\in\{-1,+1\}$ such that $\sigma_i p_i\in\cC_{2,\bR}'$. 
If $\sigma_1\sigma_2=1$, then, recalling \eqref{eq:up_lower}, 
\[
\begin{split}
&\left[p_1(Ax)-p_2(Ay)\right]^2+\left[p_2(Ax)+p_1(Ay)\right]^2\geq \left[p_2(Ax)+p_1(Ay)\right]^2\\
&\geq p_2(Ax)^2+p_1(Ay)^2\geq \kappa^2 c^2 (\|Ax\|^2  + \|Ay\|^2 ) e^{-2\Delta_{\bR}-2B}.
\end{split}
\]
while, If $\sigma_1\sigma_2=-1$, then similarly,
\[
\begin{split}
&\left[p_1(Ax)-p_2(Ay)\right]^2+\left[p_2(Ax)+p_1(Ay)\right]^2\geq \left[p_1(Ax)-p_2(Ay)\right]^2\\
&\geq p_1(Ax)^2+p_2(Ay)^2\geq \kappa^2 c^2 (\|Ax\|^2  + \|Ay\|^2  ) e^{-2\Delta_{\bR}-2B}.
\end{split}
\]

\medskip
\noindent
{\em Case 2: Either $p_1(\bbe_2)^2 < c^2$ or $p_2(\bbe_2)^2 < c^2$.}
We shall assume $p_1(\bbe_2)^2 <c^2$, the argument for the other alternative being analogous.
In this case then, $p_2(\bbe_2)^2 > 1-c^2$ due to the normalization of $p$.

We further refine into subcases.

{\em Subcase A: $\| Ay \| b^{-1} \le \| Ax \| \le b \|Ay \|$.}
Then by choice of $c$ and using \eqref{eq:up_lower},
\[
\begin{split}
 \left[p_1(Ax)-p_2(Ay)\right]^2 &\ge \left[ \sqrt{1-c^2} \kappa e^{-\Delta_\bR-B} \| Ay \| - c e^{\Delta_\bR+B} b \| Ay \| \right]^2 \\
 &= \| Ay \|^2 \left( \sqrt{1-c^2} \kappa e^{-\Delta_\bR-B} - b c e^{\Delta_\bR+B} \right)^2 \, .
  \end{split}
\]
Adding the analogous estimate for the second term yields the lower bound,
\[
 \left[p_1(Ax)-p_2(Ay)\right]^2 + \left[p_2(Ax)+p_1(Ay)\right]^2 \ge  \left( \| Ax \|^2 + \|A y \|^2 \right) \left( \sqrt{1-c^2} \kappa e^{-\Delta_\bR-B} - b c e^{\Delta_\bR+B} \right)^2 \, .
\]

{\em Subcase B: Either $\| Ay \| < b^{-1} \| A_x \|$ or $\| Ax \| < b^{-1} \| Ay \|$.}
Again, it is sufficient to argue only one of the alternatives.  Suppose that $\| Ay \| < b^{-1} \| Ax \|$.
Then,
\[
\begin{split}
 &\left[p_1(Ax)-p_2(Ay)\right]^2 + \left[p_2(Ax)+p_1(Ay)\right]^2 \ge
\left[p_2(Ax) + p_1(Ay)\right]^2 \\
&\hskip1cm\ge \| A x \|^2 \left( \sqrt{1-c^2} \kappa e^{-\Delta_\bR-B} - b^{-1} c e^{\Delta_\bR+B} \right)^2\\
&\hskip1cm\ge \frac{b^2}{b^2+1} \left( \| Ax \|^2 + \| Ay \|^2 \right) \left( \sqrt{1-c^2} \kappa e^{-\Delta_\bR-B} - b^{-1} c e^{\Delta_\bR+B} \right)^2 .
\end{split}
\]
We proceed to optimize our choice of $b$ and $c$ to obtain a common lower bound from Cases
1 and 2.  Let us choose $b=\sqrt{2}$.  Then set 
\[
c = \frac{\kappa e^{- 2\Delta_{\bR}-2B}}{\sqrt{\kappa^2 e^{-4\Delta_\bR-4B} + 4b^2}} \implies
\sqrt{c^{-2} - 1} = 2 b \kappa^{-1} e^{2 \Delta_{\bR}+2B} \, .
\]
This implies that the lower bound in Subcase A is at least,
\begin{equation}
\label{eq:subcase A}
\left( \| Ax \|^2 + \| Ay \|^2 \right) \frac{ \kappa^2 c^2}{4} \left( c^{-2} -1 \right) e^{-2 \Delta_\bR - 2B} \, .
\end{equation}
Moreover, the lower bound in Subcase B is at least as large since with $b=\sqrt{2}$,
\[
\frac{b}{\sqrt{b^2+1}} \left(\sqrt{c^{-2} - 1} - \kappa^{-1} b^{-1} e^{2 \Delta_\bR + 2B} \right)
 =  \sqrt{c^{-2}-1} \left( 1 - \frac{1}{2b^2} \right) \frac{b}{\sqrt{b^2+1}} > \frac 12 \sqrt{c^{-2} -1} .
\]
It remains to compare the lower bounds from Subcase A and Case 1.  Due to \eqref{eq:subcase A},
the lower bound from Case 1 is smaller if
\[
c^2  < c^2 \frac{1}{4} (c^{-2} -1) \implies c < \frac{1}{\sqrt{5}} \, ,
\]
which is true given our choice of $c$ whenever $b \ge 1$.

Finally, putting together this common lower bound with our upper bound for the numerator in
\eqref{eq:top_bottom}, yields,
\begin{equation}
\label{eq:narrow bound}
 \left|\frac{z \tilde\ell(A(x+iy)) \tilde p(\bbe_2)}{\tilde\ell(\bbe_2)z \tilde p(A(x+iy))}\right|\leq \sqrt 2\kappa^{-1}  c^{-1} e^{2\Delta_{\bR} +2B} 
 \le 3\sqrt{2} \kappa^{-2} e^{4 \Delta_\bR +4B} \, . 
\end{equation}
Since
\begin{equation}
\label{eq:rad double}
 \diam_{\delta_{\cC_{2,\bC}}}(A(\cC_{1,\bC}))\leq 2\sup_{h\in\cC_{1,\bC}}\delta_{\cC_{2, \bC}}(A h,\bbe_2),
\end{equation}
we have the lemma.
\end{proof}

We finally provide the key result that is needed in the applications to verify the hypothesis on the finiteness of the complex diameter in Theorem \ref{thm:complex_contraction}.
\begin{thm} \label{thm:needed} 
Let, for $i\in\{1,2\}$, $\cB_{i,\bR}$, $\cB_{i,\bC}$, $\cS_i$, $\bbe_i$ and $\bbm_i$ be as in the introduction to this section (in particular, the $\bbe_i$ and $\bbm_i$ satisfy equations \eqref{eq:archimedean} and \eqref{eq:maybetoomuch}, respectively).
Let $\cL\in L(\cB_{1,\bR},\cB_{2,\bR})$ and $\cL_\bC\in L(\cB_{1,\bC},\cB_{2,\bC})$.
Assume that $\cL(\cC_{1,\bR})\subset \cC_{2,\bR}$ and $\diam_H(\cL(\cC_{1,\bR})):=\Delta_\bR<\infty$. 
If there exists $\ve \in (0, \frac{\kappa^2}{12 \sqrt 2}e^{ -2\Delta_\bR -2B})$, such that
for all $\ell \in \cS_2$ and all $h \in \cC_{1,\bR}$,
\begin{equation}\label{eq:complex_bound}
\left|\ell(\cL_{\bC} h)-\ell(\cL h)\right|\leq \ve \ell(\cL h).
\end{equation}
Then $\cL_{\bC}(\cC_{1,\bC})\subset \cC_{2,\bC}$, and we have
\[
\diam_{\delta_{\cC_{2, \bC}}} (\cL_{\bC}(\cC_{1,\bC}))\leq 8\Delta_\bR + { 8B} +
 2 \ln[3\sqrt 2\kappa^{-2}]+ \tfrac{\sqrt{2}}{3}  \left(2+e^{\frac{\Delta_\bR}2}\right)\kappa^2e^{ -2\Delta_\bR -2B} . 
\]
\end{thm}
\begin{proof}
For all $x,y\in \cC_{1,\bR}$ and $\ell\in\cS_2$ we have using \eqref{eq:complex_bound},
\begin{equation}\label{eq:real_comp_compa}
\begin{split}
\left|\ell(\cL_{\bC} (x+iy))-\ell(\cL (x+iy))\right|\leq& \left|\ell(\cL_{\bC} x)-\ell(\cL x)\right|
+\left|\ell(\cL_{\bC} y)-\ell(\cL y)\right|\\
\leq& \ve [\ell(\cL x)+ \ell(\cL y)]\leq \ve  \sqrt 2 |\ell(\cL (x+iy))|.
\end{split}
\end{equation}
It follows that for all $h\in\cC_{1,\bC}$ and $\ell,p\in\mathring\cC_{2,\bR}'$, provided $\ve\leq \frac{1}{\sqrt 2}$,
\[
\begin{split}
\Re(\ell(\cL_{\bC} h)\overline{p(\cL_{\bC} h)})
& =
\Re\Bigg(\ell(\cL h)\overline{p(\cL h)}  + \ell((\cL_{\bC}-\cL) h )\overline{p((\cL_{\bC}-\cL) h )}  \\
& +\ell((\cL_{\bC}-\cL) h)\overline{p(\cL h )}+\ell(\cL h )\overline{p((\cL_{\bC}-\cL) h )} \Bigg)\\
&\geq \Re\left(\ell(\cL h )\overline{p(\cL h)}\right)
- 3 \sqrt 2\ve|\ell(\cL h )||{p(\cL h )}|\\
&\geq \Re\left(\ell(\cL h )\overline{p(\cL h)}\right)
-6\sqrt 2\ve \ell(\bbe_2)p(\bbe_2)\|\cL h\|\|\cL h \|
\end{split}
\]
where, in the first inequality we have used \eqref{eq:real_comp_compa} and in the second inequality we have
used Lemma \ref{lem:dual_comp}(c) and Lemma~\ref{lem:upper_dual}. 
Since we can assume $h=x+iy$ with $x,y\in \cC_{1,\bR}$, recalling \eqref{eq:up_lower} we have
\[
\begin{split}
\Re\left(\ell(\cL h )\overline{p(\cL h )}\right) & = \ell(\cL(x))p(\cL x )+\ell(\cL y )p(\cL y)\\
& \geq \ell(\bbe_2)p(\bbe_2)\kappa^2 e^{ -2\Delta_\bR -2B} \left( \|\cL x\|^2 + \|\cL y \|^2 \right) \\
& \geq  2^{-1} \ell(\bbe_2)p(\bbe_2)\kappa^2 e^{-2\Delta_\bR-2B} \| \cL h \|^2 \, , 
\end{split}
\]
where we have used \eqref{eq:norms_equiv}.  Hence, $\Re(\ell(\cL_{\bC}(h))\overline{p(\cL_{\bC}(h))})>0$, provided $\ve < \frac{\kappa^2}{ 12 \sqrt 2}e^{-2\Delta_\bR - 2B}$. 

Note that if $\ell \in \cC_{2,\bR}'$, then $\ell+\alpha\bbm\in\mathring\cC_{2,\bR}'$ for each $\alpha>0$.  Thus the above implies $\Re(\ell(\cL_{\bC}(h))\overline{p(\cL_{\bC}(h))}) \ge 0$ for
all $\ell, p \in \cC_{2, \bR}'$. 
Applying Lemma~\ref{lem:dual_comp}(a), 
 we have $\cL_{\bC}(\cC_{1,\bC})\subset \cC_{2,\bC}$.

Finally, for all $\ell,p\in\cC_{2,\bR}'$ and $h\in\cC_{1,\bC}$, it follows by \eqref{eq:real_comp_compa} that,
\begin{equation}\label{eq:real_ratios}
\begin{split}
\left|\frac{\ell(\cL_{\bC}h)p(\cL h)}{\ell(\cL h)p( \cL_{\bC}h)}\right|&\leq 
\frac{|\ell(\cL h)p(\cL h)|+|\ell([\cL_{\bC}- \cL  ]h)p(\cL h)|}{|\ell(\cL h)p( \cL h)|
-|\ell(\cL h)p([\cL_{\bC}- \cL  ]h)|}\\
&\leq \frac{1+  \sqrt 2\ve}{1- \sqrt 2\ve}  \leq 1 + 4 \ve, 
\end{split}
\end{equation}
since we assumed $\ve < \frac 1{ 12 \sqrt{2}}$.

The above holds for $\ell,p\in\cC_{2,\bR}'$, while it is needed for $\ell,p\in\hat\cC'_{\bC}$. Let $q\in\hat\cC'_{\bC}$, then, recalling Lemma \ref{lem:dual_comp}, $q=\ell_1+ip_1$, with $\ell_1,p_1\in\mathring\cC'_{2,\bR}\cup(-\mathring\cC'_{2,\bR})$. We treat explicitly only the case $\ell_1,p_1\in\mathring\cC'_{2,\bR}$, the other being identical. Set $a_1=\ell_1(\cL x)$, $b_1=\ell_1(\cL y)$, $a_2=p_1(\cL x)$, $b_2=p_1(\cL y)$, where $h=x+iy$, $x,y\in\cC_{1,\bR}$. Then $q(\cL h)=(a_1-b_2)+i(b_1+a_2)$ and, by \eqref{eq:complex_bound}, $|q((\cL_\bC-\cL)h)|\leq \ve(a_1+b_1+a_2+b_2)$. On the other hand, $\diam_H(\cL(\cC_{1,\bR}))\leq\Delta_\bR$ implies $a_1b_2\leq e^{\Delta_\bR}a_2b_1$, hence
\[
|q(\cL h)|\geq b_1+a_2\geq 2\sqrt{a_2b_1}\geq 2e^{-\frac{\Delta_\bR}2}\sqrt{a_1b_2}\geq 2e^{-\frac{\Delta_\bR}2}\min\{a_1,b_2\},
\]
while $\max\{a_1,b_2\}\leq |a_1-b_2|+\min\{a_1,b_2\}\leq |q(\cL h)|(1+\frac 12 e^{\frac{\Delta_\bR}2})$. Thus, 
\[
a_1+b_1+a_2+b_2\leq (2+e^{\frac{\Delta_\bR}2})|q(\cL h)|.
\] 
Consequently, equation \eqref{eq:real_ratios} holds for all $\ell,p\in\hat \cC'_{2,\bC}$ with $\ve$ replaced by
\[
\ve_\bC:=\ve\left(2+e^{\frac{\Delta_\bR}2}\right).
\]

The above, by Definition \eqref{eq:gauge_dist} and Lemma \ref{lem:are_the_same}, implies 
\begin{equation}
\label{eq:needed later}
\begin{split}
\delta_{\cC_{2, \bC}}(\cL_{\bC}h,\bbe_2)&=\sup_{\ell,p \in \hat\cC'_{\bC}}\ln \frac{|\ell(\cL_{\bC}h)p(\bbe_2)|}{|\ell(\bbe_2)p( \cL_{\bC}h)|} \\
&\leq \sup_{\ell,p}\ln \frac{|\ell(\cL h)p(\bbe_2)|}{|\ell(\bbe_2)p( \cL h)|}+\sup_{\ell,p}\ln \frac{|\ell(\cL_{\bC}h)p(\cL h)|}{|\ell(\cL h)p( \cL_{\bC}h)|}\\
&\leq \delta_{\cC_{2, \bC} }(\cL h,\bbe_2)+ 4 \ve_\bC
\leq 4\Delta_\bR + 4B + \ln[3 \sqrt 2\kappa^{-2}] + 4  \ve_\bC \\
&\leq 4\Delta_\bR + 4B + \ln[3\sqrt 2\kappa^{-2}]+   \tfrac{1}{ 3 \sqrt{2} }
\left(2+e^{\frac{\Delta_\bR}2}\right) \kappa^2e^{-2\Delta_\bR-2B},
\end{split}
\end{equation}
where, in the next to last inequality, we have applied \eqref{eq:narrow bound}. We conclude using
\eqref{eq:rad double}.
\end{proof}
We conclude with a comment on  \eqref{eq:complex_bound} that may facilitate checking it. 

For each $h\in\cC_{1,\bR}$ let $a(h)=\Re(\cL_{\bC} h), b(h)=\Im(\cL_{\bC} h)$. By definition $a,b$ are real linear operators on $\cB_{2,\bR}$ and $\cL_{\bC} h =a(h)+i b(h)$. 

\begin{lem}\label{lem:pert_cond}
If, for each $h\in\cC_{1,\bR}$,
\[
\begin{split}
&\ve\cL h\pm2\left[\cL h- a(h)\right]\in\cC_{2,\bR}\\
&\ve\cL h \pm 2 b(h) \in\cC_{2,\bR}
\end{split}
\]
 then condition \eqref{eq:complex_bound} is satisfied.
\end{lem}
\begin{proof}
By hypothesis, for each $\ell\in\cC_{2,\bR}'$ we have
\[
\begin{split}
&\ve\ell(\cL h)\geq 2\left|\ell(\cL h- a(h))\right|\\
&\ve\ell(\cL h)\geq 2\left|\ell(b(h) )\right|.
\end{split}
\]
Accordingly,
\[
\ve\ell(\cL h)\geq\left(|\ell(\cL h- a(h))|+|\ell(- ib(h))|\right)\geq |\ell(\cL h- \cL_{\bC}(h))|.
\]
\end{proof}


\subsection{Cone Contraction Implies Loss of Memory}
\label{sec:memory loss}\ \\

In this section, we show that cone contraction implies loss of memory.
This provides a tool to easily verify, in some cases, condition (O-\ref{cond:lambda_multiplier2}).

Consider complex Banach spaces $\bB_k$, cones $\cC_k = \cC_{\bC, k} \subset \bB_k$, vectors $\bbe_k\in\bB_k$ as described in Sections~\ref{sec:setting}  and \ref{sec:complexcones}.
Suppose there exist functionals $\bbm_k\in\bB_k'$ 
satisfying \eqref{eq:maybetoomuch} with
\begin{enumerate}
  \item[a)] $\displaystyle \inf_{k \ge 0} \kappa_k = \bar \kappa > 0$. 
\end{enumerate}
In addition, 
there are operators $\bL_{k}:\bB_k\to\bB_{k+1}$ satisfying the following:
\begin{enumerate}
\item[b)] there exists $K>0$ such that, for all $k,j\in\bN$, we have 
\[
 \left|\bbm_{k+1}\left(\bL_{k}\cdots \bL_{j} \bbe_j\right)\right|\leq K;
\]
  \item[c)] the operators $\bL_{k}$ satisfy $\bL_{k}(\cC_{k}) \subset \cC_{k+1}$, and the
  diameter in the complex projective metric of $\bL_{k}(\cC_{k})$ in $\cC_{k+1}$
  is uniformly bounded by $\Delta< \infty$.
\end{enumerate}
Define $\bL_{k,j} :=\bL_{k}\cdots \bL_{j}$ and  
\begin{equation}\label{eq:alpha-h-ell}
\begin{split}
&\boldsymbol{\alpha}_{k,j}=\bbm_k(\bL_{k-1,j} \bbe_j)\\
&\boldsymbol{h}_{k,j}=\frac{\bL_{k-1,j} \bbe_j}{\bbm_k(\bL_{k-1,j}\bbe_j)}\\
&\boldsymbol{\ell}_{k,j}(h)=\frac{\bbm_k\left(\bL_{k-1,j} h\right)}{\bbm_k(\bL_{k-1,j} \bbe_j)}.
\end{split}
\end{equation}
The next lemma will be helpful to check condition (O-\ref{cond:lambda_multiplier}).

\begin{lem}[Uniform Boundedness]
\label{lem:boundedness}
For each $h\in\bB_j$, we have
\[
\begin{split}
&|\boldsymbol{\alpha}_{k,j}|\leq K,\quad \|\boldsymbol{h}_{k,j}\|_k \leq \frac{\sqrt 2}{\bar\kappa},\quad \|\boldsymbol{\ell}_{k,j}\|_j'\leq \sqrt 2\\
&\|\bL_{k,j}h\|_{k+1} \leq \frac{2(1+2\sqrt 2)K}{\bar\kappa}\|h\|_j.
\end{split}
\]
\end{lem}
\begin{proof}
By assumption $|\boldsymbol{\alpha}_{k,j}|\leq K$.
Since, by assumption, $\bbm_k$ satisfies \eqref{eq:m_lower}, by Definition~\ref{def:complex cone}, $\bbm \in \cC_{\bC}'$ so that, for all $h \in \cC_{\bC}$,
\begin{equation}\label{eq:preliminary_bound}
 \left|\bbm_k(\bL_{k-1,j} h)\right| \geq \frac{\bar \kappa}{\sqrt 2}\|\bL_{k-1,j} h\|_k > 0.
\end{equation}
Hence, $\|\boldsymbol{h}_{k,j}\|_k \leq \frac{\sqrt 2}{\bar\kappa}$. In addition, $\boldsymbol{\ell}_{k,j}(h) \ne 0$ for all $h \in \cC_{\bC}$.  Accordingly,
$\boldsymbol{\ell}_{k,j} \in \cC_{\bC}'$ and so by Lemma~\ref{lem:upper_dual} $\|\boldsymbol{\ell}_{k,j}\|'\leq \sqrt 2$. Hence, recalling \eqref{eq:m_lower}, for each $h\in\cC_{\bC,j}$,
\[
\|\bL_{k-1,j}h\|_k \leq \frac{\sqrt 2}{\bar\kappa} \left| \bbm_k\left(\bL_{k-1,j}h\right) \right|=\frac{\sqrt 2}{\bar\kappa} \left| \boldsymbol{\ell}_{k,j}\left(h\right)\bbm_k\left(\bL_{k-1,j}\bbe_j\right) \right| \leq \frac{2K}{\bar\kappa}\|h\|_j.
\]
Finally, for an arbitrary $h\in\bB_j$, $h\neq 0$, and $\eta>0$, set $\bar h=h+(1+\eta)\sqrt 2\|h\|_j\bbe_j$. As in Lemma \ref{lem:interior}, $\bar h, \bbe_j\in\cC_{\bC,j}$, hence, applying the above to $\bar h$ and to $\bbe_j$, and letting $\eta\downarrow 0$,
\[
\|\bL_{k-1,j}h\|_k\leq \|\bL_{k-1,j}\bar h\|_k+\sqrt 2\|h\|_j\|\bL_{k-1,j}\bbe_j\|_k\leq \frac{2K}{\bar\kappa}\left[(1+\sqrt 2)+\sqrt 2\right]\|h\|_j,
\]
which proves the last inequality.
\end{proof}
The following lemma will be helpful to check condition (O-\ref{cond:lambda_multiplier2}).

 \begin{lem}[Loss of Memory] 
 \label{lem:memory loss}
 Let $\theta = \tanh(\Delta/4)$, then for each $k,j\in\bN$, $h \in \bB_j$,
 \[
\begin{split}
\| \bL_{k-1,j} h - \boldsymbol{\alpha}_{k,j} & \boldsymbol{\ell}_{k,j}(h) \boldsymbol{h}_{k,j} \|_{k}
\le  \frac{2(1+\sqrt{2}) }{\bar\kappa \theta^2} \Delta \theta^{k-j} |\boldsymbol{\alpha}_{k,j} |\, \| h\|_j. 
\end{split}
\]
 \end{lem}
 \begin{proof}
 Applying Lemma \ref{lem:compare} and Theorem \ref{thm:complex_contraction}, we have for all $h,g\in\cC_{j}$,
\[
\begin{split}
&\left\|\frac{ \bL_{k,j} h}{\bbm_{k+1}(\bL_{k,j} h)}-\frac{ \bL_{k,j} g}{\bbm_{k+1}(\bL_{k,j}g)}\right\|_{k+1} \leq \frac{\sqrt 2}{\bar \kappa}\delta_{\bC,k+1}(\bL_{k,j} h,\bL_{k,j} g)\\
&\leq \frac{\sqrt 2}{\bar\kappa} (\tanh(\Delta/4))^{k-j-1}\delta_{\bC, j+1}( \bL_{j} h, \bL_{j} g)\\
& =\frac{\sqrt 2}{\bar\kappa} \theta^{k-j-1}\delta_{\bC,j+1}( \bL_{j}  h, \bL_{j} g)\leq \frac{ \sqrt 2}{\bar\kappa} \theta^{k-j-1}\Delta.
\end{split}
\]

By definition $\boldsymbol{\ell}_{k,j}(\bbe_j)=1$ and $\bbm_k(\boldsymbol{h}_{k,j})=1$.
Moreover, 
by Lemma~\ref{lem:boundedness} and Lemma~\ref{lem:upper_dual}, 
$|\boldsymbol{\alpha}_{k,j}| \le K$, $\| \boldsymbol{h}_{k,j} \|_k \le \sqrt{2}/\bar \kappa$
and 
$\|\boldsymbol{\ell}_{k,j} \|_j' \le \sqrt{2}$.
Hence, recalling \eqref{eq:alpha-h-ell}, for all $h\in\cC_{j}$ we have
\begin{equation}
\label{eq:prelim contract 6}
\begin{split}
\left\| \bL_{k-1,j} h -\boldsymbol{\alpha}_{k,j}\boldsymbol{\ell}_{k,j}(h)\boldsymbol{h}_{k,j}\right\|_k&=\left\| \bL_{k-1,j} h - \bbm_k(\bL_{k-1,j} h)\boldsymbol{h}_{k,j}\right\|_k\\
&\leq \frac{\sqrt 2}{\theta^2\bar\kappa} \theta^{k-j}\Delta\left|\bbm_k(\bL_{k-1,j}h)\right|\\
&\leq \frac{2}{\theta^2\bar\kappa} \theta^{k-j}\Delta |\boldsymbol{\alpha}_{k,j} | \|h\|_j.
\end{split}
\end{equation}
For $h \in \bB_j$,  $h \ne 0$  (the case $h=0$ being trivial), and $\eta>0$, define $\bar h := h + (1+\eta) \sqrt{2} \| h \|_j \bbe_j$.  As in Lemma~\ref{lem:interior}, 
$\bar h, \| h \|_j \bbe_j \in \cC_{\bC,j}$. 
Lemma~\ref{lem:interior} requires a strict inequality, whence the factor $1+\eta$; the estimate below follows by letting $\eta\downarrow 0$.
Recalling $\boldsymbol{\ell}_{k,j}(\| h \|_j \bbe_j) = \| h \|_j$, as well as 
$\bbm_k(\bL_{k-1,j} \bar h) = \bbm_k(\boldsymbol{\alpha}_{k,j} \boldsymbol{\ell}_{k,j}(\bar h) \boldsymbol{h}_{k,j})$,
we apply \eqref{eq:prelim contract 6} to obtain,
\[
\begin{split}
\| \bL_{k-1,j} h - \boldsymbol{\alpha}_{k,j} & \boldsymbol{\ell}_{k,j}(h) \boldsymbol{h}_{k,j} \|_k
 = \| \bL_{k-1,j} \bar{h} - \boldsymbol{\alpha}_{k,j}\boldsymbol{\ell}_{k,j}(\bar{h}) \boldsymbol{h}_{k,j} \|_k \\
&  \le \| \bar h\|_j \frac{ 2}{\bar\kappa \theta^2} \Delta \theta^{k-j}  |\boldsymbol{\alpha}_{k,j}| \leq \| h\|_j \frac{2(1+\sqrt{2})}{\bar\kappa \theta^2} \Delta \theta^{k-j} |\boldsymbol{\alpha}_{k,j} | .
\end{split}
\]
\end{proof}


\subsection{A Useful Lower Bound}
\label{sec:lower bound}

In this section, we show under the assumptions of Theorem~\ref{thm:needed}, that the lower
bound $|\boldsymbol{\ell}_{k,j}(\boldsymbol{h}_{j,l}) | \ge K^{-1}$ holds uniformly in $j,k,l$, with 
$\boldsymbol{\ell}_{k,j}, \boldsymbol{h}_{j,l}$ defined as in \eqref{eq:alpha-h-ell}.  In addition to properties
(a)-(c) of Section~\ref{sec:memory loss}, we assume that there exist real operators $\cL_k : \cB_k \to \cB_{k+1}$
satisfying
\begin{itemize}
  \item[d)] $\cL_k(\cC_{k, \bR}) \subset \cC_{k+1, \bR}$ and $\diam_H(\cL_k(\cC_{k,\bR})) \le \Delta_{\bR} < \infty$ for all $k \ge 1$;
  \item[e)]  $\exists B<\infty$ such that $d_H(\cL_k(\bbe_k), \bbe_{k+1}) \le B$ for all $k \ge 1$; 
  \item[f)] for all $\ell \in \cS_{k+1}$ and all $h \in \cC_{k, \bR}$,
  $| \ell(\bL_k h) - \ell(\cL_k h)| \le \ve \ell(\cL_k h)$ for some $\ve \in \big(0, \frac{\bar\kappa^2}{12 \sqrt{2}} e^{-2\Delta_\bR -2B}\big)$.
\end{itemize}
Thus we are in the setting of Theorem~\ref{thm:needed}.

In this setting, we prove the lemma,
\begin{lem}
\label{lem:lower bound}
 Under the assumption (a)--(f), there exists $K >0$ such that $| \boldsymbol{\ell}_{k,j}(\boldsymbol{h}_{j,l}) | \ge K^{-1}$ for all $0 \le l \le j \le k$.  
 \end{lem}
 
 \begin{proof}
 By definition,
 \[
\boldsymbol{ \ell}_{k,j}(h) = \frac{\bbm_{ k}(\bL_{k-1,j} h)}{ \boldsymbol{\alpha}_{k,j} } \, .
 \]
Since $\bL_{k,j} \cC_{j,\bC} \subset \cC_{ k+1, \bC}$ and invoking \eqref{eq:m_lower} and property (a) of Section~\ref{sec:memory loss}, we have
$\boldsymbol{\ell}_{k,j}(h) \ne 0$ for all $h \in \cC_{j,\bC}$.  Thus $\boldsymbol{\ell}_{k,j} \in \cC_{j,\bC}'$ and since 
$\boldsymbol{h}_{j,l} \in \cC_{j, \bC}$, we may
apply Lemma~\ref{lem:partial converse} to obtain,
\[
|\boldsymbol{\ell}_{k,j}(\boldsymbol{h}_{j,l})| \ge e^{-\delta_{\cC_{j, \bC}}(\boldsymbol{h}_{j,l}, \bbe_j)} \tfrac{\bar\kappa}{\sqrt{2}} \| \boldsymbol{h}_{j,l} \|_j \, ,
\]
where we have used $\boldsymbol{\ell}_{k,j}(\bbe_j) = 1$.  
Recalling \eqref{eq:alpha-h-ell} and using Lemma \ref{lem:upper_dual}, and equation \eqref{eq:m-normalization}, we have
\[
1=\bbm_j(\boldsymbol{h}_{j,l})\leq \|\bbm_j\|_j'\|\boldsymbol{h}_{j,l}\|_j\leq \sqrt 2\|\boldsymbol{h}_{j,l}\|_j.
\]

It remains to show that $\delta_{\cC_{j, \bC}}(\boldsymbol{h}_{j,l}, \bbe_j) \le C$, where $C<\infty$ is independent of $j$ and $l$.
This follows from the cone contraction provided for each $\bL_i$.  In particular, note that
$\bL_{j-1, l} \cC_{l, \bC} \subset \bL_{j-1} \cC_{j-1, \bC}$, so that 
$\boldsymbol{h}_{j, l} \in \bL_{j-1} \cC_{j-1, \bC}$.  Thus, it suffices to show that 
$\delta_{j, \cC}(h, \bbe_j) \le C$ for each $h \in \bL_{j-1}\cC_{j-1, \bC}$ for some $C$ independent of $j$. 

To see this, apply Theorem~\ref{thm:needed}, and in particular the second line of \eqref{eq:needed later}, with 
$\cL_{\bC} = \bL_{j-1}$ and $\cL = \cL_{j-1}$ to obtain for each $h \in \cC_{j-1, \bC}$,
\[
\delta_{\cC_{j, \bC}}(\bL_{j-1} h, \bbe_j) \le \delta_{\cC_{j, \bC}}(\cL_{j-1} h, \bbe_j) + 4 \ve_\bC \, ,
\]
with $\ve_\bC=\ve(2+e^{\Delta_\bR/2})$ as in the proof of Theorem \ref{thm:needed}, applying properties (d) - (f).
Finally, applying Lemma~\ref{lem:diam_rel} and in particular \eqref{eq:narrow bound} with $A = \cL_{j-1}$ yields,
\[
\delta_{\cC_{j, \bC}}(\bL_{j-1} h, \bbe_j) \le 4 \Delta_{\bR} + 4B + \ln \big(3 \sqrt{2} \bar\kappa^{-2} \big) + 4 \ve_\bC \, ,
\]
completing the proof of the claim.
\end{proof}

\section{Application: Smooth Expanding Maps}
\label{sec:expanding}

\subsection{Sequential Expanding Maps}\ \\
\label{sec:expanding intro}
We consider the setting illustrated in Section~\ref{sec:expanding_p} and prove Theorem \ref{thm:expanding CLT}. To this end, we first check the conditions of Theorem~\ref{thm:main_bis}.
Theorem \ref{thm:expanding CLT} will then follow trivially as explained in Section~\ref{sec:CLT expanding}.

\subsubsection{ Verifying Conditions (C-\ref{cond:positivity})-(C-\ref{cond:powerbounded})} \ \\
Let $M$ be a smooth, connected $d$-dimensional  Riemannian manifold.  Without loss of generality, we can rescale the distance $d$ so that the diameter of $M$ is one and the Riemannian volume so that the volume of $M$ is also one.  We set $M_k=M$ and $\mu_k=$ Lebesgue. First of all, for $a> 0$, we consider the cones  $C_k = \cC_a$ for all $k \in \bN$,
\[
\cC_a =\left\{h\in C^1(M,\bR)\;:\; \|\nabla h(x)\| \leq a h(x) \,\, \forall x \in M \right\} \setminus \{ 0 \}.
\]
Note that (C-\ref{cond:positivity}) is trivially satisfied. Then setting $\bbe(x)=1\in \cC_a$, 
if $h\in C^1$, then $(a^{-1} \|\nabla h\|_\infty+\|h\|_\infty)\bbe+h\in \cC_a$;
that is, (C-\ref{cond:norm}) holds true.  
Indeed, it follows immediately from \eqref{eq:real cone norm} that
for $h \in C^1(M, \bC)$,
\[
\|g\|:=\|g\|_k= \sup_{x \in M} \{ a^{-1}\|\nabla g(x)\| + |g(x)| \}, 
\]
thus $\cB_k = C^1(M, \bR)$ and 
$\bB_k= C^1(M, \bC)$. 

Next, suppose 
$g\in \cC_a$ and $x_0\in M$ is such that\footnote{ Here and throughout this section, the integral
is taken with respect to the volume measure on $M$, which has been normalized to be a probability measure.} $g(x_0)=\int _M g$.\\
If $\gamma\in C^1([0,1],M)$ is a geodesic (parametrized by arclength) that connects $x$ and $x_0$, then
\[
\left| \frac{d}{dt}g(\gamma(t)) \right| \leq \|\nabla g(\gamma(t))\|\leq a g(\gamma(t)),
\]
hence 
\begin{equation}\label{eq:average bound}
e^a\int_M g\geq g(x)\geq e^{-a}\int_Mg \, , \quad \mbox{for all $x \in M$}.
\end{equation}
Defining the transfer operators $\cL_k$ as in \eqref{eq:transfer_op} with respect to the 
volume measure yields,
 \[
 \cL_k h(x)=\sum_{y\in f_k^{-1}(x)}\frac{h(y)}{\left|\det(D_yf_k)\right|} .
 \]
Hence, for some $\bar D\leq \Const \sup_{k \in \bN} \|D^2f_{k} \|_\infty = \Const A$,
\begin{equation}\label{eq:LY_exp}
\|\nabla \cL_k h\|\leq   \vartheta^{-1}\cL_k\|\nabla h\|+ \bar D \cL_k h\leq (\vartheta^{-1}a+ \bar D) \cL_k h
\end{equation}
which implies $\cL_k\cC_a\subset \cC_{\nu a}$, $\nu\in(\vartheta^{-1},1)$, for all $k \in \bN$, provided $a\geq \bar D(\nu-\vartheta^{-1})^{-1}$.

From this point forward, we fix $\nu \in (\vartheta^{-1}, 1)$ and $a>1$ such that $a > \bar D(\nu - \vartheta^{ -1})^{-1}$.

Now for each $h\in C^1$, and using \eqref{eq:average bound} since both
$\| h\|$ and $h + \| h \|$ are in $\cC_a$, 
\begin{equation}\label{exp_basic}
\begin{split}
| \cL_{k}\cdots \cL_j h|&\leq |\cL_{k}\cdots \cL_j (h+\|h\|)|+\|h\||\cL_{k}\cdots \cL_j 1|\\
&\leq e^a\int \cL_{k}\cdots \cL_j (h+2\|h\|)=  3e^a\|h\|,
\end{split}
\end{equation}
 which, iterating \eqref{eq:LY_exp}, proves (C-\ref{cond:powerbounded}).

One can easily compute the Hilbert metric $d_H$ and prove that for all $g\in\cC_{\nu a}$, $\nu\in (0,1)$,
\begin{equation}\label{eq:Delta_exp}
d_H(h,1)\leq 2 a+\ln\frac{1+\nu}{1-\nu}=:\Delta/2.
\end{equation}
So (C-\ref{cond:cone_cont}) holds true.

\subsection{Conditions (O-\ref{cond:multiplier}), (O-\ref{cond:lambda_multiplier}) (Applicability of Theorem~\ref{thm:main})}\ \\
Assume that we have a sequence of observables $g_k$ satisfying $\sup_k (\| g_k \|_\infty + \| \nabla g_k \|_\infty) = K < \infty$.  This implies in particular that $\sup_k \| g_k \| \le (1+a^{-1}) K$.

Let $\rho \in C^1$, $\rho \ge 0$ with $\int \rho \, d\mu_0=1$. 
Define $\hat g_k$ as in \eqref{eq:centering0}. For $n \in \mathbb{N}$,
define $\hat S_n$ as in
\eqref{eq:center sum}  and the variance $\sigma_n^2$ as in \eqref{eq:variance0}. 
Remark that $\sigma_n \le \Const \sqrt{n}$ by Lemma~\ref{lem:variance_up}.

With $n$ fixed and $k \le n$, define the weighted operators for $\lambda \in \mathbb{R}$ as in
\eqref{eq:twisted},
\[
\cL_{k,\lambda} h := \cL_{k, \lambda, n} h = \cL_k ( e^{i \sigma_n^{-1} \lambda \hat g_k} h) \, .
\]
Note that, $\|\varphi h\|\leq \|\varphi\|\,\|h\|$, and, for each $\lambda_0>0$ and $|\lambda|\leq \lambda_0\sigma_n$,
\[
\left\|e^{i \sigma_n^{-1} \lambda \hat g_k}\right\|=\sup_{x\in M}\left\{1+a^{-1}\sigma_n^{-1}|\lambda|\,\|\nabla g_k(x)\|\right\}\leq 1+a^{-1}\lambda_0 K .
\]
Accordingly, by (C-\ref{cond:powerbounded}) and Lemma \ref{lem:complexification}, for each $h\in\bB_k=C^1(M,\bC)$ and $j\in\bN$,
\[
\left\|\cL_{k,\lambda,n}(g_k^j h)\right\|\leq \sqrt 2 C_*\left\|g_k^j\, e^{i \sigma_n^{-1} \lambda \hat g_k}\, h\right\|\leq \sqrt 2 C_*\left(1+a^{-1}\lambda_0 K\right)\left[(1+a^{-1})K\right]^j\|h\| .
\]
The case $j=0$ shows that $\cL_{k,\lambda,n}$ is an operator in $L(\bB_k,\bB_{k+1})$, while the cases $j\geq 1$ prove the inequality in (O-\ref{cond:multiplier}) with $j_0=\infty$ and $K$ replaced by $\sqrt 2 C_*\left(1+a^{-1}\lambda_0 K\right)(1+a^{-1})K$.

Taking the derivative as in \eqref{eq:LY_exp} we have 
\[
\|\nabla \cL_{k,\lambda}h\|\leq \vartheta^{-1}\cL_{k,0} \|\nabla h\|+(\bar D+|\lambda| \sigma_n^{-1}  \|\nabla g_{k}\|_\infty) \cL_{k,0}|h|.
\]
Which can be iterated, yielding
\[
\begin{split}
\|\nabla \cL_{k,\lambda}\cdots \cL_{ j, \lambda} h\|
&\leq \vartheta^{-k+j+1}\cL_{k}\cdots \cL_j \|\nabla h\|+\sum_{l=0}^{k-j}\vartheta^{-l}(\bar D+|\lambda|  \sigma_n^{-1}  \|\nabla g_{k-l}\|_\infty )  \cL_{k-l}\cdots \cL_j |h|\\
&\leq \|h \| \frac{3 e^a(\bar D+|\lambda|  \sigma_n^{-1}   \sup_j\|\nabla g_{j}\|_\infty )}{1-\vartheta^{-1}}.
\end{split}
\]
where we have used \eqref{exp_basic} in the last step. Applying \eqref{exp_basic} again
proves (O-\ref{cond:lambda_multiplier}).

Accordingly, Theorem \ref{thm:main} applies. 
Note that Theorem \ref{thm:main} provides nontrivial information only if $\sigma_n\geq \Const n^\alpha$ for some $\alpha>\frac 13$.
\subsection{Complex Cones for Expanding Maps}\ \\
To obtain the sharper results of Theorem~\ref{thm:main_bis}, we have to verify condition (O-\ref{cond:lambda_multiplier2}). In order to do that, we use the complex cone theory introduced by Rugh and further developed by  Dubois \cite{Ru10, Du09, Du11}, 
as summarized in Section~\ref{sec:complex}.

First we need to verify that we are in the setting of Section~\ref{sec:complex}.
As noted earlier, with $\bbe =1$, the norm associated the cone is $\|h\|=\sup_x\{|h(x)|+a^{-1}\|\nabla h(x)\|\}$ which is equivalent to the $C^1$ norm, so $\cB_\bR= C^1(M, \bR)$.

For each $x\in M$ and $v\in\bR^d$, {$\| v\|  \le 1$, }  we define 
\[
\ell_{x,v}(h)=ah(x)-\langle v, \nabla h(x) \rangle .
\]
Clearly $\ell_{x,v}\in\cB_\bR'$.
In addition, it is easy to check that 
\[
\begin{split}
&\cC_a=\{h\in\cB_\bR\;:\; \ell (h)\geq 0\; \forall \ell\in\cS\}\\
&\cS=\{\ell_{x,v}\;:\; x\in M, \|v\|  \le  1\}
\end{split}
\]

This shows that \eqref{eq:com_cone_def0} and \eqref{eq:com_cone_def} are satisfied.
The complex cone $\cC_{\bC}$ is defined precisely as in Definition~\ref{def:complex cone}.

Next, set $\bbm(h)=\int h$, then, recalling \eqref{eq:average bound}, we have, for all $h\in\cC_a$,
\begin{equation}\label{eq:maybe_exp}
\bbm(h)\geq e^{-a} \|h\|_\infty\geq \frac 12e^{-a}\|h\|:=\kappa \|h\| \, ,
\end{equation}
which verifies \eqref{eq:maybetoomuch}, and by extension \eqref{eq:m_lower}.

To conclude checking the hypotheses of Theorem \ref{thm:complex_contraction} we must to prove that the diameter of the image of the complex cone is finite. Since \eqref{eq:Delta_exp} states that the diameter of the real cone is finite, we can apply Theorem \ref{thm:needed}, provided we check \eqref{eq:complex_bound}.
Recall that $\nabla \hat g_k=\nabla g_k$ and $|\hat g_k|_\infty\leq 2|g_k|_\infty$, hence $\|\hat g_k\|_{\cC^1}\leq 2\|g_k\|_{\cC^1}$.

To check  \eqref{eq:complex_bound} let  $\ell_{x,v}\in \cS$, $h \in \cC_a$,
\begin{equation}
\label{eq:upper ell}
\begin{split}
\left|\ell_{x,v}\left((\cL_{k,\lambda}-\cL_k) h\right)\right|&\leq \left|(\cL_{k,\lambda}-\cL_k) \langle (D f)^{-1}v,\nabla h\rangle\right|+ \left|(\cL_{k,\lambda}-\cL_k) \cD_v h\right|\\
&+\left|\cL_{k,\lambda}(i\lambda \sigma_n^{-1} \langle (D f)^{-1}v,\nabla \hat g_k\rangle h)\right| 
+ a \left| (\cL_{k, \lambda} - \cL_k)h \right| , 
\end{split}
\end{equation}
where $\cD_v=-\sign (\det Df_k)\frac{\langle Df_k^{-1}v,\nabla \det Df_k\rangle}{(\det Df_k)^2}$.

Since $\cL_k 1\in\cC_a$, \eqref{eq:average bound} implies the following bound on the first term above, 
\[
\begin{split}
\vartheta^{-1} \left|(\cL_{k,\lambda}-\cL_k) \| \nabla h \| \right|
& \leq \vartheta^{-1} |\lambda| \sigma_n^{-1}  | \hat g_k |_\infty |\cL_k 1 |_\infty \| \nabla h \|_\infty
\leq \vartheta^{-1} a e^a |\lambda| \sigma_n^{-1} |\hat g_k|_\infty | h |_\infty \\
& \le \vartheta^{-1} a e^{2a} |\lambda| \sigma_n^{-1} |\hat g_k|_\infty \int h 
\le \vartheta^{-1} a e^{2a} |\lambda| \sigma_n^{-1} |\hat g_k|_\infty \int \cL_k h \\
& \le \vartheta^{-1} a e^{3a} |\lambda| \sigma_n^{-1} |\hat g_k|_\infty \cL_k h(x) \, ,
\end{split}
\]
where we have used $h, \cL_k h \in \cC_a$.  Similarly, we bound the second, third and fourth terms
of \eqref{eq:upper ell} using the positivity of $h$,
\[
|\lambda| \sigma_n^{-1} (\bar D + a) |\hat g_k|_\infty \cL_k h(x) + \vartheta^{-1} | \lambda| \sigma_n^{-1} \| \nabla \hat g_k \|_\infty \cL_kh(x) \, .
\]
Putting these estimates together in \eqref{eq:upper ell} implies,
\[
\left|\ell_{x,v}\left((\cL_{k,\lambda}-\cL_k) h\right)\right|
\leq |\lambda| \sigma_n^{-1} \| \hat g_k\|_{\cC^1} \left[ \vartheta^{-1} ae^{3a} + \bar D +a \right]\cL_kh (x). 
\]
On the other hand, $\cL_k h\in\cC_{\nu a}$, hence 
\[
\ell_{x,v}(\cL_k h)=a\cL_k h-\langle v, \nabla \cL_k h\rangle\geq a(1-\nu)\cL_k h(x).
\]
Accordingly,
\[
\left|\ell_{x,v}\left((\cL_{k,\lambda}-\cL_k) h\right)\right|\leq 
|\lambda| \sigma_n^{-1} a^{-1}(1-\nu)^{-1} \| \hat g_k\|_{\cC^1}\left[ \vartheta^{-1} a e^{3a}+ \bar D +a \right]\ell_{x,v}(\cL_k h).
\]
which satisfies  \eqref{eq:complex_bound} provided $\lambda$ is small enough.  Hence Theorem~\ref{thm:needed} proves the finite diameter of the complex cone,\footnote{Note that since $\bbe = 1 \in \cC_{\nu a}$,
we may take $B=0$ in the bound on the diameter of the complex cone in Theorem~\ref{thm:needed}.} and Theorem~\ref{thm:complex_contraction} applies.

\subsection{Verification of Condition (O-\ref{cond:lambda_multiplier2})}\ \\
\label{sec:checkO3}
We can finally check condition (O-\ref{cond:lambda_multiplier2}).
By the above discussion it follows that, for $\lambda_0$ small enough, and $|\lambda|\leq \lambda_0\sigma_n$, we have that Theorem \ref{thm:complex_contraction}  applies to each $\cL_{k,\lambda, n}$ with complex diameter $\Delta_{\bC}< \infty$ uniform in $\lambda$, $k$ and $n$. 
Moreover, since $\cC_{\bR, k} = \cC_a$ and $\bbm_k = \bbm = \mu$, the Riemannian volume,
for each $k$, we are in the setting of Section~\ref{sec:memory loss} with $\bar\kappa = \kappa = \frac 12 e^{-a}$, $K$ the constant from (O-\ref{cond:lambda_multiplier}) and $\bL_k = \cL_{k, \lambda, n}$.
We can then apply Lemma~\ref{lem:memory loss} to obtain (O-\ref{cond:lambda_multiplier2}) with
$\alpha_{k,j,\lambda} = \boldsymbol{\alpha}_{k,j}$, $h_{k,j,\lambda} = \boldsymbol{h}_{k,j}$
and $\ell_{k,j,\lambda} = \boldsymbol{\ell}_{k,j}$, as defined in Section~\ref{sec:memory loss}.

By definition $\ell_{k,j,\lambda}(1)=1$ and $\bbm(h_{k,j,\lambda})=1$.
Note that this normalization is compatible with that required in (O-\ref{cond:lambda_multiplier2})
since $\bbm(h_{k,j,\lambda}) = \int h_{k,j,\lambda}$. 
Moreover, $|\alpha_{k,j,\lambda}| \le 1$ and $\| h_{k,j,\lambda} \| \le \sqrt{2}/\kappa$,
$\| \ell_{k,j,\lambda} \|' \le \sqrt{2}$ by Lemma \ref{lem:boundedness}. 
Note that all these bounds are independent of $n$.

Then since condition (b) of Section~\ref{sec:memory loss} holds with $K=1$, 
Lemma~\ref{lem:memory loss}
implies that  the property \eqref{eq:O3 decay} of (O-\ref{cond:lambda_multiplier2}) holds with $K = \frac{ 2(1+\sqrt{2}) }{\kappa \theta^2} \Delta_{\bC}$,
$\theta = \tanh(\Delta_{\bC}/4)$
and $c =  - \log \theta$.

Since we have already verified the hypotheses of Theorem~\ref{thm:needed}, conditions (d) and (f) of Section~\ref{sec:lower bound} are satisfied. Moreover, condition (e) of Section~\ref{sec:lower bound} is satisfied by \eqref{eq:Delta_exp}.
Thus the final claim that  there exists $K >0$ such that $|\ell_{k,j,\lambda}(h_{j,l,\lambda})| \ge K^{-1}$ for all $0 \le l \le j \le k \le n$, follows immediately from Lemma~\ref{lem:lower bound}.


\subsection{Central Limit Theorem and Variance for Sequential Expanding Maps}\ \\
\label{sec:CLT expanding}
We have verified hypotheses (C-\ref{cond:positivity})-(C-\ref{cond:powerbounded}) and 
(O-\ref{cond:multiplier})-(O-\ref{cond:lambda_multiplier2}). Thus Theorem~\ref{thm:main_bis}
and Corollary~\ref{cor:better corollary} apply to the setting of sequential expanding maps. Hence, Theorem \ref{thm:expanding CLT} follows.

 As mentioned in Section~\ref{sec:exp preview}, Theorem \ref{thm:expanding CLT} is conditional on the growth of $\sigma_n$. 
To find conditions implying that the variance grows sufficiently fast, in the abstract sequential case setting, is nonobvious.
Here, we provide a partial answer for the case at hand.

To simplify matters, we assume that the maps and observables appearing in the sequence all belong to a finite set\footnote{ The finiteness allows us to obtain a condition that can be checked on a finite number of finite time trajectories; otherwise some type of compactness seems to be needed.}
\[
\begin{split}
&\cF:=\{\bar f_j\}_{j=1}^N\;\quad \Lambda =\inf_j\inf_x \|(D_x\bar f_j)^{-1}\|^{-1}>1\\
&\cG=\{\bar g_j\}_{j=1}^{N}\subset \cC^1(M,\bR).
\end{split}
\]
For any given sequences $\omega^f,\omega^g\in \{1,\dots, N\}^{\bN\cup \{0\}}$, we set $f_k=\bar f_{\omega^f_k}$, $g_k=\bar g_{\omega^g_k}$.\\
The maps $\cF$ satisfy specification: for each $\ve, L>0$, there exists $D_{\cF}(\ve,L)$ such that for each set of intervals $J_k=\{a_k,\dots, a_k+L\}$, such that $a_{k+1}-a_k-L\geq D_{\cF}(\ve,L)$, and points $x_k\in M$, there exists $z\in M$ such that 
\[
| f_{a_k+i}\circ\cdots\circ f_{0}(z)- f_{a_k+i}\circ\cdots\circ  f_{a_k}(x_k)|\leq \ve,
\]
for all $k\in\bN$ and $i\in\{0,\dots,L-1\}$.

The main result in this section is the following criterion. See Lemma \ref{lem:specif} for a quantitative version.
\begin{prop}\label{prop:variance}
Let $\ln \rho \in \cC^1(M)$. If there exists $a\in (0,1)$ and $L\in\bN$ such that for each sequence $\overline\omega=\{(\omega^1_i, \omega^2_i)\}_{i=0}^{L}$  there exists a point $x_{\overline\omega}$ such that~\footnote{ For the definition of $\hat g_{\omega_i}$ see equation \eqref{eq:centering1}.}
\[
\begin{split}
&\sum_{i=1}^{L} \hat g_{\omega^2_i}(\bar f_{\omega^1_{i-1}}\circ \cdots\circ\bar f_{\omega^1_0}(x_{\overline\omega}))\geq a L\\
&\mbox{where} \quad \ve\leq \frac{a}{4\sup_j\|D \bar g_j\|_{\cC^0}} \quad \mbox{and} \quad L\geq \frac {8} { a} D_{\cF}(\ve,L)\max_j\|\hat g_j\|_\infty.
\end{split}
\]
then there exists $B>0$ such that, for $n$ large enough, $\sigma_n^2=\bE( \hS_n^2)\geq Bn$.
\end{prop}
Note that, strictly speaking, the quantity $\hat g_{\omega^2_i}$ is not determined by the finite word $\overline\omega$, since by \eqref{eq:centering1} the centering constant depends on the whole past through $\cL_{k-1}\cdots\cL_0\rho$. This is remedied in Corollary~\ref{cor:verify}, where $\hat g_k$ is replaced by $\tilde g_k=g_k-\gamma_{S,k}$, which depends only on the last $S$ steps; it is that form of the criterion which is effectively checkable.
The rest of the section is devoted to the proof of the above Proposition. But first, it is appropriate to discuss its relevance and generality.

\begin{rem}
For simplicity, in Proposition~\ref{prop:variance}, we assume  that $\ln \rho \in \cC^1$.  This is done so that
$\rho$ belongs to a real cone $\cC_a$, as defined in Section~\ref{sec:expanding intro}, provided $a$ is chosen large enough.  
The interested reader can generalize to the case $\rho \in C^1$, $\rho \ge 0$, since then there exists a $k_0 \in \bN$ such that $\cL_k \cdots \cL_0 \rho \in \cC_a$ for all $k \ge k_0$.
\end{rem}

\begin{rem}\label{rem:comparison}
Note that in the case of a single dynamical system and a single observable (that is $\bar f_k=f$, and $g_k=g$, a zero average observable), the usual condition for the linear growth of the variance is that there exists a periodic orbit $\{p, f(p), \dots, f^{q-1}(p)\}$, $f^q(p)=p$, such that\footnote{Possibly substituting $-g$ for $g$.} $\sum_{k=0}^{q-1} g( f^k(p))\geq b$, for some $b>0$.\footnote{Which implies that $g$ is not a coboundary.} This implies the condition in Proposition \ref{prop:variance}. To see it choose $\mu_k$ to be the invariant measure of $f$, so $\hat g=g$. Then, setting $L=Kq$, $K\in\bN$,
\[
\sum_{i=1}^{L} g( f^{i-1}(p))= K\sum_{i=1}^{q} g( f^{i-1}(p))\geq Kb= L\frac{b}{q}=: aL.
\]
Accordingly, the conditions of Proposition \ref{prop:variance} are satisfied with the choice $x_\omega =p$. This shows that  Proposition \ref{prop:variance} is a natural generalization to the sequential case of the usual condition. It remains to see if the condition is explicitly checkable;  this is verified in Corollary \ref{cor:verify}, although it is unclear if our condition can be reduced to checking that the observables are not coboundaries w.r.t. some dynamics.
\end{rem}
 
Note that, by equation \eqref{eq:centering1},  $\|\hat g_k\|_\infty\leq 2\|g_k\|_\infty$.
Moreover, the $\cL_k$ have a uniform spectral gap. Hence, there exists a computable $A>0$, $\nu\in (0,1)$ such that, for all $S\in\bN$ and $k\geq S$,

\begin{equation}\label{eq:nu_decay}
\begin{split}
&\left|\int_{M} g_k \cL_{k-1}\cdots \cL_0 \rho-\gamma_{S,k}\right|\leq A\nu^S\\
&\gamma_{S,k}:=\frac{1}{\int_M 1}\int_{M} g_k \cL_{k-1}\cdots \cL_{\max\{k-S,0\}} 1.
\end{split}
\end{equation}

\begin{cor}\label{cor:verify} There exist computable constants $\kappa, b, L_0>0$ such that, if for some $L\geq L_0$ and each sequence $\overline\omega=\{\omega^1_i,\omega^2_i\}_{i=0}^{L}$, there exists a point $x_{\overline\omega}$ such that 
\begin{equation}\label{eq:checkable}
\sum_{i=1}^{L} \tilde g_{\omega^2_i}(\bar f_{\omega^1_{i-1}}\circ \cdots\circ\bar f_{\omega^1_0}(x_{\overline\omega})) \geq 2\kappa \ln L,
\end{equation}
where $\tilde g_k=g_k-\gamma_{S,k}$ and $S=b\ln L$, then the hypothesis of Proposition \ref{prop:variance} are verified. Thus, there exists $B>0$ such that $\sigma_n^2\geq Bn$.
\end{cor} 
\begin{proof}
Let $a=\kappa \frac{\ln L}L$. We use \eqref{eq:nu_decay}: if  $b > -1/\ln \nu$,  then $A\nu^S\leq a$ for all $L$ larger than some computable  $L_0$, and
\[
\sum_{i=1}^{L} \hat g_{\omega^2_i}(\bar f_{\omega^1_{i-1}}\circ \cdots\circ\bar f_{\omega^1_0}(x_{\overline\omega}))\geq a L.
\]
Set $\Gamma=\sup_j\max\{\|D g_j\|_\infty, 2\|g_j\|_\infty\}$ and choose $\ve=\frac{a}{4\Gamma}$. 
Note that \eqref{eq:nu_decay} requires $k\geq S$, while for the finitely many indices $k<S$ we have only $|\hat g_k-\tilde g_k|\leq 2\Gamma$; the hypothesis of Proposition \ref{prop:variance} is thus verified for all the words that do not contain indices smaller than $S$. This suffices, since the block containing such indices then contributes an extra $\cO(S\Gamma)=\cO(1)$ to the lower bound for $S_n(p)$ in the proof of Proposition \ref{prop:variance}, which is irrelevant for $n$ large.  Then Proposition \ref{prop:variance} applies if
\begin{equation}\label{eq:LGamma}
L\geq \frac {8 D_{\cF}(\ve,L)\Gamma}a.
\end{equation}
By Lemma \ref{lem:specif}, there exists a computable $c_\star>0$ such that  we can choose $D_{\cF}(\ve,L)=c_\star\ln a^{-1}$. Consequently, equation \eqref{eq:LGamma} is satisfied if, for $L\geq L_0$,
\[
1\geq  \frac{8c_\star\ln a^{-1}\Gamma}{a L}=8 c_\star\Gamma \frac 1{\kappa \ln L}\ln\left[\frac L{\kappa\ln L}\right]=\frac{8 c_\star\Gamma}{\kappa} \frac {\ln L-\ln\ln L-\ln\kappa}{\ln L}
\]
which is satisfied if $\kappa= 8c_\star\Gamma$. 
\end{proof}

Note that condition \eqref{eq:checkable} is checkable by trial and error. This is similar to the problem of finding a periodic orbit with non-zero average in the nonsequential case (see Remark \ref{rem:comparison}).

To prove Proposition \ref{prop:variance}, we will use a martingale decomposition, 
similar to \cite{DH24}, but we must first introduce some notation.
Let
\begin{equation}\label{eq:hdef}
h_k=\cL_{k-1}\cdots\cL_0 \rho
\end{equation}
and $\widehat \cL_k(g)=h_{k+1}^{-1}\cL_k (h_k g)$. Note that $\widehat \cL_k\cdots\widehat\cL_i(g)=h_{k+1}^{-1}\cL_k\cdots\cL_i(h_i g)$, moreover $\widehat\cL_k 1=1$. Also, recalling \eqref{eq:centering0},
\begin{equation}\label{eq:centering1}
\hat g_k=g_k-\int_{M} g_k \cL_{k-1}\cdots \cL_0 \rho 
\end{equation}
We would like to find $\phi_k, Y_k\in \cC^1$ such that
\begin{equation}\label{eq:coboundary0}
\begin{split}
&\hat g_k=\phi_{k+1}\circ f_k-\phi_k+Y_k\\
&\phi_0=0\\
&\widehat \cL_k Y_k=0.
\end{split}
\end{equation}
This decomposition would imply,
\[
\wcL_k \hat g_k=\phi_{k+1}-\wcL_k\phi_k \, .
\]
With this motivation, let us define $\phi_0 =0$ and for each $k \ge 1$,
\begin{equation}\label{eq:coboundary}
\phi_k=\sum_{j=0}^{k-1}\wcL_{k-1}\cdots \wcL_j \hat g_j,
\end{equation}
and set $Y_k = \hat g_k - \phi_{k+1} \circ f_k + \phi_k$.  Then indeed,
\[
\begin{split}
-  \wcL_k Y_k=\wcL_k\left[\phi_{k+1}\circ f_k-\phi_k-\hat g_k\right]&=\sum_{j=0}^{k}\left[\wcL_{k}\cdots\wcL_j \hat g_j \right]-\sum_{j=0}^{k-1}\wcL_{k}\cdots\wcL_j \hat g_j -\wcL_k \hat g_k=0.
\end{split}
\]
Note that, by usual arguments, there exists $C_\star>0$ such that
\begin{equation}\label{eq:standard_bound}
\begin{split}
&\|\phi_k\|_{\cC^1}\leq C_\star\,,\;\;\|Y_k\|_{\cC^1}\leq C_\star\\
& C_\star^{-1} \leq h_k \leq C_\star.
\end{split}
\end{equation}
Let $\hS_n=\sum_{k=0}^{n-1}\hat g_k\circ f_{k-1}\circ \cdots\circ f_0$. We want to compute\footnote{ As usual, by $\bE(\vf)$ we mean $\int_M\vf(x) \, \rho \, dx$.}

\[
\bE(\hS_n^2)=\bE\left(\left[\phi_{n+1}\circ  f_{n}\circ \cdots\circ f_0+\sum_{k=0}^{n-1}Y_k\circ f_{k-1}\circ \cdots\circ f_0\right]^2\right)
\]
Note that, recalling \eqref{eq:hdef} and also \eqref{eq:coboundary0}, for all $k<n$,
\[
\begin{split}
\bE&  \left( \phi_{n+1} \circ f_n \circ \cdots \circ f_0 \cdot Y_k \circ f_{k-1} \circ \cdots \circ f_0 \right)
 = \int \phi_{n+1} \circ f_n \circ \cdots \circ f_k \cdot Y_k  \cdot \cL_{k-1} \cdots \cL_0 \rho \, dx \\
& = \int \phi_{n+1} \circ f_n \circ \cdots \circ f_k \cdot Y_k h_k \, dx \\
& = \int \phi_{n+1} \circ f_n \circ \cdots \circ f_{k+1} \cdot h_{k+1} \cdot \wcL_k (Y_k) \, dx = 0
\end{split}
\]
Accordingly, recalling \eqref{eq:standard_bound} and \eqref{eq:coboundary0},
\begin{equation}\label{eq:martingale}
\begin{split}
\bE(\hS_n^2)&=\bE\left(\left[\phi_{n+1}\circ  f_{n}\circ \cdots\circ f_0+\sum_{k=0}^{n-1}Y_k\circ f_{k-1}\circ \cdots\circ f_0\right]^2\right)\\
&=\bE\left(\left[\sum_{k=0}^{n-1}Y_k\circ f_{k-1}\circ \cdots\circ f_0\right]^2\right)+\cO(1)\\
&=\sum_{k=0}^{n-1} \int_M h_{k}Y_k^2+2 \sum_{k>i}\bE\left(h_i Y_k\circ f_{k-1}\circ \cdots\circ f_{i} Y_i\right)+\cO(1)\\
&=\sum_{k=0}^{n-1} \int_M h_{k}Y_k^2+2 \sum_{k>i}\bE\left(h_{k} Y_k\wcL_{k-1}\cdots\wcL_i Y_i\right)+\cO(1)\\
&=\sum_{k=0}^{n-1} \int_M h_{k}Y_k^2+\cO(1).
\end{split}
\end{equation}
Next, we want to connect the growth of $\hS_n$ along a single orbit with the growth of its $L^2$ norm.

\begin{lem}\label{lem:l2-bound} There exists $B>0$ such that if, for some $p\in M$,  $|\hS_n(p)|\geq 3An^\alpha$, $\alpha> \frac{d+1}{d+2}$, and $n$ is large enough, depending on $A$, $\alpha$ and $C_\star$, then we have $\bE(\hS_n^2)\geq A^{d+2}B n^{(d+2)\alpha-d-1}$.
\end{lem}
\begin{proof}
Equation \eqref{eq:coboundary0} allows us to write
\[
\begin{split}
|\hS_n(p)|&=\left|\phi_{n+1}\circ  f_{n}\circ \cdots\circ f_0(p)+\sum_{k=0}^{n-1}Y_k\circ f_{k-1}\circ \cdots\circ f_0(p)\right|\\
&\leq \sum_{k=0}^{n-1}|Y_k|\circ f_{k-1}\circ \cdots\circ f_0(p)+\cO(1).
\end{split}
\]
Note that, for each $q\in M$, $|Y_k(x)|\geq |Y_k(q)|-C_\star |x-q|$. If follows that, for all $x\in \{y\in M\;:\; |y-q|\leq |Y_k(q)|/(2C_\star)\}$, one has $|Y_k(x)| \ge \frac 12 |Y_k(q)|$. 
Thus,  recalling \eqref{eq:standard_bound},
\[
 \int_M h_k Y_k^2\,dx \geq  C_\star^{-1}\int_M Y_k^2\,dx \geq  c_d|Y_k|^{d+2}\circ f_{k-1}\circ\cdots\circ f_0(p),
\]
for some constant $c_d$.
Let $J_n=\{ k \in\{0,\dots,n-1\}\;:\; |Y_k|\circ f_{k-1}\circ\cdots\circ f_0(p)\geq An^{-1+\alpha}\}$, then
using \eqref{eq:martingale},
\[
\begin{split}
3An^{\alpha}&\leq |\hS_n(p)|\leq \frac{  n^{(d+1)(1-\alpha)}}{A^{d+1}}\sum_{k\in J_n}|Y_k|^{d+2}\circ f_{k-1}\circ \cdots\circ f_0(p)+An^\alpha+\cO(1)\\
&\leq \frac{ n^{(d+1)(1-\alpha)}}{A^{d+1}c_d}\sum_{k=0}^{n-1} \int_M h_k Y_k^2+An^\alpha +\cO(1)\leq   \frac{ n^{(d+1)(1-\alpha)}}{A^{d+1}c_d}\bE(\hS_n^2)+An^\alpha+\cO(n^{(d+1)(1-\alpha)}+1).
\end{split}
\]
For $n$ large enough we have $An^\alpha+\cO(n^{(d+1)(1-\alpha)}+1)\leq 2An^\alpha$, consequently
\[
\bE(\hS_n^2)\geq  \frac{A^{d+1}c_d}{ n^{(d+1)(1-\alpha)}}An^\alpha,
\]
from which the lemma follows.
\end{proof}
 
\begin{proof}[\bfseries Proof of Proposition \ref{prop:variance}]
Let $D:=D_{\cF}(\ve,L)$. For any $\omega^f,\omega^g\in\{1,\dots N\}^{\bN\cap\{0\}}$, let $f_k=\bar f_{\omega^f_k}$ and $g_k =\bar g_{\omega^g_k}$ and $r_k=k(L+D)$.  We assume $n\geq \frac{4(L+D)^2}{D}$. By hypothesis there exists $p\in M$ such that
\[
\begin{split}
\hS_n(p)&\geq \sum_{k=0}^{ \lfloor \frac{n}{L+D}\rfloor-1} \sum_{i=0}^{L-1}\hat g_{r_k+i+1}(f_{r_k+i}\circ \cdots\circ f_{0}(p))-\left[\frac{n D}{L+D}+ L\right] \sup_j\| \hat g_j\|_\infty\\
&\geq \sum_{k=0}^{ \lfloor \frac{n}{L+D}\rfloor-1} \sum_{i=0}^{L-1}\hat g_{r_k+i+1}( f_{r_{k}+i}\circ \cdots\circ f_{r_k}(x_{k}))-\frac{\ve Ln}{L+D} \sup_j\|D \hat g_j\|_{\cC^0} - \frac{2n D}{L+D} \sup_j\| \hat g_j\|_\infty \\
&\geq  aL\left\lfloor \frac{n}{L+D}\right\rfloor-\frac{a Ln}{4(L+D)}-\frac{anL}{4(L+D)}\geq  \frac{aLn}{2(L+D)}-aL\geq \frac{aLn}{4(L+D)}.
\end{split}
\]
We can then apply Lemma \ref{lem:l2-bound} with $\alpha=1$ and $3A=\frac{aL}{4(L+D)}$, which proves the proposition for $B$ small enough.
\end{proof}

It remains to discuss the Lemma used to obtain Corollary \ref{cor:verify}.
\begin{lem}\label{lem:specif}  There exists $\ve_0\in(0,1)$ such that, for all $\ve\leq \ve_0$, we can choose $D_{\cF}(\ve,L)=\inf\{n\in\bN\;:\; \ve \Lambda^n\geq \operatorname{diam}(M)\}$.
\end{lem}
\begin{proof}
Set $\ve_0$ so that all the maps $\bar f_j$ are locally invertible on balls of size $\ve_0$.
For $\ve \le \ve_0$, consider the set  
\[
A_k=\{y\in M\;:\; | f_{a_k+i}\circ\cdots\circ f_{a_k}(y)- f_{a_k+i}\circ\cdots\circ  f_{a_k}(x_k)|\leq \ve,\, 0\leq i\leq L-1\}.
\]
Then $ f_{a_k+L-1}\circ\cdots\circ f_{a_k}(A_k)$ is a ball of radius $\ve$. Accordingly, 
\[ 
f_{a_k+  L-1+D_{\cF}(\ve,L)}\circ\cdots\circ f_{a_k}(A_k)= M.
\]
Thus $f_{a_k+ L-1+D_{\cF}(\ve,L)}\circ\cdots\circ f_{a_k}(A_k)\supset A_{k+1}$. Taking the intersection of the preimages of the $A_k$, we have a nested sequence of closed sets and we can choose as $z$ any point in the intersection.
\end{proof}


\section{Application: Sequential Billiard Maps}
\label{sec:billiard}

In this section, we show that Theorems~\ref{thm:main} and \ref{thm:main_bis} apply to sequential billiards 
as described in \cite{DL22}.  Our strategy will be to define real and complex cones on which
the relevant operators act as strict contractions, and so apply the theorems of Section~\ref{sec:complex} to
verify conditions (O-\ref{cond:multiplier})-(O-\ref{cond:lambda_multiplier2}).
We begin by recalling the setting of \cite{DL22}.

Fixing constants $K \in \bN$ and $\ell_i >0$, $i = 1, \ldots K$,  we consider billiard configurations in 
$\mathbb{T}^2 = \mathbb{R}^2/ \mathbb{Z}^2$
formed by placing $K$ scatterers $B_i$ having $C^3$ boundaries and with arclengths given by $\ell_i$
and curvatures $\cK >0$.  We denote such a billiard table by $Q = \mathbb{T}^2 \setminus \cup_i B_i$.  The billiard map $T$
is defined by the motion of a point particle undergoing elastic collisions at the boundaries
and traveling with unit speed in straight lines between collisions.  We adopt the standard collision
coordinates, $x= (r,\vf)$, where $r$ is the arclength coordinate on $\partial Q$ and $\vf$ is the angle
made by the post-collision velocity with the normal to the boundary.
Thus the phase space for $T$ is given by $\cM = \cup_{i=1}^K I_i \times [- \frac \pi 2, \frac \pi 2]$, where for each $i$, $I_i = [0, \ell_i] / \sim$ is an interval
of length $\ell_i$ with endpoints identified.  Since we have fixed $K$ and $\ell_i$, $\cM$ is the same
for each such table $Q$.

Let $\cK_{\min}(Q) \le \cK_{\max}(Q)$ denote the minimum and maximum curvatures of scatterers
in the table $Q$, and let $\tau_{\min}(Q) \le \tau_{\max}(Q)$ denote the minimum and maximum
distances between consecutive collisions.

These scatterer configurations are subject to the following three constraints:
Fix $\tau_*, \cK_*, E_*>0$.  Then for each configuration $Q$, we assume:
$\tau_* \le \tau_{\min}(Q) \le \tau_{\max}(Q) \le \tau_*^{-1}$, 
$\cK_* \le \cK_{\min}(Q) \le \cK_{\max}(Q) \le \cK_*^{-1}$,
the $C^3$ norm of $\partial Q$ is at most $E_*$.
Let $\cQ(\tau_*, \cK_*, E_*)$ denote the set of billiard tables with $K$ scatterers satisfying these conditions
and let   
$\cF(\tau_*, \cK_*, E_*)$ denote the corresponding set of billiard maps.
As noted above, each $T \in \cF(\tau_*, \cK_*, E_*)$ acts on the same phase space $\cM$.

Fixing these constants ensures that maps in $\cF(\tau_*, \cK_*, E_*)$ enjoy the uniform
properties listed as {\bf (H1)}-{\bf (H5)} in \cite{DL22}: a common set of invariant cones in 
the tangent space, uniform distortion bounds, uniform growth lemma, etc.

We choose as our reference measure on $\cM$, $\musrb := \frac{1}{2 |\partial Q|} \cos \vf \, dr d\vf$, 
which is preserved by each $T \in \cF(\tau_*, \cK_*, E_*)$.  The associated transfer operator
$\cL_T$ acting on measurable functions is defined by,
\[
\cL_T h = h \circ T^{-1} \, .
\]
In order to generate sufficient hyperbolicity to strictly contract the cones $\cC_{c,A,L}(\delta)$ defined 
in Section~\ref{sec:cone def} below,
we require a notion of distance in $\cQ(\tau_*, \cK_*, E_*)$.  For a table $Q \in \cQ(\tau_*, \cK_*, E_*)$,
the boundary $\partial B_i$ of each scatterer can be parametrized according to arclength by a function
$u_{i, \theta}$, where $\theta \in [0, \ell_i)$ represents the tagged point in $\partial B_i$ 
where the parametrization begins.  Given two tables $Q, \tilde Q \in \cQ(\tau_*, \cK_*, E_*)$ with
respective parametrizations $u_{i,\theta}$, $\tilde u_{i, \theta}$, let $\Pi_K$ denote the set of permutations
$\pi$ on $\{ 1, \ldots, K \}$ satisfying $\tilde \ell_{\pi(i)} = \ell_i$.  Define the distance,
\[
\mathbbm{d}(Q, \tilde Q) = \min_{\pi \in \Pi_K} \min_{\theta \in [0, \ell_i)} \sum_{i=1}^K |u_{i,\theta} - \tilde u_{\pi(i), \theta}|_{C^2(I_i, \mathbb{R}^2)} .
\]
For $Q_0 \in \cQ(\tau_*, \cK_*, E_*)$ and $\ve \le \frac 12 \min \{ \tau_*, \cK_* \}$, define 
\[
\cQ(Q_0, E_*; \ve) = \{ Q \in \cQ(\tau_*/2, \cK_*/2, E_*) : \mathbbm{d}(Q, Q_0) < \ve \} .
\]
Let $\cF(Q_0, E_*; \ve)$ denote the corresponding set of billiard maps.

\begin{defin}
\label{def:admissible}
Fix $Q_k \in \cQ(\tau_*, \cK_*, E_*)$, $\ve>0$ from \cite[Lemma~6.6(b)]{DL22} and $N_{\cF} \in \mathbb{N}$
from \cite[Theorem~2.3]{DL22}.
We define each of our maps $f_k$ (in the notation of Section~\ref{sec:setting}) by
$f_k = T_{N_\cF} \circ \cdots \circ T_1$, where $T_j \in \cF(Q_k, E_*; \ve)$.

With this definition of $f_k$, the sequence $f_{k,0} = f_{k-1} \circ \cdots \circ f_0$ is necessarily an
{\em $N_{\cF}$-admissible sequence} in the terminology of \cite[Definition~2.4]{DL22}.  
To such sequences, we may the apply the uniform properties {\bf (H1)}--{\bf (H5)} of \cite{DL22} corresponding to
the family with chosen constants $\cF(\tau_*/2, \cK_*/2, E_*)$.
\end{defin}

Note that each
$f_k$ is comprised of a collection of $N_\cF$ billiard maps all of which are close with respect to the
distance $\mathbbm{d}$, but the maps comprising $f_{k+1}$ are not necessarily close to the maps
comprising $f_k$.\footnote{The notion of $N_{\cF}$-admissible sequence can also be simplified by choosing $f_k = T_k^{N_{\cF}}$,
where $T_k \in \cF(\tau_*, \cK_*, E_*)$, i.e. simply requiring $N_{\cF}$ iterates of the same map before changing
to a different (but not necessarily close) map.}

In the notation of Section~\ref{sec:setting} then, we have $M_k = \cM$, $\mu_k = \musrb$, while 
$f_k = T_{N_\cF} \circ \cdots \circ T_1$, and $\cL_k = \cL_{T_{N_\cF}} \cdots \cL_{T_1}$, where $T_j \in \cF(Q_k, E_*; \ve)$
as in Definition~\ref{def:admissible}.

\subsection{Definition and Contraction of the Real Cone $C_{\bR} = \cC_{c,A,L}(\delta)$}\ \\
\label{sec:cone def}
Next, we must define precisely the cone of functions
constructed in \cite{DL22} which is strictly contracted by $\cL_k$.

According to \cite[Section~3.1]{DL22}, there exist stable and unstable cones $C^s$ and $C^u$ 
in the tangent space of $\cM$
(not to be confused with the projective cones of functions defined below), that are strictly contracted by
$DT^{-1}$ and $DT$, respectively, for all $T \in \cF(\tau_*, \cK_*, E_*)$.  
For fixed $k_0 \in \mathbb{N}$, define the usual homogeneity strips,
\[
\mathbb{H}_{\pm k} = \{ (r,\vf) \in \cM : (k+1)^{-2} \le | \pm \tfrac{\pi}{2} - \vf | \le k^{-2} \} , \; \; \mbox{for all $k \ge k_0$.}
\]
For convenience, label $\mathbb{H}_0 = \cM \setminus \cup_{|k| \ge k_0} \mathbb{H}_k$.
We say a $C^1$ curve $W \subset \cM$ is stable or cone-stable if its tangent vector at each point lies in $C^{ s}$.  We call it 
homogeneous if it lies in a single homogeneity strip. 

Let $\cW^s$ denote the set of homogeneous cone-stable curves whose curvature is bounded by $B_0>0$.  We choose
$B_0$ so that $T^{-1}\cW^s \subset \cW^s$ for all $T \in \cF(\tau_*, \cK_*, E_*)$.

Since we will subdivide curves when they cross $\partial \bH_k$, we will refer to the extended 
singularity sets for $T$ (respectively $T^{-1}$) as $\cS_1^{\bH} = \cS_0 \cup T^{-1}( \cS_0 \cup_{|k| \ge k_0} \partial \bH_k)$ (respectively $\cS_{-1}^{\bH} = \cS_0 \cup T( \cS_0 \cup_{|k| \ge k_0} \partial \bH_k)$),
where $\cS_0 = \{ (r, \vf) \in \cM : \vf = \pm \pi/2 \}$.

For $\alpha \in (0,1]$, $a \ge 1$ and $W \in \cW^s$, we define the following cone of test functions,
following \cite[Section~4.1]{DL22}.  Let $d(\cdot , \cdot)$ define distance on $W$ induced by arclength.
\[
\cD_{a, \alpha}(W) := \left\{ \psi \in C^0(W) : \psi > 0, \frac{\psi(x)}{\psi(y)} \le e^{a d(x,y)^\alpha} \right\}
\]

To define the cone for $\cL_k$, we shall need a notion of distance between stable curves.  To this end, we view
each $W \in \cW^s$ as the graph of a $C^2$ function of the $r$ coordinate,
\begin{equation}
\label{eq:graph}
W = \{ G_W(r) = (r, \vf_W(r)) : r \in I_W \} .
\end{equation}
For $W^1, W^2 \in \cW^s$, if $W^1$ and $W^2$ lie in the same homogeneity strip with $|I_{W^1} \cap I_{W^2}| >0$, define the distance between them by
\[
d_{\cW^s}(W^1, W^2) = |\vf_{W^1} - \vf_{W^2}|_{C^1(I_{W^1} \cap I_{W^2})} + |I_{W^1} \bigtriangleup I_{W^2}| .
\]
Otherwise, define $d_{\cW^s}(W^1, W^2) = \infty$.  Although $d_{\cW^s}$ is not a metric, it is sufficient for our needs.

When $d_{\cW^s}(W^1, W^2) < \infty$, we will also need to measure a distance between test functions.  For
$\beta < \alpha$ as defined below and $\psi_i \in \cD_{a, \beta}(W^i)$, define,
\[
d_*(\psi_1, \psi_2) = \Big|\,\psi_1 \circ G_{W^1} \| G'_{W^1} \| - \psi_2 \circ G_{W^2} \| G'_{W^2} \|\,  \Big |_{C^\beta(I_{W^1} \cap I_{W^2})} \, ,
\]
where $\| G'_W \| = \sqrt{1 + (d\vf_W/dr)^2}$.

For a given length scale $\delta$, define 
\[
\cW^s_-(\delta) = \{ W \in \cW^s : |W| \le 2\delta \} \quad \mbox{and} \quad
\cW^s(\delta) = \{ W \in \cW^s : |W| \in [\delta, 2\delta] \} \, .
\]
Next, let $\cA_*$ denote the set of functions on $\cM$ whose restriction to each $W \in \cW^s$
is integrable with respect to the arclength measure $dm_W$.  Define,
\[
\tri h \tri_+^\sim = \sup_{\stackrel{W \in \cW^s_-(\delta)}{\psi \in \cD_{a, \beta}(W)}} \frac{| \int_W h \psi \, dm_W |}{\int_W \psi \, dm_W} , \quad \mbox{for } h \in \cA_* .
\]
Set $\cA_0 = \{ h \in \cA_* : \tri h \tri_+^\sim < \infty \}$, and note that $\tri \cdot \tri_+^\sim$
defines a seminorm on $\cA_0$.  Thus identifying two functions $g$ and $h$ in $\cA_0$ if
$\tri g - h \tri_+^\sim=0$, we define $\cA$ to be the resulting normed vector space of equivalence classes.
Note that if $g \sim h$, then $g=h$ almost everywhere with respect to both the Lebesgue measure and
$\musrb$. 

Next, for $h \in \cA$, define the following two quantities,
\begin{equation}
\label{eq:tri def}
\tri h \tri_+ = \sup_{\stackrel{\scriptstyle W \in \cW^s(\delta)}{\psi \in \cD_{a,\beta}(W)}} \frac{\left|\int_W h \psi \, dm_W\right|}{\int_W \psi \, dm_W} ,
\hskip1.3cm
\tri h \tri_- = \inf_{\stackrel{\scriptstyle W \in \cW^s(\delta)}{\psi \in \cD_{a,\beta}(W)}} \frac{\int_W h \psi \, dm_W}{\int_W \psi \, dm_W} .
\end{equation}

Denote the average value of 
$\psi$ on $W$ by $\fint_W \psi\, dm_W=\frac{1}{|W|}\int_W \psi \, dm_W$.  
Since all of our integrals on $W \in \cW^s$
will be taken with respect to the arc-length $dm_W$, to keep our notation concise, we will drop the measure from our 
integral notation in the following.

For exponents, $\alpha, \beta, \gamma, q \in (0,1)$  and constants $a, c, A, L >1$, $\delta >0$, we define the cone
\begin{align}
\cC_{c,A, L}(\delta)  =  \Bigg\{ & h \in \cA\setminus \{0\} : \nonumber\\
&\tri h \tri_+\leq L\tri h \tri_- ;
\label{eq:cone 1} \\
& \sup_{W \in \cW^s_-(\delta)} \sup_{\psi \in \cD_{a, \beta}(W)} |W|^{-q}\frac{|\int_W h \psi|}{\fint_W\psi}  \le  A \delta^{1-q} \tri h \tri_- ;
\label{eq:cone 2} \\
&\forall\, W^1, W^2 \in \cW^s_-(\delta):  d_{\cW^s}(W^1, W^2) \le \delta,  \forall \psi_i \in \cD_{a, \alpha}(W^i): d_*(\psi_1, \psi_2)=0, \nonumber\\
&\left|\frac{\int_{W^1} h \psi_1}{\fint_{W^1}\psi_1}  - \frac{\int_{W^2} h \psi_2}{\fint_{W^2}\psi_2} \right|\leq
d_{\cW^s}(W^1, W^2)^\gamma \, \delta^{1-\gamma}   c A \tri h \tri_-  \Bigg\} .
\label{eq:cone 3}
\end{align}
Let $\delta_0$ be small enough so that  \cite[eq.~(3.8)]{DL22} holds true. Next,
we assume the parameters of the cone satisfy the constraints detailed in \cite[Section~5.3]{DL22}. In particular
\begin{equation}\label{eq:const_cond}
\begin{split}
&q \in (0, 1/2),\quad 0<\beta < \alpha \le 1/3, \quad \gamma \le \min \{ \alpha - \beta, q \}\\
&3\delta<\delta_0, \quad e^{a 2(\delta_0)^\beta}   \le 2 , \quad
4A \bar C_0 \delta \le \delta_0/4,\quad c > 16 C_s^q 
\end{split}
\end{equation}
where $\bar C_0>1$ is defined in \cite[Lemma~3.3]{DL22} and $C_s := \sqrt{1 + (\cK_*^{-1} + \tau_*^{-1})^2  }$ is the maximum slope of cone-stable curves in $\cW^s$.
The fundamental property of the above cone is the following.
\begin{thm}{\cite[Theorem~6.12]{DL22}}\label{thm:cone_contraction}
Fix constants $c, L$ and $A$ satisfying the conditions of \cite[Section~5.3]{DL22}.  Then
there exists $\delta_* > 0$ such that for all $\delta \in (0, \delta_*)$ 
there exists
$N_{\cF} = N_{\cF}(\delta)$ and $\ve >0$ such that $\cL_{T_{N_\cF}} \cdots \cL_{T_1} (\cC_{c,A,L}(\delta))\subset \cC_{c,A,L}(\delta)$ with uniformly finite diameter in $\cC_{c,A,L}(\delta)$ for all choices of $Q_k \in \cQ(\tau_*, \cK_*, E_*)$ with $T_j \in \cF(Q_k, E_*; \ve)$, $j=1, \ldots , N_{\cF}$. 
\end{thm}
Since all the parameters are fixed from here forward, to simplify notation and to match notation with
Section~\ref{sec:complex}, we will refer to $\cC_{c,A,L}(\delta)$ simply as $C_{\bR}$ in what follows.

\subsection{Conditions (C-\ref{cond:positivity})-(C-\ref{cond:powerbounded})}
\label{sec:condC_bill}\ \\
Having defined the cone we can now check the conditions on the dynamics.
\begin{prop}
\label{prop:C-check}
Fix $\tau_*, \cK_*, E_*>0$, and let $\cF(\tau_*, \cK_*, E_*)$ be as defined above.  Then
conditions (C-\ref{cond:positivity})-(C-\ref{cond:powerbounded}) of Section~\ref{sec:setting} hold
with $M_k = \cM$, $\mu_k = \musrb$, $C_k = \cC_{c,A,L}(\delta)$
and $f_k = T_{N_\cF} \circ \cdots \circ T_1$,  $T_j \in \cF(Q_k, E_*; \ve)$.
\end{prop}

\begin{proof}
Remark that since $C_k = \cC_{c,A,L}(\delta)$ is the same for all $k$, then $\|\cdot\|_k=\|\cdot\|$ is independent of $k$.

Condition (C-\ref{cond:positivity}) follows from \cite[Remark 7.2]{DL22}, while Condition (C-\ref{cond:cone_cont}) follows from \cite[Theorems 6.12 and 6.13]{DL22}. 

Condition (C-\ref{cond:norm}) instead follows from the beginning of the proof of \cite[Lemma 8.15]{DL22} where $\bbe$ is chosen to be the constant function $1$.

To prove condition (C-\ref{cond:powerbounded}) note that by choice of our reference measure, we have $\cL_k1=1$. Hence, 
\[
-\|h\|\bbe=-\|h\|\cL_k \bbe\preceq \cL_k h \preceq \|h \| \cL_k \bbe \preceq \|h\|\bbe
\]
which implies $\|\cL_k h\|\leq \|h\|$, which proves (C-\ref{cond:powerbounded}) with $C_*=1$.
\end{proof}
Before introducing the observables for which we will prove the CLT and verifying conditions 
(O-\ref{cond:multiplier})-(O-\ref{cond:lambda_multiplier2}), we need to present a description of the cone adapted to the theory developed in Section~\ref{sec:complex}.
\subsection{Alternative Description of the Cone $C_{\bR}$}\ \\
\label{sec:alt cone}
So that our notation for the cone $C_{\bR}$ coincides with the formulation in Section~\ref{sec:complex}, we 
present the following equivalent description of $C_{\bR}$ in terms of the positivity of a set of linear functionals.
With the parameters $a, \alpha, \beta, \gamma, q, c, A, L, \delta$ fixed as above, 
define a corresponding set of linear functionals $\cS$ on $\cA$ as follows:
\begin{equation}
\ell_{W^1, W^2, \psi_1, \psi_2, L}(h)  =L\frac{\int_{W^1} h\psi_1}{\int_{W^1} \psi_1}-\frac{\int_{W^2} h\psi_2}{\int_{W^2} \psi_2}, \quad W^i \in \cW^s(\delta), \psi_i \in \cD_{a, \beta}(W^i) .  \label{eq:first linear} 
\end{equation}
\begin{align}
\ell_{W^0, W^1, \psi_0, \psi_1, q, \pm}(h) & =A\delta^{1-q}\frac{\int_{W^0} h\psi_0}{\int_{W^0} \psi_0}\pm |{W^1}|^{1-q}\frac{\int_{W^1} h\psi_1}{\int_{W^1} \psi_1}, \label{eq:second linear} \\
& \quad W^0 \in \cW^s(\delta), W^1 \in \cW^s_-(\delta),  \psi_i \in \cD_{a,\beta}(W^i) . \nonumber
\end{align}
\begin{align}
\ell_{W^0, W^1, W^2, \psi_0, \psi_1, \psi_2, \gamma, \pm}(h) & =d_{\cW^s}(W^1,W^2)^\gamma \delta^{1-\gamma}cA \frac{\int_{W^0} h\psi_0}{\int_{W^0} \psi_0} \label{eq:third linear} \\
& \qquad \pm \left[\frac{|W^1|\int_{W^1} h \psi_1}{\int_{W^1} \psi_1}-\frac{|W^2|\int_{W^2} h \psi_2}{\int_{ W^2} \psi_2}\right],
\nonumber \\
& \hskip-2.5cm W^0 \in \cW^s(\delta), \psi_0 \in \cD_{a, \beta}(W^0), W^1, W^2 \in \cW^s_-(\delta), \psi_i \in \cD_{a, \alpha}(W^i) , \nonumber \\
& \hskip-2.5cm d_{\cW^s}(W^1, W^2) \le \delta, d_*(\psi_1, \psi_2) = 0 .   \nonumber
\end{align}
Let $\cS$ denote the collection of all such linear functionals as $W^i$ and $\psi_i$ range over the stated 
sets of stable curves and test functions. 

\begin{lem}
\label{lem:alt}
Let $\cS$ be as defined above.  Then,
\[
C_\bR=\{h\in \cA \setminus \{0\}\;:\; \ell(h)\geq 0, \forall \ell\in\cS\}.
\] 
\end{lem}

\begin{proof}
The equivalences between the three classes of linear functionals defined in \eqref{eq:first linear}-\eqref{eq:third linear}
and the three cone conditions \eqref{eq:cone 1}-\eqref{eq:cone 3} is immediate.

Consider, for example, \eqref{eq:cone 1}.  If $\tri h \tri_+ \le L \tri h \tri_-$, then by definition of $\tri \cdot \tri_+$ and $\tri \cdot \tri_-$ it follows that
$\ell_{W^1, W^2, \psi_1, \psi_2, L}(h) \ge 0$ for all $W^1, W^2 \in \cW^s(\delta)$ and $\psi_i \in \cD_{a, \beta}(W^i)$.
On the other hand if $\ell_{W^1, W^2, \psi_1, \psi_2, L}(h) \ge 0$ for all such $W^i$ and $\psi_i$, then
it follows that 
\[
L \inf_{\stackrel{W^1 \in \cW^s(\delta)}{\psi_1 \in \cD_{a,\beta}(W^1)}} \frac{\int_{W^1} h \psi_1}{\int_{W^1} \psi_1}
\ge \sup_{\stackrel{W^2 \in \cW^s(\delta)}{\psi_2 \in \cD_{a, \beta}(W^2)}} \frac{\int_{W^2} h \psi_2}{\int_{W^2} \psi_2} ,
\] 
which is precisely \eqref{eq:cone 1}.

The other conditions are similarly parallel, with \eqref{eq:cone 2} corresponding to the set of functionals in \eqref{eq:second linear} and \eqref{eq:cone 3} corresponding to the set of functionals in \eqref{eq:third linear}. 
\end{proof}

We now define the norm $\| h \|$ precisely as in \eqref{eq:b_norms}, letting $\bbe = 1$.  Then $\cB_{\bR}$ is the completion
of $\cA$ in this norm,  It follows that $\cS \subset \cB'_{\bR}$.

Finally, following \eqref{eq:com_cone_def} we define the larger cone $\cC_{\bR} \supset C_{\bR}$ by
\[
\cC_{\bR} = \{ h \in \cB_{\bR} \setminus \{ 0 \} : \ell(h) \ge 0, \forall \ell \in \cS \} .
\] 


\subsection{Introduction of Observables  and Some Basic Properties }\ \\
\label{sec:check O1_a}
Although we work with sequential maps of the form $f_k = T_{N_\cF} \circ \cdots \circ T_1$,
$T_j \in \cF(Q_k, E_*; \ve)$ according to Proposition~\ref{prop:C-check},
we would like to prove the Central Limit Theorem for observables sampled at each billiard 
collision rather than at each $N_\cF$ collisions.  To this end, we will work with observables
of the form
\begin{equation}
\label{eq:observe}
g_k = \sum_{j=0}^{N_{\cF}-1} \tilde g_{k,j} \circ T_j \circ \cdots \circ T_1,
\end{equation}
where each $\tilde g_{k,j} \in C^\alpha(M)$,\footnote{ We use the norm $|g|_{C^\alpha}=\|g\|_\infty+\sup_{x,y}\frac{|g(x)-g(y)|}{d(x,y)^\alpha}$, so that $C^\alpha$ is a Banach algebra.} with $\alpha \in (0, 1/3)$
from the definition of $C_{\bR}$, and with $|\tilde g_{k,j}|_{C^\alpha}\leq \tilde K$ for some $\tilde K>0$
and all $k \in \mathbb{N}$, $j \le N_\cF -1$. 

For this class of observables, we will prove (O-\ref{cond:multiplier}) in several steps.  
The first step is the following Lemma, whose proof is postponed to Section~\ref{sec:lemmaO1}. 
\begin{lem}
\label{lem:norm est}
There exists $C_5>0$ such that for all $g \in C^\alpha(M)$ and $h \in \cC_\bR$,
$\| g h \| \le C_5 | g |_{C^\alpha} \tri h \tri_-$.
In particular, setting $C_\star = 4A(1+2a)\left(14C_s^q+9\delta^\gamma+5c+15\delta\right)$, the constant obtained in the proof, if
\begin{equation}
\label{eq:cop bound}
\begin{split}
&\cop >\Xi_* |g|_{C^\alpha} \tri h \tri_- \, , \, \mbox{ with } \\
&\Xi_*= \max \left\{ \frac{(3+2L)L(1+2a)}{L-1},  \frac{(5+2L) A(1+2a)}{A-2^{1-q}}  , \frac{ C_\star + 2cA(1+2a)L }{cA - 2C_s} \right\},
\end{split}
\end{equation}
then $\cop \pm gh \in \cC_\bR$ and so $\| gh \| \le \cop$.
\end{lem}

\begin{rem}
\label{rem:equiv tri}
Note that if $\cop >0$, $h \in \cC_{\bR}$ and $\cop - h \in \cC_\bR$, then necessarily,
\[
\frac{ \int_W (\cop - h) \psi}{\int_W \psi} > 0 \quad \implies \quad \cop > \frac{ \int_W h \psi}{\int_W \psi} \, , 
\]
for each $W \in \cW^s(\delta)$ and $\psi \in \cD_{a,\beta}(W)$.  Taking the supremum over $W$ and $\psi$, we
conclude that $\| h \| \ge \tri h \tri_+$.

On the other hand, according to Lemma~\ref{lem:norm est} with $g \equiv 1$, $\| h \| \le C_5 \tri h \tri_-$.  Thus
for $h \in \cC_{\bR}$, the two quantities are equivalent, i.e.
\[
\tri h \tri_- \le \tri h \tri_+ \le \|h \| \le C_5 \tri h \tri_- \, .
\] 
\end{rem}

\begin{cor}
\label{cor:in space}
If $h \in \cB_{\bR}$ and $g \in C^\alpha(M)$, then $gh \in \cB_{\bR}$ and $\| gh\| \le 3 C_5 |g|_{C^\alpha} \|h \|$, where $C_5>0$ from Lemma~\ref{lem:norm est} is independent of $g$ and $h$.
\end{cor}
\begin{proof}
For $h \in \cB_{\bR}$, $h \neq 0$, by \eqref{eq:b_norms}, $h + \| h \|, \| h \| \in \cC_{\bR} \cup \{ 0 \}$
(recalling that in the present setting $\bbe = 1$).  Then writing
$gh = g(h+\| h \|) - g\| h \|$ implies 
\[
\| gh \| \le \| g (h+\| h \|) \| + \| g \| h \| \, \| \leq C_5 |g|_{C^\alpha} \| h + \| h \| \, \| + C_5 |g|_{C^\alpha} \| h \|
\]
where we have used the fact that $\| 1\| =1$ and applied Lemma~\ref{lem:norm est} twice. 
\end{proof}

\begin{rem}\label{rem:one_step}
Theorem \ref{thm:cone_contraction} provides a bound on the norm only for sequences of operators of length $N_{\cF}$. To obtain finer information (for example, a CLT for observables sampled at each
billiard collision rather than sampled once every $N_{\cF}$ collisions) requires us  to estimate uniformly the norms of the operators $\cL_T$ for each $T \in \cF(\tau_*/2, \cK_*/2, E_*)$. Note that this does not follow from \cite{DL22}. In fact, the estimates in  \cite{DL22} cannot even be used to prove that $\cL_T\in L(\cB_\bR, \cB_\bR)$.
\end{rem}
The issue pointed out in Remark \ref{rem:one_step} is solved by the next Proposition proven in Section~\ref{sec:proof of bounded L}.
\begin{prop}
\label{lem:bounded L}
There exists $C_6 >0$ such that for all $T \in \cF(\tau_*/2, \cK_*/2, E_*)$ and all $h \in \cB_{\bR}$,
$\| \cL_T h \| \le C_6 \| h \|$.
\end{prop}

\begin{cor}
\label{cor:bounded twist}
Fix $\lambda_0>0$ and $T \in \cF(\tau_*/2, \cK_*/2, E_*)$.  For $n \in \mathbb{N}$, $|\lambda| \le \lambda_0 \sigma_n$, and $g \in C^\alpha(M)$ with $|g|_{C^\alpha}\leq \tilde K$, define
$\cL_{T,\lambda} h  = \cL_T\big(e^{i \lambda \sigma_n^{-1} g} h \big)$ for $h \in \cB_{\bR}$.  Then 
\[
\| \cL_{T,\lambda} h \| \le (1+\lambda_0) C_7 \| h\| \, ,
\]
where $C_7:=3C_5C_6 \tilde K$.
\end{cor}
\begin{proof}
Using Proposition~\ref{lem:bounded L} and then Corollary~\ref{cor:in space}, we estimate
\[
\| \cL_{T,\lambda} h  \| \le C_6\, 3C_5 | e^{i \lambda \sigma_n^{-1} g}|_{C^\alpha} \| h \| \, .
\]
Finally, using $|e^{i \theta} -1| \le |\theta|$ yieds 
$|e^{i \lambda \sigma_n^{-1} g} |_{C^\alpha} \le 1 + |\lambda| \sigma_n^{-1} |g|_{C^\alpha}\leq 1+\lambda_0\tilde K$.
\end{proof}

Note that Corollaries \ref{cor:in space} and \ref{cor:bounded twist} are stated for $h\in\cB_\bR$, while below they are applied to complex valued functions. By Lemma \ref{lem:complexification} both extend to $\cB_\bC$ at the price of a factor $\sqrt 2$, which we absorb in the constants.


\subsection{Verification of Hypothesis (O-\ref{cond:multiplier})}\ \\ \label{sec:check O1}

{We define  $\cL_k = \cL_{T_{N_\cF}} \cdots \cL_{T_1}$, $T_j \in \cF(Q_k, E_*; \ve)$.
For fixed $n \in \mathbb{N}$, $\lambda_0>0$, and $j \le k \le n$, $|\lambda| < \lambda_0 \sigma_n$, define the complex operators 
$\cL_{k,\lambda}$ and their concatenations
 as in \eqref{eq:transfer_gen},
\[
\cL_{k ,\lambda} h = \cL_k (e^{i \lambda \sigma_n^{-1}  \hat g_k} h), \quad
\cL_{k,j,\lambda} = \cL_k \cdots \cL_j( e^{i \lambda \sigma_n^{-1} \sum_{m=j}^k \hat g_m \circ f_{m,j}} h) \, ,
\]
where $\hat g_k$ are the centered versions as in \eqref{eq:centering0}, relative to the initial density $\rho$ of Theorem \ref{thm:billiard CLT}.

To verify (O-\ref{cond:multiplier}), we begin with $j=1$ and disregard the centering of the
observables, since as remarked after (O-\ref{cond:lambda_multiplier}), the centering does not
affect the estimate of the norm.
  
We must bound the norm of $\cL_{k, \lambda} (g_k h)$ for $h \in \cB_{\bR}$.
Since $\cL_{k, \lambda} = \cL_{T_{N_\cF}, \lambda} \cdots \cL_{T_1, \lambda}$, where $\cL_{T_i,\lambda}$ stands for the operator of Corollary \ref{cor:bounded twist} with $g=\tilde g_{k,i-1}$. We have
\[
\cL_{ k, \lambda} (g_k h)  = \sum_{j=0}^{N_\cF-1} \cL_{ k, \lambda}( \tilde g_{k,j} \circ T_j \circ \cdots \circ T_1 \cdot h)  = \sum_{j=0}^{N_\cF -1} \cL_{T_{N_\cF}, \lambda}  \cdots \cL_{T_{j+1}, \lambda} ( \tilde g_{k,j} \cL_{T_j, \lambda} \cdots \cL_{T_1, \lambda} h) . 
\]
Using the triangle inequality and applying Corollary~\ref{cor:bounded twist} $N_{\cF} - j$ times to the $j$th term yields,
\begin{equation}
\label{eq:L_k norm}
\begin{split}
\| \cL_{k, \lambda}(g_k h) \| & \le \sum_{j=0}^{N_\cF-1} (C_7(1+\lambda_0))^{N_\cF-j} \| \tilde g_{k,j} \cL_{T_j, \lambda} \cdots \cL_{T_1, \lambda} h \| \\
& \le \sum_{j=0}^{N_\cF-1}  ( C_7 (1+\lambda_0))^{N_\cF} 3 C_5 | \tilde g_{k,j}|_{C^\alpha} \| h \| 
\le 3 N_{\cF} (C_7(1+\lambda_0))^{N_\cF}  C_5 \tilde K \| h \| \, ,
 \end{split}
\end{equation}
where we have applied Corollary~\ref{cor:in space} followed by Corollary~\ref{cor:bounded twist} again in the second line.
This proves (O-\ref{cond:multiplier}) for $j=1$ with $K = 3 C_5 \tilde K N_{\cF} (C_7(1+\lambda_0))^{N_\cF}$.

For $j=2$, we proceed similarly.  Now,
\[
\begin{split}
g_k^2 & = \sum_{i,j=0}^{N_\cF-1} \tilde g_{k,j} \circ T_j \circ \cdots \circ T_1 \cdot \tilde g_{k,i} \circ T_i \circ \cdots T_1 \\
& = \sum_{j=0}^{N_\cF-1} \tilde g_{k,j}^2 \circ T_j \circ \cdots \circ T_1 + 2 \sum_{j=1}^{N_\cF-1} \sum_{i=0}^{j-1}
\tilde g_{k,j} \circ T_j \circ \cdots \circ T_1 \cdot \tilde g_{k,i} \circ T_i \circ \cdots T_1 \, .
\end{split}
\]
Multiplying by $h$ and applying $\cL_{k, \lambda} = \cL_{T_{N_\cF}, \lambda} \cdots \cL_{T_1, \lambda}$ yields,
\[
\begin{split}
\cL_{k, \lambda}(g_k^ 2 h) & =  \sum_{j=0}^{N_\cF-1} \cL_{T_{N_\cF}, \lambda \cdots \cL_{T_{j+1}, \lambda} ( \tilde g_{k,j}^2 \cL_{T_j, \lambda} \cdots \cL_{T_1, \lambda} h) } \\
& \quad + 2 \sum_{j=1}^{N_\cF-1} \sum_{i=0}^{j-1}
\cL_{T_{N_\cF}, \lambda} \cdots \cL_{T_{j+1}, \lambda} ( \tilde g_{k,j} \cL_{T_j, \lambda} \cdots \cL_{T_{i+1}, \lambda} ( \tilde g_{k,i} \cL_{T_i, \lambda} \cdots \cL_{T_1, \lambda} h)) \, . 
\end{split}
\]
And applying Corollary~\ref{cor:bounded twist} and Corollary~\ref{cor:in space} as in
\eqref{eq:L_k norm} yields,
\[
\begin{split}
\| \cL_{k, \lambda}(g_k^2 h) \| & \le N_\cF ( C_7(1+\lambda_0))^{N_\cF} 3 C_5 \tilde K^2 \| h \| + N_{\cF}(N_\cF-1) ( C_7 (1+\lambda_0))^{N_\cF} 9 C_5^2 \tilde K^2 \| h \| \\
& 
\le 9 C_5^2 N_\cF^2 \tilde K^2  (C_7(1+\lambda_0))^{N_\cF} \| h \| ,
\end{split}
\]
where we have used the fact that $| \tilde g_{k,j}^2 |_{C^\alpha} \le |\tilde g_{k,j} |_{C^\alpha}^2$.
This proves (O-\ref{cond:multiplier}) for $j=2$ since the constant is $\le K^2$.

For $j = 3$ the process is similar, using the fact that
\[
g_k^3 = \sum_{i,j,\ell=0}^{N_\cF-1}  \tilde g_{k,\ell} \circ T_\ell \circ \cdots \circ T_1 \cdot \tilde g_{k,j} \circ T_j \circ \cdots \circ T_1 \cdot \tilde g_{k,i} \circ T_i \circ \cdots T_1 \, .
\]
Then, estimating precisely as in the case $j=2$, we obtain,
\[
\| \cL_{k, \lambda} (g_k^3 h) \| \le N_\cF^3 27 C_5^3 \tilde K^3 ( C_7(1+\lambda_0))^{N_\cF} \| h\| \, ,
\]
which is the required estimate for $j=3$.  The same strategy works for any $j \in \bN$, completing the verification
of (O-\ref{cond:multiplier}) with $j_0 = \infty$.

\subsection{Contraction of the Complex Cone $\cC_{\bC}$}\ \\
\label{sec:cone contract}
With $\cB_{\bR}$ and $\cC_{\bR}$ defined as in Section~\ref{sec:alt cone}, we proceed to define
their complexifications $\cB_{\bC}$ and $\cC_{\bC}$ precisely as in Section~\ref{sec:complex},
specifically Definition~\ref{def:complex cone}.

In light of Remark~\ref{rem:equiv tri}, we define 
\begin{equation}
\bbm(h) := \int_M h \, d\musrb \qquad \mbox{for $h \in C_{\bR}$.}
\end{equation}
By \cite[Remark~7.2]{DL22}, there exists $\bar{C} \ge 1$ such that
\begin{equation}
\label{eq:m upper lower}
\bar{C}^{-1}\tri h \tri_- \le \bbm(h) \le \bar{C} \tri h \tri_+ ,
\end{equation}
and since $\tri \cdot \tri_-$ is equivalent
to $\| \cdot \|$ by Remark~\ref{rem:equiv tri}, $\bbm$ extends to a linear functional on $\cC_{\bR}$ and thus to $\cB_{\bR}$ and $\cB_{\bC}$.  Indeed, disintegrating $\musrb$ as in the proof of
\cite[Lemma~7.1]{DL22} it follows that $\bbm \in \cS_*$, where, as in Section~\ref{sec:complex}
(see definition \eqref{eq:maybetoomuch}), $\cS_*$ denotes the weak-* closure of the convex hull of 
$\{ \cop \ell : \cop \in \bR_+, \ell \in \cS \}$.
Now $\bbm(1) = 1$ and combining \eqref{eq:m upper lower} with Remark~\ref{rem:equiv tri}
implies that $\bbm$ satisfies  
\eqref{eq:maybetoomuch} and therefore \eqref{eq:m_lower} with $\kappa = C_5^{-1} \bar{C}^{-1}$.

The main result of this section is the following.

\begin{thm}
\label{thm:complex contract}
Let $\cF(\tau_*, \cK_*, E_*)$ and $f_k = T_{N_\cF} \circ \cdots \circ T_1$ be as in Proposition~\ref{prop:C-check}.
Let $\cC_{\bR}$ be as defined in \eqref{eq:cone 1}-\eqref{eq:cone 3} and Lemma~\ref{lem:alt},
and let $\cC_{\bC}$ be its complexification as in Definition~\ref{def:complex cone} with
associated projective metric $\delta_{\cC}$ as in \eqref{eq:gauge_dist}.

There exists $\lambda_0 > 0$ such that for all $k , n \in \mathbb{N}$, $k \le n$, and 
all $\lambda \in \mathbb{R}$ such that $|\lambda| \sigma_n^{-1} < \lambda_0$, 
$\cL_{k, \lambda}( \cC_{\bC} ) \subset \cC_{\bC}$ and $\diam_{\delta_{\cC}}(\cL_{k,\lambda}(\cC_{\bC})) =: \Delta_{\bC}  < \infty$.
\end{thm}

\begin{proof}
We will apply Theorem~\ref{thm:needed} to the operators $\cL_{k, \lambda}$.
First, by Proposition~\ref{prop:C-check}, the real transfer operators $\cL_k$ satisfy
$\cL_k \in L(\cB_{\bR})$ and $\diam_H(\cL(\cC_{\bR})) := \Delta_{\bR} < \infty$;
since $\bbe =1$ and $\cL_k 1 = 1$, we may take $B=0$ in the statement of 
Theorem~\ref{thm:needed} and in item (e) of Section~\ref{sec:lower bound}, see footnote~\ref{foot:B}.
Moreover, Lemma~\ref{lem:norm est} with $g \equiv 1$ implies that $\bbe = 1$ satisfies
\eqref{eq:archimedean}, while \eqref{eq:m upper lower} and Remark~\ref{rem:equiv tri}
imply that $\bbm( \cdot) = \musrb(\cdot)$ satisfies \eqref{eq:maybetoomuch}.

In the following, we use the uncentered twist, that is $\cL_{k,\lambda}h=\cL_k(e^{i\lambda\sigma_n^{-1}g_k}h)$ which differs from the centered one only by a constant phase.

Next, for $h \in \cB_{\bR}$,
\begin{equation}
\label{eq:real im}
\cL_{k, \lambda} h = \cL_k (e^{i \lambda \sigma_n^{-1} g_k} h) = \cL_k( \cos(\lambda \sigma_n^{-1} g_k) h) + i \cL_k( \sin(\lambda \sigma_n^{-1} g_k) h ) \, .
\end{equation} 
Then temporarily labelling $t = \lambda \sigma_n^{-1}$ for brevity, we can write, 
\[
\cos (t (\tilde g_{k,1} \circ T_1 + \tilde g_{k,0}) )
= \cos (t \tilde g_{k,1}) \circ T_1 \cos(t \tilde g_{k,0})  - \sin (t \tilde g_{k,1}) \circ T_1 \sin (t \tilde g_{k,0}) \, ,
\]
so that inductively, using Corollary~\ref{cor:in space} and estimating as in \eqref{eq:L_k norm},
it follows that $\cL_{k, \lambda}$ defines a bounded operator on $\cB_{\bC}$
by Lemma~\ref{lem:complexification}.

It remains to verify condition \eqref{eq:complex_bound} of Theorem~\ref{thm:needed}. We will do this via
Lemma~\ref{lem:pert_cond}, which implies \eqref{eq:complex_bound}.  Recalling
\eqref{eq:real im}, Lemma~\ref{lem:pert_cond} with $a(h) = \cL_k(\cos(\lambda \sigma_n^{-1} g_k) h)$ and
$b(h) = \cL_k( \sin(\lambda \sigma_n^{-1} g_k) h )$
follows immediately from the following lemma.

\begin{lem}
\label{lem:small pert check}
Let $g_k$ be as in \eqref{eq:observe}.
For all $\ve > 0$ there exists $t_0 > 0$ such that for all $h \in \cC_{\bR}$, all $f_k$ and all $t \in (-t_0, t_0)$,
\[
\begin{split}
& \cL_k  [\ve h \pm 2 h (1-\cos(t g_k) ) ] \in \cC_{\bR} \, ,\\
& \cL_k [\ve h \pm 2 h \sin(t g_k) ] \in \cC_{\bR} \, .
\end{split}
\]
\end{lem}

\begin{proof}[Proof of Lemma~\ref{lem:small pert check}]
Let $\psi_t = 2(1-\cos(t g_k))$ and let $h \in \cC_{\bR}$.  
First we claim that
\begin{equation}
\label{eq:O t}
\cL_k(\psi_t h) \in \cB_{\bR} \quad \mbox{and} \quad \| \cL_k (\psi_t h) \| \le C | t| \| h\| \, , 
\end{equation}
for some $C>0$ independent of $t$, $k$ and $h$.

To prove the claim, expand $\cos(t g_k)$ using \eqref{eq:observe} to obtain,
\begin{equation}
\label{eq:cos expand}
\cos(t g_k) = \prod_{j=0}^{N_\cF-1} \cos (t \tilde g_{k,j} \circ T_j \circ \cdots \circ T_1) + P_k(t) \, ,
\end{equation} 
where $P_k(t)$ is a polynomial of degree $N_\cF$ in terms involving $\cos(t \tilde g_{k,j})$ and 
$\sin(t \tilde g_{k,j})$.  Indeed, every term in $P_k(t)$ has at least one factor of the form
$\sin(t \tilde g_{k,j} \circ T_j \circ \cdots \circ T_1)$.

First we estimate the leading term of $\cL_k(\psi_t h)$,
\begin{equation}
\label{eq:telescope}
\begin{split}
&\cL_k  \Big( \Big(1 - \prod_{j=0}^{N_\cF-1} \cos (t \tilde g_{k,j} \circ T_j \circ \cdots \circ T_1) \Big)h \Big)
=  \cL_k \big( \big( 1 - \cos(t \tilde g_{k,0}) \big) h \big) \\
& + \sum_{\ell = 1}^{N_\cF-1} \cL_k \Big( \prod_{j=0}^{\ell-1} \cos (t \tilde g_{k,j} \circ T_j \circ \cdots \circ T_1)
- \prod_{j=0}^{\ell} \cos (t \tilde g_{k,j} \circ T_j \circ \cdots \circ T_1) \Big) h \Big)  .
\end{split}
\end{equation}
By Lemma~\ref{lem:norm est} and Proposition~\ref{lem:bounded L}, since $h\in\cC_\bR$ and $1-\cos(t \tilde g_{k,0})\in C^\alpha$,
\[
\| \cL_k \big( \big( 1 - \cos(t \tilde g_{k,0}) \big) h \big) \| \le C_6^{N_\cF}  C_5 |1- \cos(t \tilde g_{k,0})|_{C^\alpha} \| h\|
\le C' |t| \| h \| \, .
\]
Next, for each $\ell$, using that $\cL_{T_j} 1 = 1$, 
\begin{equation}
\label{eq:split L}
\begin{split}
\cL_k & \Big( \prod_{j=0}^{\ell-1} \cos (t \tilde g_{k,j} \circ T_j \circ \cdots \circ T_1)
- \prod_{j=0}^{\ell} \cos (t \tilde g_{k,j} \circ T_j \circ \cdots \circ T_1) \Big) h \Big) \\
& = \cL_k \Big( \big( 1 - \cos(t \tilde g_{k,\ell} \circ T_\ell \circ \cdots \circ T_1 \big) \prod_{j=0}^{\ell-1} \cos (t \tilde g_{k,j} \circ T_j \circ \cdots \circ T_1) h \Big) \\
& = \cL_{T_{N_\cF}} \cdots \cL_{T_{\ell+1} } \Big( \big(1 - \cos(t \tilde g_{k, \ell}) \big) 
\cL_{T_\ell} \big (\cos (t \tilde g_{k,\ell-1})\\
&\phantom{= \;}
\times  \cL_{T_{\ell -1} } \big( \cos(t \tilde g_{k, \ell-2}) \cdots  \cL_{T_1} \big( \cos ( t \tilde g_{k,0}) h \big) \big) \big) \Big).
\end{split}
\end{equation}
Again using Corollary~\ref{cor:in space} and Proposition~\ref{lem:bounded L}, the norm of this expression is bounded by
\[
\begin{split}
 C_6^{N_\cF - \ell -1} & 3C_5 | 1 - \cos (t\tilde g_{k,\ell}) |_{C^\alpha} \left\| \cL_{T_\ell} \left (\cos (t \tilde g_{k,\ell-1}) \cL_{T_{\ell -1} } \big(  \cdots  \cL_{T_1} \big( \cos ( t \tilde g_{k,0}) h \big) \big) \right) \right\| \\
& \le C_6^{N_\cF} (3C_5)^{\ell+1} \tilde K^\ell  | 1 - \cos (t\tilde g_{k,\ell}) |_{C^\alpha} \| h\|
\le C'' |t| \| h\| \, .
\end{split}
\]
Combining these estimates in \eqref{eq:telescope} and using the triangle inequality proves the required bound
for the leading term of $\cL_k(\psi_t h)$ from \eqref{eq:cos expand}.  According to \eqref{eq:cos expand},
the remaining terms of $\cL_k(\psi_t h)$ are given by $\cL_k( P_k(t) h)$.  These terms can be estimated 
in a similar manner as the product of cosines.  
Since each term of $\cL_k( P_k(t) h)$ contains one factor of the form
$\sin (t \tilde g_{k,j})$ and $| \sin(t \tilde g_{k,j})|_{C^\alpha} = \mathcal{O}(|t|)$, the claim \eqref{eq:O t} is proved.

Since 
$\cL_k(\psi_t h) \in \cB_{\bR}$, there exists $E_t \in \bR^+$ with 
$E_t \le C |t| \tri h \tri_-$ such that
$E_t \pm \cL_k(\psi_t h) \in \cC_{\bR}$.  Then,
\[
\cL_k (\ve h \pm h \psi_t) = \cL_k (\ve h - E_t) + E_t \pm \cL_k (h \psi_t)
= \ve (\cL_k h - \ve^{-1} E_t) + E_t \pm \cL_k (h \psi_t) \, ,
\]
where we have used that $\cL_k 1 =1$.  The second term is an element of $\cC_{\bR}$ by choice of $E_t$.

For the first term, we use a stronger fact: that $\diam_H( \cL_k(\cC_{\bR})) \le \Delta_{\bR} < \infty$.  
In fact, by \cite[Prop. 6.13]{DL22} for $h \in \cC_{\bR}$, we have
$\alpha(1, \cL_k h ) \ge \frac{(1-\chi) A}{A+1} \tri \cL_k h \tri_-$, where $A$ is the cone parameter from \eqref{eq:cone 2}
and $\chi <1$ is the contraction
in cone parameters.

It follows that
\[
\cL_k h - \frac{E_t}{\ve} \in \cC_{\bR} \quad \mbox{if} \quad \frac{E_t}{\ve} < \alpha(1, \cL_k h) , \quad
\mbox{which is implied by} \quad \frac{E_t}{\ve} < \frac{(1-\chi) A}{A+1} \tri \cL_k h \tri_- \, .
\]
By choice of $E_t$, \eqref{eq:m upper lower} together with \eqref{eq:cone 1} applied to $\cL_k h\in\cC_\bR$ yields $\tri h \tri_-\leq \bar C\bbm(\cL_k h)\leq \bar C^2\tri \cL_k h\tri_+\leq \bar C^2 L\tri \cL_k h\tri_-$. Hence, $E_t \le C |t| \tri h \tri_- \le \bar C^2 L C |t| \tri \cL_k h \tri_-$,
so the above condition is satisfied if
\[
|t| < \ve \frac{(1-\chi) A}{ C\bar C^2 L (A+1)} \, .
\]
Taking this as our value of $t_0$ proves the first statement of the lemma.

Similarly, setting $\tilde \psi_t = 2 \sin(t g_k)$, we can expand $2 \sin(t g_k) = \tilde P_k(t)$ where $\tilde P_k(t)$
is a polynomial in $\sin(t \tilde g_{k,j} \circ T_j \circ \cdots \circ T_1)$ and $\cos(t \tilde g_{k,j} \circ T_j \circ \cdots \circ T_1)$ and each term of $\tilde P_k(t)$ has a factor of the form $\sin(t \tilde g_{k,j} \circ T_j \circ \cdots \circ T_1)$.
Then estimating as in \eqref{eq:split L}, and using the fact that $| \sin (t \tilde g_{k,j}) |_{C^\alpha} = \mathcal{O}(|t|)$,
we conclude,
\[
\cL_k(\tilde \psi_t h) \in \cB_{\bR} \quad \mbox{and} \quad \| \cL_k(\tilde \psi_t h) \| \le \tilde C |t| \|h \| \, ,
\]
for some $\tilde C>0$ independent of $t$, $k$ and $h$.  From this point, the same proof with $\tilde \psi_t$ replacing
$\psi_t$ yields the second statement of the lemma. 
\end{proof}
Fixing\footnote{ Recall that we take $B=0$ when applying Theorem~\ref{thm:needed},
since $\bbe =1 = \cL_k 1 \in \cL_k(\cB_{\bR})$.} $\ve < \frac{\kappa^2}{12 \sqrt{2}} e^{-2 \Delta_{\bR}}$, let $t_0>0$ be from Lemma~\ref{lem:small pert check}.
Set $\lambda_0 = t_0$.  Then for all $\lambda \in \mathbb{R}$ such that
$|\lambda| \sigma_n^{-1} < \lambda_0$, setting $t = \lambda \sigma_n^{-1}$, we have $|t| < t_0$ so that Lemma~\ref{lem:small pert check} implies that
Lemma~\ref{lem:pert_cond} and therefore \eqref{eq:complex_bound} holds for $\cL_{k, \lambda}$.
This completes the proof of the theorem.
\end{proof}


\subsection{Verification of (O-\ref{cond:lambda_multiplier}) and (O-\ref{cond:lambda_multiplier2})}\ \\
\label{sec:check O3}
Throughout this section, we assume $|\lambda| < \lambda_0 \sigma_n$, where $\lambda_0$ is from Theorem~\ref{thm:complex contract}.
From our choices of $\bbe_k = 1$, $\bbm_k(h) = \int_M h \, d\musrb$, property (a)
of Section~\ref{sec:memory loss} holds with $\bar\kappa = \kappa = C_5^{-1}\bar C^{-1}$ from
Section~\ref{sec:cone contract}.  Moreover, property (b) of Section~\ref{sec:memory loss} is trivially
satisfied with $K=1$, and property (c)
follows from Theorem~\ref{thm:complex contract} with $\bL_k = \cL_{k, \lambda, n}$.

 Hence setting $\alpha_{k,j,\lambda} = \boldsymbol{\alpha}_{k,j}$, $h_{k,j,\lambda} = \boldsymbol{h}_{k,j}$ and $\ell_{k,j,\lambda} = \boldsymbol{\ell}_{k,j}$,
 Lemma \ref{lem:boundedness} implies
\[
\|h_{k,j,\lambda}\|\leq \frac{\sqrt 2}{\kappa},\quad \|\ell_{k,j,\lambda}\|'\leq \sqrt 2
\]
while $|\alpha_{k,j,\lambda}| \le 1$, by the invariance of the $SRB$ measure. This implies,
\[
\| \cL_{k,j,\lambda} h \| \le  \frac{2(1+2\sqrt 2)}{\kappa} \| h \|,
\]
which verifies (O-\ref{cond:lambda_multiplier}) with $K = 2(1+2\sqrt 2)/\kappa$.

Condition (O-\ref{cond:lambda_multiplier2}) is verified exactly as in Section \ref{sec:checkO3}
 since Theorem~\ref{thm:complex contract} implies that Theorem~\ref{thm:complex_contraction} applies to the 
operators $\cL_{k,\lambda,n}$.


\subsection{Central Limit Theorem for Sequential Dispersing Billiards}
\label{sec:CLT billiard}

We collect our previous results on the contraction of the relevant operators in order to state the culminating
Central Limit Theorem for sequential dispersing billiards.  This is essentially a restatement of 
Theorem~\ref{thm:sample}, but now all the relevant objects have been properly defined.  

\begin{thm}
\label{thm:billiard CLT}
Fix $\tau_*, \cK_*, E_*>0$ and let $\cF(\tau_*, \cK_*, E_*)$ be the associated family of billiard maps.  Let
$\ln\rho \in \cC^1$, $\int \rho d\musrb=1$.
For a sequence of observables $(\tilde g_{k,j})$ as defined in \eqref{eq:observe}, choose an $N_\cF$-admissible sequence
$(T_i)_{i \in \mathbb{N}}$ as defined in Definition~\ref{def:admissible}, and define the centered observables 
$\hat{g}_k$ as in \eqref{eq:centering0}.  Then, recalling the distribution function $F_n(x)$ from \eqref{eq:Fn def},
for all $c_\star \in (0,1)$  there
exists $\bar C>0$ such that, for all $n \in \mathbb{N}$ such that $\sigma_n \ge \max \{ 1, c_\star n^{1/3} \ln(n+1) \}$, we have
\[
\left|F_n(x)-\frac{1}{\sqrt{2\pi}}\int_{-\infty}^xe^{-\frac {y^2}2}dy\right|\leq \bar C \sigma_n^{-3}(\ln \sigma_n)^2 n.
\]
\end{thm}

\begin{proof}
We have verified hypotheses (C-\ref{cond:positivity})-(C-\ref{cond:powerbounded}) in Section \ref{sec:condC_bill} and 
conditions (O-\ref{cond:multiplier})-(O-\ref{cond:lambda_multiplier2}) in Section \ref{sec:check O3}, accordingly
the theorem is an immediate application of Theorem~\ref{thm:main_bis} and Corollary~\ref{cor:better corollary}.
\end{proof}

The remaining issue is to verify the condition for the growth of the variance. \\
We start with the following general estimate, recalling $h_k = \cL_{k-1}\cdots \cL_0 \rho$,
\begin{equation}\label{eq:baseVarS_n}
\begin{split}
&\sigma_n^2:=\int [\hat S_n](x)^2\rho(x) d \mu\\
&=\sum_{k=0}^{n-1}\int \hat g_k^2 h_k+2\sum_{k=0}^{n-1}\sum_{j=k+1}^{n-1}\int \hat g_j\cL_{j-1}\cdots\cL_k (\hat g_k h_k)\\
&=\sum_{k=0}^{n-1}\int \hat g_k^2 h_k+2\sum_{k=0}^{n-1}\sum_{j=1}^{\infty}\int \hat g_{j+k}\cL_{j+k-1}\cdots\cL_k (\hat g_k h_k)
+\cO\left(\sum_{k=0}^{n-1}\sum_{j=n}^{\infty}\nu^{j-k}\right)\\
&=\sum_{k=0}^{n-1}\left[\int \hat g_k^2 h_k+2\sum_{j=1}^{\infty}\int \hat g_{j+k}\cL_{j+k-1}\cdots\cL_k (\hat g_k h_k) \right]
+\cO(1),
\end{split}
\end{equation}
where we have used  hypotheses (O-\ref{cond:multiplier}) together with \eqref{eq:decay_real}, with $\nu=e^{-c}$. Note that \eqref{eq:decay_real} applies due to the centering given
by \eqref{eq:centering}. 
In the case of billiards, choosing $\rho=1$ yields  $h_k=1$. Also, for simplicity, we 
set $d\mu = d\musrb$ and choose  $g_k$ such that $\int g_k=0$, thus $g_k=\hat g_k$.

\begin{lem}\label{lem:poor-man-growth} 
If there exists $ C_g>0$ such  that $ C_g\int \hat g_k^2 \geq  K\|\hat g_k\|_k$ for all $k\in\bN$, $K$ being the constant in (O-\ref{cond:multiplier}), then, for $\nu$ small enough (depending only on $C_g$)  there exists $C_8>0$, independent of $n$, such that, for all $n\in\bN$,
\[
\sigma_n^2\geq \frac 34\sum_{k=0}^{n-1} \int \hat g_k^2-C_8 .
\]
In particular, if $\inf_k\int \hat g_k^2=:m>0$, then $\sigma_n^2\geq \frac 12 m n$ for all $n\geq 4C_8m^{-1}$.
\end{lem}
\begin{proof}
We use condition (O-\ref{cond:multiplier}) and \eqref{eq:decay_real} (recall that we have set $e^{-c}=\nu$) to compute
\[
\begin{split}
&\int \hat g_k^2 +2\sum_{j=1}^{\infty}\int \hat g_{j+k}\cL_{j+k-1}\cdots\cL_k \hat g_k=\int \hat g_k^2 +2\sum_{j=1}^{\infty}\int\cL_{j+k} \hat g_{j+k}\cL_{j+k-1}\cdots\cL_k \hat g_k\\
 &\geq \int\hat g_k^2- 2K\sum_{j=1}^{\infty}\|\cL_{j+k-1}\cdots\cL_k \hat g_k\|_{j+k}\geq \int \hat g_k^2 -\frac{6K \nu }{1-\nu} \|\hat g_k\|_k\\
&\geq \left[1- \frac{6 C_g\nu}{1-\nu} \right]\int \hat g_k^2.
\end{split}
\]
The formula \eqref{eq:baseVarS_n} implies the lemma, $C_8$ being the constant implied by the $\cO(1)$ in \eqref{eq:baseVarS_n}, which is bounded by $\Const\nu(1-\nu)^{-2}$, uniformly in $n$.
\end{proof}
The above Lemma, although effective, is not very satisfactory as it is not easy to check and applies only if $\nu$ is small. It would be nice to have the equivalent of Proposition \ref{prop:variance}; unfortunately, it is not clear how to proceed. It seems natural to introduce a stable foliation, as in \cite{Li95c}, but there are some technical problems that require extra ideas.

In the literature, there are some relevant results only in the random case  \cite{DFGV18, DFGV18bis, DH21}, described in the next section. For completeness, we show that the random case can be treated also in our setting.


\subsection{Variance Growth in the Random Case}\label{sec:random}\ \\
Let $\Omega=\fA^\bN$ be a probability space with a measure $\Pe$. Assume that $\Pe$ is shift invariant and ergodic with respect to the shift and that $\Pe(\{\omega_1=\fa\})>0$ for each $\fa\in\fA$.

Let $\{f_{\fa}\}_{\fa\in\fA}$ be a set of  maps of a compact Riemannian manifold $M$ satisfying conditions 
(C-\ref{cond:positivity})-(C-\ref{cond:powerbounded}) of Section~\ref{sec:abstract} with constants independent of $\fa\in\fA$. In the case of billiards, the index $\fa$ labels blocks of the form $T_{N_\cF}\circ\cdots\circ T_1$ of Section \ref{sec:billiard}, since (C-\ref{cond:cone_cont}) is verified for such blocks rather than for the single maps.  To simplify the exposition for this calculation of the variance, we take the manifold $M$ and cone $C$ to be fixed,
while the maps are random. Also, we assume that the maps $f_{\fa}$ preserve a common invariant measure 
$\mu \in \cB$. 
We define $\cL_{\fa}$ to be the transfer operator of $f_{\fa}$ with respect to $\mu$.  This implies
in particular that $\cL_{\fa} 1 = 1$.

Let  $\{g_{\fa}\}_{\fa\in\fA}\subset L^2(M,\bR) \cap \cB'$ be a set of observables satisfying (O-\ref{cond:multiplier})
with constants uniform for $\fa\in\fA$.  Note that in Section \ref{sec:check O1} it is proven that this is true if the $f_{\fa}$ are billiard maps in the family $\cF(\tau_*, \cK_*, E_*)$.
Suppose that, for each $\fa\in\fA$,
\[
\int_M g_{\fa} d\mu =0.
\]
 Let  $\bbg(\omega,x)=g_{\omega_0}(x)$ for each $(\omega,x)\in\Omega\times M$ and consider the dynamics $F:\Omega\times M\to\Omega\times M$ defined by
\begin{equation}\label{eq:skew}
F(\omega,x)=(\tau\omega, f_{\omega_0}(x)).
\end{equation}
For each $\omega\in\Omega$, we have the sequential observables 
\[
\{\bbg\circ F^k\}=\{g_{\omega_k}\circ f_{\omega_{k-1}}\circ\cdots \circ f_{\omega_0}\},
\]
hence, for a given $\omega\in\Omega$, we can write
\[
\begin{split}
&\hat S_n(\omega,x)=\sum_{k=0}^{n-1}g_{\omega_k}\circ f_{\omega_{k-1}}\circ\cdots \circ f_{\omega_0}(x)=\sum_{k=0}^{n-1}\bbg\circ F^k(\omega, x)\\
&\sigma_n^2(\omega)=\int_M \hat S_n(\omega,x)^2 d\mu.
\end{split}
\]
The reader can check \cite{DFGV18, DFGV18bis, DH21} and references therein for a discussion of a similar model (but limited to expanding maps) and \cite{ALS09} (for the case of toral automorphisms) using a related approach. The following result is a generalization of \cite{ALS09}; similar ideas can be applied to the case of expanding maps.

\begin{thm}\label{thm:random_var} $\Pe$ is i.i.d.. Then $\Pe$-a.s. the variance $\sigma_n$ either grows linearly or there exists $\hat \psi\in L^2(M)$ such that $g_{\fa}=  \hat\psi-\hat\psi\circ f_{\fa}$ for all $\fa\in\fA$, and hence $\sigma_n$ is uniformly bounded.
\end{thm}
\begin{proof}
Let $\beta:\Omega\to \bR$ be defined by
\[
\beta(\omega)=\int_{M} g_{\omega_0}^2 d\mu +2 \sum_{j=1}^\infty\int_{M} g_{(\tau^{j}\omega)_0}\cL_{(\tau^{j-1}\omega)_0}\cdots \cL_{\omega_0} g_{\omega_0} d\mu ,
\]
where $\tau$ is the shift. Then by \eqref{eq:baseVarS_n}  (taking $\rho=1$ there),
applied to the sequence determined by $\omega\in\Omega$,
\begin{equation}\label{eq:sigmabound}
\sigma_n^2(\omega)=\sum_{k=0}^{ n-1} \beta(\tau^k\omega)+\cO(1).
\end{equation}
\begin{lem}\label{lem:alternative}
Either $\sigma_n^2$ grows linearly $\Pe$-a.s. or $\bE_{\Pe}(\sigma_n^2)=\cO(1)$.
\end{lem}
\begin{proof}
Integrating \eqref{eq:sigmabound} yields $\bE_{\Pe}(\sigma_n^2)=n\bE_{\Pe}(\beta)+\cO(1)$.
It follows, by Birkhoff's ergodic theorem, that  $\Pe$ almost surely
\[
\lim_{n\to\infty} \frac 1n\sigma_n^2=\bE_{\Pe}(\beta).
\]
Thus, $\sigma_n^2$ can grow sublinearly only if $\bE_{\Pe}(\beta)=0$, but then $\bE_{\Pe}(\sigma_n^2)=\cO(1)$.
\end{proof}
Thus, a sublinear growth can happen only if $\bE_{\Pe}(\sigma_n^2)=\cO(1)$.
If so, consider the annealed measure $\bP:=\bP_e\times \mu$ and the related expectation $\bE$. 

The sequence $\hat S_n$ is then uniformly bounded in $L^2(\Omega\times M,\bP)$, hence it is weakly-compact. Let $\psi\in L^2(\Omega\times M,\bP)$ be 
a weak accumulation point of the Ces\`aro averages $A_N=\frac 1N\sum_{n=0}^{ N-1}\hat S_n$.
Then, for each $\vf\in L^2(\Omega\times M,\bP)$, we have
\[
\begin{split}
\bE(\vf\psi)&=\lim_{j\to\infty}\bE\Big(\vf \frac{1}{N_j}\sum_{n=0}^{N_j-1}\sum_{k=0}^{n-1} g_{\omega_k}\circ f_{\omega_{k-1}}\circ\cdots\circ f_{\omega_0} \Big)\\
&=\lim_{j\to\infty} \bE\Big(\cL_{\omega_0}\vf  \frac{1}{N_j}\sum_{n=0}^{N_j-1}\sum_{k=1}^{n-1}g_{\omega_k}\circ f_{\omega_{k-1}}\circ\cdots\circ f_{\omega_1}\Big)+\bE(\vf g_{\omega_0})\\
&=\lim_{j\to\infty} \bE\Big(\cL_{\omega_0}\vf \frac{1}{N_j}\sum_{n=0}^{N_j-1}\sum_{k=0}^{n-1}g_{(\tau\omega)_k}\circ f_{(\tau\omega)_{k-1}}\circ\cdots\circ f_{(\tau\omega)_0}\Big)+\bE(\vf g_{\omega_0})\\
&=\bE\left(\cL_{\omega_0}\vf \cdot \psi\circ\tau\right)+\bE(\vf g_{\omega_0})=\bE(\vf \psi\circ F)+\bE(\vf g_{\omega_0}),
\end{split}
\]
where, in the third equality, we have used that the normalized difference between the two sums is of order $N_j^{-1}$,
since the sequence $\hat S_n$ is bounded in $L^2$.
It follows that $g_{\omega_0}(x)=\psi(\omega,x)- \psi\circ F(\omega,x)$, that is $\psi=\bbg+\psi\circ F$.

Next, let us define $\hat \psi(x)=\bE_{\bP_e}(\psi)$, $\bar g=\bE_{\bP_e}(\bbg)$ and note that
\[
\bar g(x)=\hat\psi-\bE_{\bP_e}(\hat\psi\circ f_{\omega_0}).
\]
Furthermore, define
\[
\bbg\circ F^k=g_{\omega_k}\circ f_{\omega_{k-1}}\circ \cdots\circ f_{\omega_0}(x)=\hat\psi\circ F^{k}-\hat\psi \circ F^{k+1}+Z_k.
\]
Note that
\[
\begin{split}
&\bE_{\bP_e}(Z_k\;|\;\omega_{0},\dots,\omega_{k})=Z_k\\
&\bE_{\bP_e}(Z_k\;|\;\omega_{0},\dots,\omega_{k-1})=0,
\end{split}
\]
where the second identity follows from the independence of the $\omega_k$.
Hence, $\{Z_k\}$ is a martingale difference sequence for the filtration generated by $\omega_0,\dots,\omega_k$, and $M_n=\sum_{k=0}^{n-1}Z_k$ is the associated martingale. It follows that
\[
\bE(\hat S_n^2)=\bE\left([ \hat\psi-\hat\psi\circ F^n +M_n]^2\right)\geq \bE(M_n^2)-2\bE([\hat\psi\circ F^n-\hat \psi]^2)^{\frac12}\bE(M_n^2)^{\frac 12}
\]
Since $\psi\in L^2$ it follows that there exists a constant $C>0$ such that, for all $n\in\bN$,
\[
C\geq \bE(M_n^2)=\sum_{k=0}^{n-1}\bE(Z_k^2)=n\bE(Z_0^2).
\] 
Thus, it must be that
\[
0=\bE(Z_0^2)=\bE([\bbg+\hat\psi\circ F-\hat \psi]^2).
\]
That is, $g_{\fa} -\hat\psi+\hat\psi\circ f_{\fa}=0$ for each $\fa\in\fA$. 
\end{proof}
\begin{rem}\label{rem:L2cob} Note that if the $f_{\fa}$ are smooth Anosov maps and the observables $g_{\fa}$ are smooth, then by Liv\v{s}ic's theory, see \cite{LMM86}, we have that $\hat \psi\in\cC^0$. Hence, the $g_{\fa}$ must be continuous coboundaries, a condition that is checkable. On the contrary, it is not obvious how to check that the $g_{\fa}$ are not $L^2$ coboundaries, so although the present section applies to billiards, the result is still not totally satisfactory. Nevertheless, our result is much sharper than the usual one, e.g. see \cite[Equation (2.10)]{Ki98} where one has an $L^2$ coboundary only with respect to the map \eqref{eq:skew}. A similar result is proven in \cite{ALS09}, but limited to Anosov automorphisms.
\end{rem}

\subsection{Proof of Lemma~\ref{lem:norm est}}\label{sec:lemmaO1}\ \\
To estimate $\| gh \|$, we will find $\cop > 0$ such that $\cop \pm gh \in \cC_\bR$.  To this end, it will be convenient
to estimate the quantity $\frac{\int_W hg  \psi \, dm_W}{\int_W \psi \, dm_W}$ from above and below, where
$W \in \cW^s(\delta)$ and $\psi \in \cD_{a,\beta}(W)$.  For brevity, we will omit $dm_W$ when writing the integral.

For $g \in C^\alpha(M)$ and $\psi \in \cD_{a, \beta}(W)$, we choose
$B > 0$ so that $g \psi + B \in \cD_{a, \beta}(W)$.
Letting $H^\beta(g \psi)$ denote the H\"older constant of $g\psi$ with exponent $\beta$, 
let $B > a^{-1} H^\beta(g \psi ) + |g \psi|_\infty$.  It follows that 
\[
\frac{g(x)\psi(x) +B}{g(y)\psi(y)+B} \le \frac{|g(x)\psi(x) - g(y)\psi(y)|}{g(y)\psi(y)+B} + 1
\le e^{a d(x,y)^\beta}
\]
for all $x,y \in M$, so that in particular, $g\psi+B \in \cD_{a,\beta}(W)$. 
Note that $\psi \in \cD_{a, \beta}(W)$ implies, for $x,y \in W$,
\begin{equation}
\label{eq:psi distortion}
\left| \frac{\psi(x)}{\psi(y)} - 1 \right| \le e^{a d(x,y)^\beta} -1 \le a d(x,y)^\beta  e^{a(2\delta)^\beta }
\le 2 a d(x,y)^\beta \, ,
\end{equation}
so that $H^\beta(\psi) \le 2a |\psi|_{C^0}$, where we have used $e^{a (2\delta)^\beta} \le 2$
by \cite[eq.~(4.9)]{DL22}.
Since $a>1$, we may choose 
\begin{equation}
\label{eq:B choice}
B \le |g|_{C^\beta} |\psi|_{C^\beta} \le |g|_{C^\alpha} (1+2a) |\psi|_{C^0} \, .
\end{equation}

Now, since $|\psi|_{C^0} \le 2 \psi$,
\begin{equation}
\label{eq:upper}
\begin{split}
\frac{\int_W h g \, \psi }{\int_W \psi }
& = \frac{\int_W h (g \psi + B) }{\int_W \psi } - B \frac{\int_W h   }{\int_W \psi } \\
& \le \frac{\int_W h (g \psi + B) }{\int_W g \psi + B } \frac{\int_W g \psi + B }{\int_W \psi } 
 \\
& \le 3(1+2a) | g|_{C^\alpha} \tri h \tri_+ \le 3(1+2a) L |g|_{C^\alpha} \tri h \tri_- \, , 
\end{split}
\end{equation}
where we have dropped the second term since $\tri h \tri_- > 0$ for $h \in \cC_\bR$.
Similarly, we estimate the lower bound, dropping the first term, which is positive,
\begin{equation}
\label{eq:lower}
\begin{split}
\frac{\int_W h g \, \psi }{\int_W \psi }
& = \frac{\int_W h (g \psi + B) }{\int_W \psi } -  \frac{\int_W h B }{\int_W \psi } \\
& \ge - 2(1+2a) |g|_{C^\alpha} \tri h \tri_+ \ge - 2(1+2a) L |g|_{C^\alpha} \tri h \tri_- \, .
\end{split}
\end{equation}

Next we prove that if $\cop$ satisfies \eqref{eq:cop bound}, then $\cop + gh \in \cC_\bR$. 
In order for $\cop + gh$ to satisfy the first cone condition \eqref{eq:cone 1}, we need
\[
\sup_{W, \psi} \frac{\int_W (\cop + gh) \psi}{\int_W \psi }\le L \inf_{W, \psi} \frac{\int_W (\cop + gh) \psi}{\int_W \psi} \, .
\]
This is equivalent to
\[
\cop + \sup_{W, \psi} \frac{ \int_W gh \psi}{\int_W \psi} \le L \Big(\cop + \inf_{W, \psi} \frac{\int_W gh \psi}{\int_W \psi} \Big) \, ,
\]
which according to \eqref{eq:upper} and \eqref{eq:lower} is satisfied if
\[
\cop + 3(1+2a) L |g|_{C^\alpha} \tri h \tri_- \le L (\cop -  2(1+2a) L |g|_{C^\alpha} \tri h \tri_- )\, .
\]
Thus it suffices to choose $\cop$ so that
\begin{equation}
\label{eq:first norm condition}
\cop \ge \frac{3+2L }{L-1} L (1+2a) |g|_{C^\alpha} \tri h \tri_-  \, .
\end{equation}

In order to verify the second cone condition \eqref{eq:cone 2}, $\cop>0$ should satisfy,
\[
|W|^{-q} \frac{ | \int_W (\cop + gh) \psi |}{\fint_W \psi} \le A \delta^{1-q} \tri \cop + gh \tri_- \, ,
\]
for all $W \in \cW^s_-(\delta)$, $\psi \in \cD_{a,\beta}(W)$.

The left hand side is bounded above by 
\begin{equation}
\label{eq:useful g}
\begin{split}
& |W|^{1-q} \cop + |W|^{-q}\frac{|\int_W h (g\psi + B) |}{\fint_W g \psi + B} \frac{\int_W g\psi + B}{\int_W \psi}
+  |W|^{-q} \frac{|\int_W h B |}{\fint_W \psi}  \\
& \le 2^{1-q} \delta^{1-q} \cop + 5(1+2a)|g|_{C^\alpha} A \delta^{1-q} \tri h \tri_- \, ,
\end{split}
\end{equation}
while the right hand side is bounded below using \eqref{eq:lower},
\[
A \delta^{1-q} \tri \cop + gh \tri_- \ge A \delta^{1-q} (\cop - 2(1+2a) L |g|_{C^\alpha} \tri h \tri_- ) \, .
\]
Putting these estimates together, it suffices to choose $\cop>0$ so that
\[
2^{1-q} \delta^{1-q} \cop + 5(1+2a) |g|_{C^\alpha} A \delta^{1-q} \tri h \tri_- \le A \delta^{1-q} (\cop -  2(1+2a) L |g|_{C^\alpha} \tri h \tri_- )
\]
This holds when
\begin{equation}
\label{eq:second norm condition}
\cop \ge \frac{ (5+2L) A   }{A - 2^{1-q}}  (1+2a) |g|_{C^\alpha} \tri h \tri_-  \, .
\end{equation}

Finally, we find $\cop$ to satisfy the third cone condition \eqref{eq:cone 3},
\begin{equation}
\label{eq:third set up}
\left| \frac{\int_{W^1} (\cop + gh) \psi_1}{\fint_{W^1} \psi_1} - \frac{\int_{W^2} (\cop + gh) \psi_2}{\fint_{W^2} \psi_2} \right|
\le d_{\cW^s}(W^1, W^2)^\gamma \delta^{1-\gamma} cA \tri \cop + gh \tri_- \, ,
\end{equation}
for $W^1, W^2 \in \cW^s_-(\delta)$ with $d_{\cW^s}(W^1, W^2) \le \delta$ and
$\psi_i \in \cD_{a,\alpha}(W^i)$ satisfying $d_*(\psi_1, \psi_2) =0$.

Without loss of generality, we may assume $|W^2| \ge |W^1|$ and $\fint_{W^1} \psi_1 = 1$.
Also, we assume  
\begin{equation}
\label{eq:long W2}
|W^2|^q \ge \delta^{q-\gamma} d_{\cW^s}(W^1, W^2)^\gamma \tfrac c2 \, ,
\end{equation}
otherwise, applying the second cone condition to $\cop + gh$ with $\cop$ satisfying \eqref{eq:second norm condition} yields,
\begin{equation}
\label{eq:trivial small}
\begin{split}
\left| \frac{\int_{W^1} (\cop + gh) \psi_1}{\fint_{W^1} \psi_1} - \frac{\int_{W^2} (\cop + gh) \psi_2}{\fint_{W^2} \psi_2} \right|
& \le (|W^1|^q + |W^2|^q) A \delta^{1-q} \tri \cop + gh \tri_- \\
&\hskip-.5cm \le 2 A \delta^{1-\gamma} d_{\cW^s}(W^1, W^2)^\gamma \tfrac c2 \tri \cop + gh \tri_- \, ,
\end{split}
\end{equation}
which is precisely the needed estimate with no additional condition on $\cop$.

We proceed to estimate both sides of \eqref{eq:third set up} under the assumption \eqref{eq:long W2}.
As before, the right-hand side of \eqref{eq:third set up} is bounded below by,
\begin{equation}
\label{eq:lower third}
d_{\cW^s}(W^1, W^2)^\gamma \delta^{1-\gamma} cA (\cop - 2(1+2a)L |g|_{C^\alpha} \tri h \tri_- ) \, .
\end{equation}
To bound the left-hand side of \eqref{eq:third set up} from above, 
we first split up the differences, using that $\fint_{W^1} \psi_1 = 1$, 

\begin{equation}
\label{eq:split third}
\begin{split}
& \quad \; \, \cop \left| |W^1| - |W^2| \right| + \left| \frac{\int_{W^1} hg \psi_1}{\fint_{W^1} \psi_1} - \frac{ \int_{W^2} hg \psi_2}{ \fint_{W^2} \psi_2} \right| \\
& \le \cop \left| |W^1| - |W^2| \right| 
+ \frac{\left| \int_{W^2} hg \psi_2 \right| }{\fint_{W^2} \psi_2 } \left| \fint_{W^2} \psi_2 - 1 \right| 
+ \left| \int_{W^1} hg \psi_1 - \int_{W^2} hg \psi_2 \right| \, .
\end{split}
\end{equation}
The first term above is bounded by \cite[eq.~(5.8)]{DL22},
\begin{equation}
\label{eq:length diff}
\left| |W^1| - |W^2| \right| \le 2C_s d_{\cW^s}(W^1, W^2) \, ,
\end{equation}
where $C_s = \sqrt{1 + (\cK_*^{-1} + \tau_*^{-1})^2}$ denotes the maximum absolute value of the slope of stable curves in $\cW^s_-(\delta)$,
while the difference in the second term above is bounded by \cite[eq.~(5.10)]{DL22},
\begin{equation}
\label{eq:norm almost}
\left| \fint_{W^2} \psi_2 - 1\right| = |W^2|^{-1} \left| \int_{W^2} \psi_2 - |W^2| \right| 
\le \frac{6 C_s d_{\cW^s}(W^1, W^2)}{|W^2|}\, .
\end{equation}
Moreover, by \eqref{eq:useful g},
\[
\frac{\left| \int_{W^2} hg \psi_2 \right|}{\fint_{W^2} \psi_2 } \le 5 (1+2a) |W^2|^q |g|_{C^\alpha} A \delta^{1-q} \tri h \tri_- \, .
\]
Combining these estimates in \eqref{eq:split third} yields the following upper bound for the 
left side of \eqref{eq:third set up}
\begin{equation}
\label{eq:simpler split}
\begin{split}
& \left| \frac{\int_{W^1} (\cop + gh) \psi_1}{\fint_{W^1} \psi_1} - \frac{\int_{W^2} (\cop + gh) \psi_2}{\fint_{W^2} \psi_2} \right|
\le \cop 2 C_s d_{\cW^s}(W^1, W^2) \\
& + \frac{\delta^{1-q} d_{\cW^s}(W^1, W^2)}{|W^2|^{1-q}} 5(1+2a) 6 C_s A |g|_{C^\alpha} \tri h \tri_-
+  \left| \int_{W^1} hg \psi_1 - \int_{W^2} hg \psi_2 \right| \\
& \le \cop 2 C_s d_{\cW^s}(W^1, W^2)
+ 4C_s^q (1+2a) A \delta^{1-\gamma} d_{\cW^s}(W^1, W^2)^\gamma |g|_{C^\alpha} \tri h \tri_- \\
& \; + \left| \int_{W^1} hg \psi_1 - \int_{W^2} hg \psi_2 \right|  \, ,
\end{split}
\end{equation}
where we have used \eqref{eq:long W2} together with $d_{\cW^s}(W^1, W^2) \le \delta$,
$\gamma \le q < 1/2$ and the constraint $c \ge 16 C_s^q$ from \cite[eq.~(5.7)]{DL22}.

To estimate the last term in \eqref{eq:simpler split},
we recall the notation of matched and unmatched pieces.
Recall that each curve $W^i$ is expressed as the graph of a function over an arclength interval $I_i$,
\[
W^i = \{ G_i(r) = (r, \vf_i(r) : r \in I_i \}
\]
Since $d_{\cW^s}(W^1, W^2) \le \delta$,  we have $| I_1 \cap I_2 | > 0$.  
Let
$U^i = \{ G_i(r) : r \in I_1 \cap I_2 \}$ denote the matched pieces and let $V^i$ denote the at most two unmatched pieces.
Remark that $|V^i| \le C_s d_{\cW^s}(W^1, W^2)$.
Choose $B >0$ such that $g \psi_2 + B \in \cD_{\frac{a}{2}, \alpha}(U^2)$.  Following \eqref{eq:B choice}
we may choose $B \le 2 (1+2a) |g|_{C^\alpha} |\psi_2|_{C^0}$.
Define the following functions on $U^1$,
\[
\begin{split}
\tg & = g \circ G_2 \circ G_1^{-1} \cdot \frac{\| G'_2 \|}{\| G'_1 \|} \circ G_1^{-1} \, ; \quad
\tpsi_2 = \psi_2 \circ G_2 \circ G_1^{-1} \\
\tB & = B  \cdot \frac{\| G'_2 \|}{\| G'_1 \|} \circ G_1^{-1}  \, .
\end{split}
\]
Then by definition, $d_*(\tg\tpsi_2 + \tB, g\psi_2 + B) = 0$.
Moreover,
in the proof of \cite[Lemma~5.5(c)]{DL22}, it is shown that for $r,s \in I_1 \cap I_2$, $x = G_k(r)$,
$y = G_k(s)$,
\begin{equation}
\label{eq:graph control}
\frac{\| G'_k(r) \|}{\| G'_k(s) \|} \le e^{B_* d(x,y) } \,\, 
\mbox{ and } \,\,
\frac{\| G'_2(r) \|}{\| G'_1(r) \|} \le e^{d_{\cW^s}(U^1, U^2)} \le e^{\delta} \le 2 \, ,
\end{equation}
for some constant $B_*$ depending only on the maximum curvature of curves in
$\cW^s_-(\delta)$.  This implies in particular that
$\tilde B \le 2 B$.
and that $\tg \tpsi_2 + \tilde B \in \cD_{a, \alpha}(U^2)$
by \cite[eq.~(5.31)]{DL22}, as long as $2B_* (2\delta)^{1-\alpha} \le \frac{a}{2}$, which
is compatible with the restriction on $\delta$ from \cite[eq.~(5.32)]{DL22}.

With these preparations, we are ready to estimate the difference of integrals in \eqref{eq:simpler split}
by splitting into matched and unmatched pieces,
\begin{equation}
\label{eq:final split?}
\begin{split}
\Big| \int_{W^1}& hg\psi_1 - \int_{W^2} hg\psi_2 \Big| 
 \le  \sum_{i=1}^2 \left| \int_{V^i} h g \psi_i \right| + \left| \int_{U^1} hg\psi_1 - \int_{U^2} hg\psi_2 \right| \\
& \le  \sum_{i=1}^2 \left| \int_{V^i} h g \psi_i \right|
+ \left| \int_{U^1} h (g \psi_1 - \tg \tpsi_2) \right| 
 + \left| \int_{U^1} h \tg \tpsi_2 - \int_{U^2} h g \psi_2 \right| \, . 
\end{split}
\end{equation}

The following sublemma allows us to estimate the principal differences in \eqref{eq:final split?}.
\begin{sublem}
\label{sub:main diff}
The following estimates hold true
\begin{itemize}
  \item[a)] $\displaystyle \left| \int_{U^1} h (g \psi_1 - \tg \tpsi_2) \right|  \le 12 (3 + 2a) d_{\cW^s}(W^1, W^2)^{\alpha - \beta} \delta A |g|_{C^\alpha} \tri h \tri_-$.
  \item[b)] $\displaystyle \left| \int_{U^1} h \tg \tpsi_2 - \int_{U^2} h g \psi_2 \right| 
  \le 20(1+2a) (c+3\delta) A d_{\cW^s}(W^1, W^2)^\gamma \delta^{1-\gamma} |g|_{C^\alpha} \tri h \tri_-$.
\end{itemize}
\end{sublem}

Postponing the proof of the sublemma, we use it to complete the estimate for the third cone condition.
Using the sublemma, together with \eqref{eq:useful g} to bound the integrals on unmatched pieces, 
and recalling that $|V^i| \le C_s d_{\cW^s}(W^1, W^2)$, we 
bound \eqref{eq:final split?} by
\[
\begin{split}
\left| \int_{W^1} hg\psi_1 \right.  & - \left. \int_{W^2} hg\psi_2 \right| 
 \le 40 (1+2a) d_{\cW^s}(W^1, W^2)^q C_s^q A \delta^{1-q} |g|_{C^\alpha} \tri h \tri_- \\
& + 12 (3 + 2a) d_{\cW^s}(W^1, W^2)^{\alpha - \beta} \delta A |g|_{C^\alpha} \tri h \tri_- \\
& + 20(1+2a) (c+3\delta) A d_{\cW^s}(W^1, W^2)^\gamma \delta^{1-\gamma} |g|_{C^\alpha} \tri h \tri_- \\
& \le 4(1+2a)(10 C_s^q + 9\delta^\gamma + 5c + 15 \delta) d_{\cW^s}(W^1, W^2)^\gamma \delta^{1-\gamma} A |g|_{C^\alpha} \tri h \tri_- \, ,
\end{split} 
\] 
where we have used that $d_{\cW^s}(W^1, W^2) \le \delta$ and $\gamma \le \min \{ q, \alpha - \beta \}$.
This estimate combined with \eqref{eq:simpler split} yields our final upper bound for the left hand side of
\eqref{eq:third set up},
\[
\begin{split}
\left| \frac{\int_{W^1} (\cop + gh) \psi_1}{\fint_{W^1} \psi_1} - \frac{\int_{W^2} (\cop + gh) \psi_2}{\fint_{W^2} \psi_2} \right|
\le& 2 \cop C_s d_{\cW^s}(W^1, W^2)^\gamma \delta^{1-\gamma} \\
&+ C_\star d_{\cW^s}(W^1, W^2)^\gamma \delta^{1-\gamma} |g|_{C^\alpha} \tri h \tri_- \, ,
\end{split}\]  
where $C_\star = 4A(1+2a)(14 C_s^q + 9\delta^\gamma + 5c + 15 \delta)$.
Combining this with the lower bound \eqref{eq:lower third}, we see that \eqref{eq:third set up} will be satisfied
if
\[
\cop 2 C_s + C_\star |g|_{C^\alpha} \tri h \tri_- 
\le cA (\cop - 2(1+2a)L |g|_{C^\alpha} \tri h \tri_- ) \, .
\]
Since $cA > 2C_s$ by \cite[eq.~(5.36)]{DL22}, this in turn is satisfied if
\begin{equation}
\label{eq:third norm condition}
\cop \ge \frac{C_\star + 2cA(1+2a)L}{cA - 2C_s} |g |_{C^\alpha} \tri h \tri_- \, .
\end{equation}
Now taking together \eqref{eq:first norm condition}, \eqref{eq:second norm condition} and
\eqref{eq:third norm condition}, we conclude that $\cop + gh \in \cC_\bR$ if
$\cop$ satisfies \eqref{eq:cop bound}.

We claim that the same value of $\cop$ implies that $\cop - gh \in \cC_\bR$.  To see this, write
$\cop - gh = \cop + (-g)h$, and notice that $-g \in C^\alpha$ with $|-g|_{C^\alpha} = |g|_{C^\alpha}$.
Applying the previous argument to $-gh$, we see that $\cop + (-g)h \in \cC_\bR$ when
$\cop$ satisfies \eqref{eq:cop bound}.  This completes the proof of the lemma.

It remains to prove the sublemma used during the argument.

\begin{proof}[Proof of Sublemma~\ref{sub:main diff}]
We prove the two statements one at a time:\\
a) We would like to apply the second cone
condition \eqref{eq:cone 2} to the integral;  however, $g \psi_1 - \tilde g \tpsi_2$ may not be in $\cD_{a,\beta}(U^1)$.  To remedy this, first note that
for $x = G_1(r) \in U^1$, using the fact that $d_*(\psi_1, \psi_2)=0$,
\[
\begin{split}
|g(x)\psi_1(x) -\tilde g \tpsi_2(x)| & = \tfrac{1}{\| G'_1 \|} |g\psi_1(G_1(r)) \cdot \| G'_1(r) \| - g\psi_2(G_2(r)) \cdot \| G'_2 \| | \\
& \le |\psi_1(G_1(r))|  |g(G_1(r)) - g(G_2(r))| \\
& \le 2 H^\alpha(g) d_{\cW^s}(W^1, W^2)^\alpha \, ,
\end{split}
\]
where we have used  the bound $|\psi_1| \le e^{a(2\delta)^\alpha} \fint_{W^1} \psi_1 \le 2$.
 
So on the one hand,
\[
|(g \psi_1 - \tilde g \tpsi_2 )(x) - (g \psi_1 - \tilde g \tpsi_2)(y)| \le 4 H^\alpha(g) d_{\cW^s}(W^1, W^2)^\alpha \, , \mbox{ for } x, y \in U^1 \, .
\]
On the other hand, for $x = G_1(r) \in U^1$, let $\tilde x = G_2(r) \in U^2$ denote the corresponding point in $U^2$. Then,
\[
\begin{split}
d(\tilde x, \tilde y) =& \int_{r_1}^{r_2} \sqrt{1 + (d\vf_2/dr)^2} \, dr \le \sup_{I_1 \cap I_2} \sqrt{ \frac{1 + (d\vf_2/dr)^2}{1 + (d\vf_1/dr)^2} } 
\int_{r_1}^{r_2} \sqrt{1 + (d\vf_1/dr)^2 } \, dr\\
& \le 2 d(x,y) \, ,
\end{split}
\]
where we have used the estimate before \cite[eq. (5.9)]{DL22} to estimate the ratio of Jacobians.  Then
recalling \eqref{eq:B choice},
\[
\begin{split}
|(g\psi_1 - \tilde g \tpsi_2)(x) - (g \psi_1 - \tilde g \tpsi_2)(y)| 
& \le |g \psi_1(x) - g \psi_1(y)| + |g (\tilde x) \psi_1(x) - g (\tilde y) \psi_1(y)| \\
&  \le |g|_{C^\alpha} |\psi_1|_{C^\alpha} (d(x,y)^\alpha + d(\tilde x, \tilde y)^\alpha) \\
& \le 6 (1+2a) |g|_{C^\alpha} d(x,y)^\alpha \, .
\end{split}
\]
Putting these estimates together, we see that
\[
\frac{|(g\psi_1 - \tilde g \tpsi_2)(x) - (g\psi_1 - \tilde g \tpsi_2)(y)| }{d(x,y)^\beta} \le \frac{|g|_{C^\alpha}}{d(x,y)^\beta} \min \{ 4 d_{\cW^s}(W^1, W^2)^\alpha , 6(1+2a) d(x,y)^{\alpha} \} \, ,
\]
and the expression is maximized when the two quantities are equal, i.e.~when $d_{\cW^s}(W^1, W^2) = (3(1+2a)/2)^{1/\alpha} d(x,y)$.  Thus,
\[
H^\beta_{U^1}(g\psi_1 - \tilde g \tpsi_2) \le 6(1+2a) |g|_{C^\alpha} d_{\cW^s}(W^1, W^2)^{\alpha-\beta}  \, .
\]
So we choose a constant $\Delta \ge a^{-1} H^\beta_{U^1}(g\psi_1 - \tilde g \tpsi_2) + |g \psi_1 - \tilde g \tpsi_2|_{C^0(U^1)}$, so that
$g \psi_1 - \tilde g \tpsi_2+ \Delta \in \cD_{a,\beta}(U^1)$, and $\Delta \le |g \psi_1 - \tilde g \tpsi_2|_{C^\beta(U^1)} \le (8+6a) |g|_{C^\alpha} d_{\cW^s}(W^1, W^2)^{\alpha - \beta}$.

Now applying \eqref{eq:cone 2} completes the proof of statement (a) of the sublemma,
\[
\begin{split}
\left| \int_{U^1} h (g \psi_1 - \tg \tpsi_2) \right|
& \le \left| \int_{U^1} h (g \psi_1 - \tg \tpsi_2 + \Delta) \right| + \left| \int_{U^1} h \Delta \right| \\
& \le |U^1|^q \fint_{U^1} (g \psi_1 - \tg \tpsi_2 + 2\Delta) A \delta^{1-q} \tri h \tri_- \\
& \le 2 (18 + 12a) d_{\cW^s}(W^1, W^2)^{\alpha - \beta} \delta A |g|_{C^\alpha} \tri h \tri_-\, .
 \end{split}
\]

\medskip
\noindent
b) With $B$ and $\tilde B$ as chosen before the statement of the sublemma, write,
\begin{equation}
\label{eq:matched split}
\begin{split}
 \left| \int_{U^1}\!\!\! h \tg \tpsi_2 - \int_{U^2}\!\!\!  h g \psi_2 \right| &
\le \left| \frac{\int_{U^1} h (\tg \tpsi_2 + \tB)}{\fint_{U^1} \tg \tpsi_2 + \tB} - \frac{ \int_{U^2} h(g \psi_2+B)}{\fint_{U^2} g\psi_2 +B} \right| \fint_{U^2}\!\! (g\psi_2 + B) \\
& + \left| \int_{U^1} h(\tg \tpsi_2 + \tB) \right| \left| 1 - \frac{|U^1|}{|U^2|} \right|
+ \left| \int_{U^1} h \tB - \int_{U^2} h B \right|.
\end{split}
\end{equation}
The first term of \eqref{eq:matched split} can be estimated using \eqref{eq:cone 3} since
$d_*(\tg \tpsi_2 + \tB, g \psi_2 + B) =0$ and both are valid test functions by choice of $B$.  For the second
term of \eqref{eq:matched split}, we use \cite[eq.~(5.24)]{DL22} to estimate
\begin{equation}
\label{eq:close piece}
\left| 1 - \frac{|U^2|}{|U^1|} \right| \le d_{\cW^s}(U^1, U^2) \le d_{\cW^s}(W^1, W^2) \, .
\end{equation}
Using this together with \eqref{eq:cone 2} yields the bound,
\[
\begin{split}
& \left| \frac{\int_{U^1} h (\tg \tpsi_2 + \tB)}{\fint_{U^1} \tg \tpsi_2 + \tB} - \frac{ \int_{U^2} h(g \psi_2+B)}{\fint_{U^2} g\psi_2 +B} \right| \fint_{U^2} (g\psi_2 + B) 
+ \left| \int_{U^1} h(\tg \tpsi_2 + \tB) \right| \left| 1 - \frac{|U^1|}{|U^2|} \right| \\
& \le cA d_{\cW^s}(W^1, W^2)^\gamma \delta^{1-\gamma} \tri h \tri_- 12 (1+2a) |g|_{C^\alpha}
+ 36(1+2a) A \delta |g|_{C^\alpha} \tri h \tri_- d_{\cW^s}(W^1, W^2) \, . 
\end{split}
\]
Moreover, it is clear that the final term in \eqref{eq:matched split} can be estimated in an analogous way, with
$B$ and $\tilde B$ replacing the test functions $g\psi_2 + B$ and $\tg \tpsi_2 + \tB$.  Thus,
\[
\begin{split}
\left| \int_{U^1} h \tB - \int_{U^2} h B \right| 
& \le B cA d_{\cW^s}(W^1, W^2)^\gamma \delta^{1-\gamma} \tri h \tri_- + 3B d_{\cW^s}(W^1, W^2) \delta A \tri h \tri_- \\
& \le 8(1+2a) A d_{\cW^s}(W^1, W^2)^\gamma \delta^{1-\gamma} |g|_{C^\alpha} \tri h \tri_- (3\delta + c) \, .
\end{split}
\]
Putting this estimate together with the previous one completes the proof of statement (b), using that
$d_{\cW^s}(W^1, W^2) \le \delta$.
\end{proof}


\subsection{Proof of Proposition~\ref{lem:bounded L} (Bounding the Norm of $\cL_T$) }\ \\
\label{sec:proof of bounded L}

As mentioned in Remark \ref{rem:one_step}, the goal of this section is to obtain a bound on $\|\cL_T\|$, that is, for the action of a single operator rather than a sequence of $N_\cF$ operators. The proof is similar to the proof of Lemma \ref{lem:norm est} in Section~\ref{sec:lemmaO1}. Unfortunately, the details are different enough that we need to carry it out explicitly.

We will first bound $\cL h$ for $h \in \cC_{\bR}$ and then show how this 
extends to all $h \in \cB_{\bR}$ as in the proof of Corollary~\ref{cor:in space}.

Let $h \in \cC_{\bR}$ and $T \in \cF(\tau_*, \cK_*, E_*)$ with $\cL = \cL_T$.  We must estimate 
the norm of $\cL h$.  To this end, we first obtain bounds on $\tri \cL h \tri_+$ and $\tri \cL h \tri _-$.
Let $W \in \cW^s(\delta)$ and $\psi \in \cD_{a,\beta}(W)$.  Then,
\begin{equation}
\label{eq:first change}
\int_W \cL h \, \psi  = \sum_{i \in L_1(W)} \int_{W_i} h \, \psi \circ T \, J_{W_i}T 
+ \sum_{i \in S_1(W)} \int_{W_i} h \, \psi \circ T \, J_{W_i}T  \, ,
\end{equation}
where we have changed variables and denote by $L_1(W)$ the connected homogeneous components
of $T^{-1}W = \{ W_i \}_i$ longer than $\delta$, and by $S_1(W)$ those homogeneous components shorter than 
$\delta$.

Following \eqref{eq:B choice}, in order to transform the $\psi \circ T J_{W_i}T$ into valid test functions, we choose $B_i$ such that
\begin{equation}
\label{eq:Bi}
\begin{split}
a^{-1} H^\beta(\psi \circ T J_{W_i}T) +& |\psi \circ T  J_{W_i}T |_{C^0(W_i)}  < B_i 
 \le |J_{W_i}T \psi \circ T|_{C^\beta(W_i)} \\
& \le (1+2aC_1 +C_d) |J_{W_i}T|_{C^0(W_i)} |\psi \circ T|_{C^0(W_i)} \, ,
\end{split}
\end{equation}
where $C_d > 0$ is a distortion constant for $J_{W^i}T$, uniform for $T \in \cF(\tau_*, \cK_*, E_*)$
as in \cite[eq. (3.5)]{DL22}, and we have used \eqref{eq:psi distortion} to bound the distortion of $\psi \circ T$
together with $d(Tx, Ty) \le C_1 d(x,y)$, where $C_1 \ge 1$ is the minimum hyperbolicity 
constant\footnote{That is, setting $x_0 = 1 + 2 \cK_* \tau_*$, 
$C_1 \ge 1$ satisfies $\| DT^{-n}(x) v \| \ge C_1^{-1} \Lambda_0^n \| v\|$ for
all $T \in \cF(\tau_*, \cK_*, E_*)$ and $v$ in the stable cone. \label{foot:C1}}
from \cite[eq. (3.1)]{DL22}.
Then as in Section~\ref{sec:lemmaO1}, $J_{W_i}T \psi \circ T + B_i \in \cD_{a, \beta}(W)$.
For ease of notation, let us name the combined distortion constants,
\[
C_a := 1 + 2aC_1 + C_d \, .
\]

Now we estimate the contribution from long pieces in \eqref{eq:first change} using \eqref{eq:tri def} and \eqref{eq:cone 1}, as well as the upper bound on $B_i$ from \eqref{eq:Bi},
\begin{equation}
\label{eq:long first}
\begin{split}
\sum_{i \in L_1(W)} \int_{W_i} h \,& (\psi \circ T \, J_{W_i}T  + B_i) - B_i \sum_{i \in L_1(W)} \int_{W_i} h 
\\
& \le \tri h \tri_+ \sum_{i \in L_1(W)} \int_{W_i}  (\psi \circ T \, J_{W_i} T + 2B_i) \\
& \le \big[ 1+2C_a e^{C_d (2\delta)^{1/3} +a (2\delta)^\beta} \big] L \tri h \tri_-  \sum_{i \in L_1(W)} \int_{TW_i} \psi\\
&  \le \big[ 1 + 8C_a \big] L \tri h \tri_- \int_W \psi \, ,
\end{split}
\end{equation}
recalling that by \cite[Sect.~5.3]{DL22}, $\delta$ is small enough that 
$e^{a(2\delta)^\beta}$ and $e^{C_d (2\delta)^{1/3}}$ are each less than 2.

Next we estimate the contribution from short pieces in \eqref{eq:first change}, using \eqref{eq:cone 2},
\begin{equation}
\label{eq:short first}
\begin{split}
\sum_{i \in S_1(W)} & \int_{W_i} h \, (\psi \circ T \, J_{W_i}T  + B_i) - B_i \sum_{i \in S_1(W)} \int_{W_i} h 
\\
& \le \tri h \tri_- \sum_{i \in S_1(W)} A \delta^{1-q} |W_i|^q \fint_{W_i}  (\psi \circ T \, J_{W_i} T + 2B_i) \\
& \le \tri h \tri_- A \delta \big[ 1+2 C_a) \big] |\psi|_{C^0(W)} \sum_{i \in S_1(W)}  |J_{W_i}T|_{C^0(W_i)} . 
\end{split}
\end{equation}
Since $\psi \in \cD_{a, \beta}(W)$ and $|W| \ge \delta$, we estimate
$\delta |\psi|_{C^0(W)} \le \delta e^{a|W|^\beta} \fint_W \psi \le 2 \int_W \psi$, recalling again
that $e^{a(2\delta)^\beta} \le 2$.  
The sum over Jacobians is bounded by $C_0 \theta_0$, which is the uniform
one-step expansion estimate stated in \cite[eq. (3.4)]{DL22}, where $C_0 \ge 1$ and $\theta_0 < 1$.
Putting these bounds together yields,
\begin{equation}
\label{eq:short second}
\sum_{i \in S_1(W)} \int_{W_i} h \, \psi \circ T \, J_{W_i}T \le 
\tri h \tri_-  A \big[ 1+2 C_a \big] 2 C_0 \theta_0 \int_W \psi \, .
\end{equation}
Substituting \eqref{eq:long first} and \eqref{eq:short second} into \eqref{eq:first change}
and taking the supremum over $W \in \cW^s(\delta)$ and $\psi \in \cD_{a, \beta}(W)$ yields,
\begin{equation}
\label{eq:upper tri}
\tri \cL h \tri_+ \le \tri h \tri_- \Big( L \big[ 1 + 8 C_a \big]
+ A \big[ 1+2 C_a\big] 2 C_0 \theta_0 \Big) \, .
\end{equation}

Next, we need the analogous lower bound for $\tri \cL h \tri_-$.  Taking $W \in \cW^s(\delta)$ 
and $\psi \in \cD_{a, \beta}(W)$, we first change variables as in \eqref{eq:first change}.
For the estimate on long pieces, we use the fact that integrals of $h$ against valid test functions
on curves of length at least $\delta$ are positive by the first cone condition, \eqref{eq:cone 1}.
Then, using \eqref{eq:Bi} and the distortion bounds as before,
\begin{equation}
\label{eq:long lower}
\begin{split}
\sum_{i \in L_1(W)} &\int_{W_i} h \,  (\psi \circ T \, J_{W_i}T + B_i) - B_i \sum_{i \in L_1(W)} \int_{W_i} h 
\\
& \ge \sum_{i \in L_1(W)} \tri h \tri_- \int_{W_i}  (\psi \circ T \, J_{W_i} T + B_i) 
- \tri h \tri_+ \int_{W_i} B_i \\
& \ge \sum_{i \in L_1(W)} \tri h \tri_- 2 \int_{TW_i} \psi - L \tri h \tri_- C_a e^{C_d (2\delta)^{1/3} +a (2\delta)^\beta} \int_{TW_i} \psi \\
&  \ge \tri h \tri_- \big[ 2 - 4C_0\theta_0 -  4L C_a \big]  \int_W \psi \, .
\end{split}
\end{equation}
In the second line we used that, by \eqref{eq:Bi}, $\int_{W_i}B_i\geq \int_{W_i} |\psi \circ T \, J_{W_i}T|_{C^0(W_i)}\geq \int_{TW_i}\psi$, whence the factor $2$.  In the last line we used that $\{TW_i\}_{i\in L_1(W)\cup S_1(W)}$ is a partition of $W$ and that, by the estimates leading to \eqref{eq:short second},
\[
\sum_{i \in S_1(W)}\int_{TW_i}\psi\leq \delta|\psi|_{C^0(W)}\sum_{i\in S_1(W)}|J_{W_i}T|_{C^0(W_i)}\leq 2C_0\theta_0\int_W\psi\, ,
\]
so that $\sum_{i \in L_1(W)}\int_{TW_i}\psi\geq (1-2C_0\theta_0)\int_W\psi$.

The contribution from short pieces is estimated using \eqref{eq:cone 2} as in \eqref{eq:short second},
but with the observation that the integrals of $h$ on short pieces can be negative,
\begin{equation}
\label{eq:short lower}
\begin{split}
\sum_{i \in S_1(W)} & \int_{W_i} h \, \psi \circ T \, J_{W_i}T  + B_i) - B_i \sum_{i \in S_1(W)} \int_{W_i} h 
\\
& \ge - \tri h \tri_- \sum_{i \in S_1(W)} A \delta^{1-q} |W_i|^q \fint_{W_i}  (\psi \circ T \, J_{W_i} T + 2B_i) 
\\
& \ge - \tri h \tri_-  A \big[ 1+2 C_a\big] 2 C_0 \theta_0 \int_W \psi \, .
\end{split}
\end{equation}
Putting together \eqref{eq:long lower} and \eqref{eq:short lower} and taking the appropriate infima yields,
\begin{equation}
\label{eq:lower tri}
\tri \cL h \tri_- \ge \tri h \tri_- \Big( 2 - 4C_0\theta_0 - 4LC_a - A \big[ 1+2 C_a \big] 2 C_0 \theta_0 \Big) \, .
\end{equation}

We are now in a position to choose $\cop>0$ so that $\cop + \cL h$ satisfies the
first cone condition \eqref{eq:cone 1}.  To simplify notation, define
$H_0 = L \big[ 1 + 8 C_a\big]
+ A \big[ 1+2 C_a \big] 2 C_0 \theta_0$ and
$H_1 = -2 + 4C_0\theta_0 + 4L C_a + A \big[ 1+2 C_a \big] 2 C_0 \theta_0 >0$.
Then using \eqref{eq:upper tri} and \eqref{eq:lower tri},
\[
\frac{\tri \cop + \cL h \tri_+}{\tri \cop + \cL h \tri_-}
\le L \quad \impliedby \quad \frac{\cop + H_0 \tri h \tri_-}{\cop - H_1 \tri h \tri_-} \le L \, .
\]
This later inequality is satisfied as soon as
\begin{equation}
\label{eq:cone 1 check}
\cop \ge \tri h \tri_- \frac{H_0 + L H_1}{L-1} \, .
\end{equation}

To guarantee that $\cop + \cL h$ satisfies the second cone condition \eqref{eq:cone 2}, we need
\begin{equation}
\label{eq:cone 2 prelim}
|W|^{-q} \frac{|\int_W (\cop + \cL h) \psi |}{\fint_W \psi} \le A \delta^{1-q} \tri \cop + \cL h \tri _-\, .
\end{equation}
for all $W \in \cW^s_-(\delta)$ and $\psi \in \cD_{a, \beta}(W)$.
Using \eqref{eq:long first} and \eqref{eq:short first}, the left hand side is bounded above by
\begin{equation}
\label{eq:cone 2 first}
|W|^{1-q} \cop + |W|^{-q} \frac{| \int_W \cL h \psi |}{\fint_W \psi} \, .
\end{equation}
Changing variables as in \eqref{eq:first change}, we estimate the contribution from long pieces
precisely as in \eqref{eq:long first}.  Yet the contribution from short pieces must be done with 
care since $W$ itself may be short.  As in \eqref{eq:short first}, we estimate,
\[
\begin{split}
\sum_{i \in S_1(W)} & \int_{W_i} h \, (\psi \circ T \, J_{W_i}T  + B_i) - B_i \sum_{i \in S_1(W)} \int_{W_i} h 
\\
& \le \tri h \tri_- \sum_{i \in S_1(W)} A \delta^{1-q} |W_i|^q \fint_{W_i}  (\psi \circ T \, J_{W_i} T + 2B_i) \\
& \le \tri h \tri_- A \delta^{1-q} \big[ 1 + 2C_a \big]  
\sum_{i \in S_1(W)} |W_i|^q |\psi|_{C^0} |J_{W_i}T|_{C^0(W_i)} \, .
\end{split}
\]
Now $| \psi|_{C^0} \le e^{a(2\delta)^\beta} \fint_W \psi \le 2 \fint_W \psi$ so that,
\begin{equation}
\label{eq:short short}
\begin{split}
\sum_{i \in S_1(W)}\!\!\!\!  \frac{|\int_{W_i} h \psi \circ T J_{W_i}T|}{|W|^{q}\fint_W \psi}
&\le 2 \tri h \tri_- A  \frac{\big[ 1 + 2 C_a \big] }{ \delta^{q-1}}\hskip-6pt
\sum_{i \in S_1(W)} \frac{|W_i|^q}{|W|^q} |J_{W_i}T|_{C^0(W_i)} \\
& \hskip -1cm\le \tri h \tri_- 4 A \delta^{1-q} \big[ 1 + 2 C_a \big]  
\sum_{i \in S_1(W)} \frac{|TW_i|^{q}}{|W|^{q}} \frac{|TW_i|^{1-q}}{|W_i|^{1-q}} \\
& \hskip -1cm \le  \tri h \tri_- 4A \delta^{1-q} \big[ 1 + 2 C_a \big]  
\left( \sum_{i \in S_1(W)} \frac{|TW_i|}{|W_i|} \right)^{1-q} \\
& \hskip -1cm\le  \tri h \tri_- 4A \delta^{1-q} \big[ 1 + 2 C_a \big] \left( C_0 \theta_0 \right)^{1-q}
\end{split}
\end{equation}
where in the third line we have used the H\"older inequality together with the fact that
$\sum_i \frac{|TW_i|}{|W|} \le 1$, and in the last line we have again used the one-step
expansion from \cite[eq.~(3.4)]{DL22}.

Putting this estimate together with the estimate on long pieces yields an upper bound
for the left side of \eqref{eq:cone 2 prelim} (using the fact that $|W| \le 2\delta$),
\begin{equation}
\label{eq:upper 2}
\begin{split}
(2\delta)^{1-q}\cop & + \tri h \tri _- \Big( (2\delta)^{1-q} \big[1 + 8 C_a \big] L 
+ 4 A \delta^{1-q} \big[ 1 + 2 C_a \big] (C_0 \theta_0)^{1-q} \Big) \\
& =: 2^{1-q} \delta^{1-q} \cop + \delta^{1-q} \tri h \tri_- H_2 \, .
\end{split}
\end{equation}
On the other hand, we obtain a lower bound on the right side of \eqref{eq:cone 2 prelim}
using \eqref{eq:lower tri}
\[
A \delta^{1-q} \tri \cop + \cL h \tri_- \ge A \delta^{1-q} ( \cop - H_1 \tri h \tri_- ) \, .
\]
This together with \eqref{eq:upper 2} implies that \eqref{eq:cone 2 prelim} is satisfied provided
\[
2^{1-q} \delta^{1-q} \cop + \delta^{1-q} \tri h \tri_- H_2 \le  A \delta^{1-q} ( \cop - H_1 \tri h \tri_- ) \, .
\]
The above holds true provided,
\begin{equation}
\label{eq:cone 2 check}
\cop \ge \tri h \tri_- \frac{H_2 + A H_1}{A - 2^{1-q}} \, .
\end{equation}

\begin{rem}
\label{rem:useful bound}
Although $\cL h$ may not be in the cone, the estimates \eqref{eq:cone 2 first} and \eqref{eq:upper 2}
together prove that for all $W \in \cW^s_-(\delta)$ and all $\psi \in \cD_{a, \beta}(W)$,
\[
\left| \int_W \cL h \, \psi \right| \le \fint_W \psi |W|^q \delta^{1-q} H_2 \tri h \tri_- \, ,
\]
which will be useful in what follows.
\end{rem}
It remains to choose $\cop$ large enough to satisfy the third cone condition \eqref{eq:cone 3}.
For this we need,
\begin{equation}
\label{eq:cone 3 prelim}
\left|\frac{\int_{W^1} (\cop + \cL h ) \psi_1}{\fint_{W^1}\psi_1}  - \frac{\int_{W^2} (\cop+  \cL h) \psi_2}{\fint_{W^2}\psi_2} \right|
\leq
d_{\cW^s}(W^1, W^2)^\gamma \, \delta^{1-\gamma}   c A \tri \cop + \cL h \tri_- \, ,
\end{equation}
for $W^j \in \cW^s_-(\delta)$ and $\psi_j \in \cD_{a, \alpha}(W^j)$ with $d_*(\psi_1, \psi_2) = 0$.
As before, using \eqref{eq:lower tri} 
the right hand side of \eqref{eq:cone 3 prelim} is bounded below by
\begin{equation}
\label{eq:lower 3}
d_{\cW^s}(W^1, W^2)^\gamma \delta^{1-\gamma} cA ( \cop - H_1 \tri h \tri_- ) \, .
\end{equation}
We proceed to obtain an upper bound for the left hand side of \eqref{eq:cone 3 prelim}. 
Without loss of generality, we may assume that $\fint_{W^1} \psi_1 = 1$ and $|W^2| \ge |W^1|$.
Also, we assume
\begin{equation}
\label{eq:lower W2}
|W^2|^q \ge \tfrac 14 c \delta^{q-\gamma} d_{\cW^s}(W^1, W^2)^\gamma \, .
\end{equation}
Otherwise, applying \eqref{eq:cone 2} to both terms separately and using \eqref{eq:upper tri}, we obtain,
\[
\begin{split}
\left|\frac{\int_{W^1} (\cop + \cL h ) \psi_1}{\fint_{W^1}\psi_1}  - \frac{\int_{W^2} (\cop+  \cL h) \psi_2}{\fint_{W^2}\psi_2} \right|
& \leq 2 |W^2|^q A \delta^{1-q} \tri \cop + \cL h \tri_- \\
& \leq \tfrac{1}{2} c A d_{\cW^s}(W^1, W^2)^\gamma \delta^{1-\gamma} ( \cop + H_0 \tri h \tri_-) \, .
\end{split}
\]
This, together with \eqref{eq:lower 3} implies that \eqref{eq:cone 3 prelim} holds whenever 
\begin{equation}\label{eq:no7.60}
\cop \ge \tri h \tri_- (H_0 + 2 H_1)
\end{equation}
 
We proceed to prove the upper bound under the assumption \eqref{eq:lower W2}.  Now,
\[
\begin{split}
& \left|\frac{\int_{W^1} (\cop + \cL h ) \psi_1}{\fint_{W^1}\psi_1}  - \frac{\int_{W^2} (\cop+  \cL h) \psi_2}{\fint_{W^2}\psi_2} \right|
\leq \left| |W^1| - |W^2| \right| \cop + \left|\frac{\int_{W^1} \cL h \, \psi_1}{\fint_{W^1}\psi_1}  - \frac{\int_{W^2} \cL h \, \psi_2}{\fint_{W^2}\psi_2} \right| \\
& \quad \leq \left| |W^1| - |W^2| \right| \cop
 + \frac{\left| \int_{W^2} \cL h \, \psi_2 \right|}{\fint_{W^2} \psi_2}
\left| \fint_{W^2} \psi_2 -1 \right| 
+ \left| \int_{W^1} \cL h \, \psi_1 - \int_{W^2} \cL h \, \psi_2 \right| \, .
\end{split}
\]
Using \eqref{eq:length diff}, \eqref{eq:norm almost} and 
Remark~\ref{rem:useful bound}, we can write
\begin{equation}
\label{eq:focus}
\begin{split}
& \left|\frac{\int_{W^1} (\cop + \cL h ) \psi_1}{\fint_{W^1}\psi_1}   -  \frac{\int_{W^2} (\cop+  \cL h) \psi_2}{\fint_{W^2}\psi_2} \right|
 \le \left| \int_{W^1} \cL h \, \psi_1 - \int_{W^2} \cL h \, \psi_2 \right| \\
& \hskip24pt+ \cop 2 C_s d_{\cW^s}(W^1, W^2) +\left( \frac{\delta}{|W^2|}\right)^{1-q} 6 C_s H_2  \tri h \tri_-  d_{\cW^s}(W^1, W^2) \\
& \hskip12pt \le \left| \int_{W^1} \cL h \, \psi_1 - \int_{W^2} \cL h \, \psi_2 \right|\\
& \hskip24pt+ \delta^{1-\gamma} d_{\cW^s}(W^1, W^2)^\gamma \left( 2 C_s^q H_2 \tri h \tri_- + 2\cop C_s \right) \, ,
\end{split}
\end{equation}
where in the last line we have used \eqref{eq:lower W2} together with $d_{\cW^s}(W^1, W^2) \le \delta$,
$\gamma \le q < 1/2$ and the constraint $c \ge 16 C_s^q$ from \eqref{eq:const_cond} (see \cite[eq.~(5.7)]{DL22}for more details).

It remains to estimate the difference in integrals in \eqref{eq:focus}.  For this, we change variables as usual
and integrate on elements of $\cG_1^\delta(W^k) = \{ W^k_i \}_i$, the homogeneous connected components
of $T^{-1}W^k$, with long pieces subdivided to have length between $\delta$ and $2\delta$.  As in 
Section~\ref{sec:lemmaO1} and following \cite[Sect.~5.2.3]{DL22}, we subdivide elements of
$\cG_1(W^k)$ into matched and unmatched pieces by defining a foliation of vertical line segments 
$\{ \ell_x \}_{x \in W^1_i}$ centered at $x$ of length at most $3 C_1 d_{\cW^s}(W^1, W^2)$
such that their images under $T$ either terminate on a singularity curve in $\cS_{-1}^{\bH}$ or else
are unstable curves having length at least 
$d_{\cW^s}(W^1, W^2)$ on either side of $T(x) \in W^1$.  In the latter
case, either $T(\ell_x)$ intersects $W^2$, or by the uniform transversality of stable and unstable
curves, lies within distance $C_2 d_{\cW^s}(W^1, W^2)$ of an endpoint of $W^1$, where 
$C_2$ depends only on the minimum angle between stable and unstable curves.
When $T(\ell_x)$ intersects $W^2$, then necessarily $\ell_x$ intersects an element of 
$\cG_1^\delta(W^2)$.  We call subcurves for which all points are connected by such vertical line segments `matched.'  The rest we call `unmatched.' 

With this identification, we may label elements of $\cG_1^\delta(W^1)$ and $\cG_1^\delta(W^2)$
so that each element of $\cG_1^\delta(W^1)$ contains at most one matched subcurve and
at most 2 unmatched subcurves.  We will use the decomposition
$\cG_1^\delta(W^k) =( \cup_j U^k_j )\cup (\cup_j V_j^k)$, so that $U^1_j$ and $U^2_j$ are matched
and so recalling \eqref{eq:graph}, are defined as the graphs of functions $G_{U^k_j}$ over the same
$r$-interval $I_j$ for each $j$.  Using this decomposition, we write
\[
\int_{W^k} \cL h \, \psi_k = \sum_j \int_{U^k_j} h \, \psi_k \circ T \, J_{U^k_j} T 
+ \sum_j \int_{V^k_j} h \, \psi_k \circ T \, J_{V^k_j}T \, .
\]
We estimate the contribution to \eqref{eq:focus} from unmatched pieces first.

As noted previously, an unmatched curve in $V^k_j$ has image satisfying
$|T(V^k_j)| \le C_2 d_{\cW^s}(W^1, W^2)$.  
Choosing $B_i$ as in \eqref{eq:Bi}, we estimate the sum on $V^k_j$ using \eqref{eq:short short} to obtain,
\begin{equation}
\label{eq:cone 3 unmatched}
\begin{split}
\sum_{k=1}^2 \sum_j \int_{V^k_j} h \, \psi_k \circ T \, J_{V^k_j}T
& \le \tri h \tri_- A \delta^{1-q} \big[ 1 + 2C_a \big] \\
&\hskip12pt \times\sum_{k=1}^2 \sum_j 
|\psi_k|_{C^0} |J_{V^k_j}T|_{C^0(V^k_j)} |V^k_j|^q \\
&\hskip-24pt \le 
6 A \delta^{1-q} \tri h \tri_- C_2^q d_{\cW^s}(W^1, W^2)^q \sum_j |J_{V^2_j}T|^{1-q}_{C^0(V^2_j)} \, ,
\end{split}
\end{equation}
where we have used the fact that since $\fint_{W^1} \psi_1 = 1$, we have $|\psi_1|_{C^0} \le e^{a(2\delta)^\alpha} \le 2$, 
and also $|\psi_2 |_{C^0} \le 4$ as in  \cite[eq.~(5.9)]{DL22}.
Since $q<1/2$ and there are at most 2 curves $V^k_j$ corresponding
to each element of $\cG_1^\delta(W^k)$ the final sum converges uniformly for $W^k \in \cW^s$ 
and $T \in \cF(\tau_*, \cK_*, E_*)$
by \cite[Lemma~3.4]{DZ11}; we denote by $C_4$ a uniform bound for it.

Finally, we estimate the contribution to \eqref{eq:focus} from matched pieces $U^k_j$.
For this we will need to change test functions on the matched curves since it may be
that $d_*(\psi_1\circ T \, J_{U^1_j}T, \psi_2 \circ T \, J_{U^2_j}T) \ne 0$.  We define
the following functions on $U^1_j$,
\[
\begin{split}
\tpsi_2 = \psi_2 \circ T \circ G_{U^2_j} \circ G_{U^1_j}^{-1} \, ;
& \quad \tJ_{U^2_j}T = J_{U^2_j}T \circ G_{U^2_j} \circ G_{U^1_j}^{-1} \\
\tT_{U^2_j}(\psi_2) = \tpsi_2 \cdot \tJ_{U^2_j}T \frac{\| G'_{U^2_j} \| \circ G_{U^1_j}^{-1} }
{\| G'_{U^1_j} \| \circ G_{U^1_j}^{-1} } \, ;
& \quad \hT_{U^k_j}(\psi_k) = \psi_k \circ T \, J_{U^k_j}T \, .
\end{split}
\]
Since these functions may not belong to the cone of test functions, we choose a constant
$B_j$ as in \eqref{eq:Bi} so that 
$\hT_{U^2_j}(\psi_2) + B_j \in \cD_{\frac{a}{2},\alpha}(U^2_j)$.\footnote{Increasing
$B_j$ to $2B_j$ is sufficient to scale $a$ to $a/2$.}  Define
\[
\tilde B_j = B_j \cdot \frac{\| G'_{U^2_j} \|} 
{\| G'_{U^1_j} \| } \circ G_{U^1_j}^{-1}  \, .
\]
Then by construction, $d_*(\hT_{U^2_j}(\psi_2) + B_j, \tT_{U^2_j}(\psi_2) + \tilde B_j) = 0$. 
Moreover, the bounds of \eqref{eq:graph control} apply to $G'_{U^k_j}$ since in particular,
$d_{\cW^s}(U^1_j, U^2_j) \le C d_{\cW^s}(W^1, W^2)$ 
for some uniform $C>0$ by \cite[Lemma~3.3]{DZ13}.  This implies, as before, that
$\tilde B_j \le 2 B_j$ and that
$\tT_{U^2_j}(\psi_2) + \tilde B_j \in \cD_{a, \alpha}(U^2_j)$
by \cite[eq.~(5.31)]{DL22}.

For each $j$ we split the difference between matched pieces in \eqref{eq:focus} as follows,
\begin{equation}
\label{eq:again split}
\begin{split}
& \left| \int_{U^1_j} h \hT_{U^1_j}(\psi_1)  - \int_{U^2_j} h \hT_{U^2_j}(\psi_2) \right|
\le \left| \int_{U^1_j} h (\hT_{U^1_j}(\psi_1) - \tT_{U^2_j}(\psi_2)) \right|\\
& \qquad+ \left|  \int_{U^1_j} h \tT_{U^2_j}(\psi_2) - \int_{U^2_j} h \hT_{U^2_j}(\psi_2) \right|  \le \left| \int_{U^1_j} h (\hT_{U^1_j}(\psi_1) - \tT_{U^2_j}(\psi_2)) \right|\\
& \qquad+ \left| \frac{\int_{U^1_j} h (\tT_{U^2_j}(\psi_2) + \tilde B_j) }{\fint_{U^1_j} \tT_{U^2_j}(\psi_2) + \tilde B_j}
- \frac{\int_{U^2_j} h (\hT_{U^2_j}(\psi_2) + B_j)}{\fint_{U^2_j} \hT_{U^2_j}(\psi_2) + B_j} \right| \fint_{U^2_j} (\hT_{U^2_j}(\psi_2) + B_j)
\\
& \qquad + \frac{ \left| \int_{U^1_j} h (\tT_{U^2_j}(\psi_2) + \tilde B_j) \right|}{\fint_{U^1_j} \tT_{U^2_j}(\psi_2)+ \tilde B_j}
\left| \frac{|U^2_j| - |U^1_j|}{|U^1_j|} \right| \fint_{U^2_j} (\hT_{U^2_j}(\psi_2) + B_j) \\
& \qquad+ \left| \int_{U^1_j} h \tilde B_j - \int_{U^2_j} h B_j \right| \, .
\end{split}
\end{equation}
The following sublemma estimates the three most relevant terms of \eqref{eq:again split}.
\begin{sublem}
\label{sub:split again terms}
There exists $C_3>0$, independent of $j$, $W^k$ and $\psi_k$, such that,
\begin{itemize}
  \item[a)] $\left| \int_{U^1_j} h (\hT_{U^1_j}(\psi_1) - \tT_{U^2_j}(\psi_2)) \right|
  \le C_3 d_{\cW^s}(W^1, W^2)^{\alpha-\beta} \delta A \tri h \tri_- |J_{U^1_j}T|_{C^0(U^1_j)}$
  \item[b)] $\displaystyle \left| \frac{\int_{U^1_j} h (\tT_{U^2_j}(\psi_2) + \tilde B_j)}{\fint_{U^1_j} \tT_{U^2_j}(\psi_2) + \tilde B_j}
- \frac{\int_{U^2_j} h (\hT_{U^2_j}(\psi_2) + B_j)}{\fint_{U^2_j} \hT_{U^2_j}(\psi_2) + B_j} \right|
\le C_3 d_{\cW^s}(W^1, W^2)^\gamma \delta^{1-\gamma} cA \tri h \tri_-  $
  \item[c)] $\displaystyle \frac{ \left| \int_{U^1_j} h (\tT_{U^2_j}(\psi_2)+ \tilde B_j) \right|}{\fint_{U^1_j} \tT_{U^2_j}(\psi_2) + \tilde B_j}
\left| \frac{|U^2_j| - |U^1_j|}{|U^1_j|} \right|
\le C_3 d_{\cW^s}(W^1, W^2) \delta A \tri h \tri_-$.
  \end{itemize}
\end{sublem}

We postpone the proof of the sublemma and use it to complete the proof of the bound on the third cone condition.
Using the 3 items of Sublemma~\ref{sub:split again terms} to bound the corresponding terms of
\eqref{eq:again split} yields (recalling again $|\psi_2|_{C^0} \le 4$),
\begin{equation}
\label{eq:intermediate}
\begin{split}
&\left| \int_{U^1_j} h \hT_{U^1_j}(\psi_1)  - \int_{U^2_j} h \hT_{U^2_j}(\psi_2) \right|
\le \left| \int_{U^1_j} h \tilde B_j - \int_{U^2_j} h B_j \right| \\
&+ C_3 d_{\cW^s}(W^1, W^2)^\gamma \delta^{1-\gamma} A \tri h \tri_- \left( \delta^\gamma 
+ 4 (c+\delta) (1 + C_a)  \right)
|J_{U^2_j}T|_{C^0(U^2_j)}  \, .
\end{split}
\end{equation}
Since $d_*(B_j, \tilde B_j) = 0$, it is clear that we can use parts (b) and (c) of Sublemma~\ref{sub:split again terms}
to estimate the last difference of integrals with test functions $B_j$ and $\tilde B_j$, i.e.
\[
\begin{split}
\left| \int_{U^1_j} h \tilde B_j - \int_{U^2_j} h B_j \right| 
& \le \left| \frac{\int_{U^1_j} h  \tilde B_j }{\fint_{U^1_j} \tilde B_j}
- \frac{\int_{U^2_j} h  B_j }{\fint_{U^2_j} B_j} \right| \fint_{U^2_j} B_j
+ \frac{ \left| \int_{U^1_j} h \tilde B_j \right|}{\fint_{U^1_j}  \tilde B_j}
\left| \frac{|U^2_j| - |U^1_j|}{|U^1_j|} \right| \fint_{U^2_j} B_j \\
& \le C_3 d_{\cW^s}(W^1, W^2)^\gamma \delta^{1-\gamma} A \tri h \tri_- \left( 4cC_a + 4 \delta C_a \right) |J_{U^2_j}T|_{C^0(U^2_j)} \, .
\end{split}
\]
Using this estimate together with \eqref{eq:intermediate} completes our estimate on matched pieces,
\[
\sum_j \left| \int_{U^1_j} h \hT_{U^1_j}(\psi_1)  - \int_{U^2_j} h \hT_{U^2_j}(\psi_2) \right|
\le C_0 C_3 d_{\cW^s}(W^1, W^2)^\gamma \delta^{1-\gamma} A \tri h \tri_- \left( \delta^\gamma + 8 c(1+4C_a)  \right)  ,
\]
where we have used $\delta<c$ and the Jacobian to sum over $j$ according to \cite[Lemma~3.3]{DL22}.
Combining this estimate with the estimate on unmatched pieces from \eqref{eq:cone 3 unmatched} in \eqref{eq:focus} yields
finally the upper bound on the left hand side of \eqref{eq:cone 3 prelim}
\begin{equation}
\label{eq:finally}
\begin{split}
& \left|\frac{\int_{W^1} (\cop + \cL h ) \psi_1}{\fint_{W^1}\psi_1} -\frac{\int_{W^2} (\cop+  \cL h) \psi_2}{\fint_{W^2}\psi_2} \right|
 \le  \frac{2d_{\cW^s}(W^1, W^2)^\gamma \left(  C_s^q H_2 \tri h \tri_- + \cop C_s \right) }{\delta^{\gamma-1}}\\
 & \; + C_0 C_3 d_{\cW^s}(W^1, W^2)^\gamma \delta^{1-\gamma} A \tri h \tri_- \left( \delta^\gamma + 8c (1+4C_a)  \right) \\
 & \; + 6 A C_4 \delta^{1-\gamma} d_{\cW^s}(W^1, W^2)^q C_2^q \tri h \tri_- \\
&  =: d_{\cW^s}(W^1, W^2)^\gamma \delta^{1-\gamma} c A \left( H_3 \tri h \tri_- + 2C_s (cA)^{-1}\cop \right)  ,
\end{split}
\end{equation}
where we have used the fact that $\gamma \le q$ and $d_{\cW^s}(W^1, W^2) \le \delta$.

Using \eqref{eq:finally} as the upper bound on the left hand side of \eqref{eq:cone 3 prelim} and \eqref{eq:lower 3}
as the lower bound on the right hand side of \eqref{eq:cone 3 prelim}, we conclude that \eqref{eq:cone 3 prelim}
is satisfied if
\[
H_3 \tri h \tri_- + 2C_s (cA)^{-1} \cop \le \cop - H_1 \tri h \tri_-
\impliedby \cop \ge \frac{H_1 + H_3}{cA - 2C_s} cA \tri h \tri_- \, .
\]
This is a valid choice since $cA > 2C_s$ by \cite[eq.~(5.36)]{DL22}.

Since $\tri h \tri_- \le \| h \|$ by Remark~\ref{rem:equiv tri}, this estimate together with \eqref{eq:cone 1 check}, 
\eqref{eq:cone 2 check}, \eqref{eq:no7.60}  implies that
$\cop + \cL h \in \cC_{\bR}$ whenever $h \in \cC_{\bR}$ and 
\[
\cop \ge \bar C_6 \|h \| := \max \left\{  H_0+2H_1, \frac{H_0 + LH_1}{L-1} , \frac{H_2 + AH_1}{A - 2^{1-q}}, \frac{H_1 + H_3}{cA - 2C_s} cA \right\} \| h \| .
\]
We claim that this suffices to complete the proof of Proposition~\ref{lem:bounded L}.

As in the proof of Corollary~\ref{cor:in space}, for $h \in \cB_{\bR}$, we have 
$h + \| h \|, \| h \| \in \cC_{\bR} \cup \{ 0 \}$.  Thus with $\cop$ chosen as above for the element
$h$, and since $\| h + \| h \| \, \| \le 2 \| h\|$, we may write,
\[
3\cop + \cL h = 2\cop + \cL (h + \| h \|) + (\cop - \| h \|) \, ,
\]
where we have used $\cL 1 = 1$.  Then, since both terms on the right belong to $\cC_{\bR}$,
so does $3 \cop + \cL h$.  Now replacing $h$ with $-h$, and since $\| -h \| = \| h \|$, we have
also that $3\cop - \cL h \in \cC_{\bR}$.  We conclude that $\| \cL h \| \le 3\cop$,
which completes the claim and the proof of 
Proposition~\ref{lem:bounded L} with $C_6= 3 \bar  C_6$.
\begin{proof}[\bf \em Proof of Sublemma~\ref{sub:split again terms}]
(a)  We would like to simply apply the second cone condition \eqref{eq:cone 2} to the integral
in part (a), but $\hT_{U^1_j}(\psi_1) - \tT_{U^2_j}(\psi_2)$ is not necessarily a valid test function.
Following \eqref{eq:Bi}, it suffices to choose 
\[
a^{-1} H^\beta(\hT_{U^1_j}(\psi_1) - \tT_{U^2_j}(\psi_2))
+ |\hT_{U^1_j}(\psi_1) - \tT_{U^2_j}(\psi_2)|_{C^0} \le D_j \le |\hT_{U^1_j}(\psi_1) - \tT_{U^2_j}(\psi_2)|_{C^\beta(U^1_j)}
\]
in order to guarantee that $\hT_{U^1_j}(\psi_1) - \tT_{U^2_j}(\psi_2) + D_j \in \cD_{a, \beta}(U^1_j)$.
We proceed to estimate $|\hT_{U^1_j}(\psi_1) - \tT_{U^2_j}(\psi_2)|_{C^\beta(U^1_j)}$.

First, from \cite[eq.~(5.18)]{DL22} it follows that,
\begin{equation}
\label{eq:c0 close}
|\hT_{U^1_j}(\psi_1) - \tT_{U^2_j}(\psi_2)|_{C^0(U^1_j)} \le C d_{\cW^s}(W^1, W^2)^\alpha 
|\hT_{U^2_j}(\psi_2)|_{C^0(U^2_j)} \, .
\end{equation}
Then by standard estimates (see, for example, the proof of Sublemma~\ref{sub:main diff} or
\cite[Lemma~4.3]{DZ13}), we conclude
\[
|\hT_{U^1_j}(\psi_1) - \tT_{U^2_j}(\psi_2)|_{C^\beta(U^1_j)} 
\le 3 d_{\cW^s}(W^1, W^2)^{\alpha-\beta} \max \{ C |\hT_{U^2_j}(\psi_2)|_{C^0}, H^\beta(\hT_{U^1_j}(\psi_1) + H^\beta(\tT_{U^2_j}(\psi_2)) \} \, . 
\]
However, since $\psi_k \in \cD_{a, \alpha}(W^k)$ and $J_{U^k_j}T$ and 
$\frac{\| G'_{U^1_j} \| }{\| G'_{U^2_j} \| } \circ G_{U^1_j}^{-1}$ enjoy 
similar distortion bounds by \cite[Lemma~5.5]{DL22}, we have
\[
H^\beta(\hT_{U^1_j}(\psi_1) + H^\beta(\tT_{U^2_j}(\psi_2))
\le C' ( |\hT_{U^1_j}(\psi_1)|_{C^0} + |\tT_{U^2_j}(\psi_2)|_{C^0} )
\le C'' |\hT_{U^2_j}(\psi_2)|_{C^0} \, .
\]
Combining this estimate with \eqref{eq:c0 close} yields,
\begin{equation}
\label{eq:Holder ok}
|\hT_{U^1_j}(\psi_1) - \tT_{U^2_j}(\psi_2)|_{C^\beta(U^1_j)}
\le C d_{\cW^s}(W^1, W^2)^{\alpha - \beta} |\hT_{U^2_j}(\psi_2)|_{C^0(U^2_j)} \, ,
\end{equation}
for some uniform $C>0$.  We will use this as our upper bound on $D_j$.

With this bound established, we complete the proof of (a) using the second cone condition
\eqref{eq:cone 2},
\[
\begin{split}
& \left| \int_{U^1_j} h (\hT_{U^1_j}(\psi_1) - \tT_{U^2_j}(\psi_2) + D_j) \right| 
+ \left| \int_{U^1_j} h D_j \right| \\
& \le \left( \fint_{U^1_j} \hT_{U^1_j}(\psi_1) - \tT_{U^2_j}(\psi_2) + 2D_j \right) A |U^1_j|^q \delta^{1-q}
\tri h \tri _- \\
& \le 12 C d_{\cW^s}(W^1, W^2)^{\alpha-\beta} A 2^{1-q} \delta \tri h \tri_- |J_{U^2_j}T|_{C^0(U^2_j)} \, ,
\end{split}
\]
where we have used the previous bound $|\psi_2|_{C^0} \le 4$.

\medskip
\noindent
(b) Since have already verified that $\tT_{U^2_j}(\psi_2) + \tilde B_j \in \cD_{a, \alpha}(U^1_j)$ and $\hT_{U^2_j}(\psi_2)+ B_j \in \cD_{a, \alpha}(U^2_j)$, with $d_*( \tT_{U^2_j}(\psi_2) + \tilde B_j, \hT_{U^2_j}(\psi_2)+ B_j) = 0$, we may apply the third cone condition \eqref{eq:cone 3}
immediately to obtain,
\[
\left| \frac{\int_{U^1_j} h (\tT_{U^2_j}(\psi_2) + \tilde B_j)}{\fint_{U^1_j} \tT_{U^2_j}(\psi_2) + \tilde B_j}
- \frac{\int_{U^2_j} h (\hT_{U^2_j}(\psi_2) + B_j)}{\fint_{U^2_j} \hT_{U^2_j}(\psi_2) + B_j} \right|
\le d_{\cW^s}(U^1_j, U^2_j)^\gamma \delta^{1-\gamma} cA \tri h \tri_- \, .
\]
Statement (b) of the lemma follows using again that $d_{\cW^s}(U^1_j, U^2_j) \le C d_{\cW^s}(W^1, W^2)$ 
for some uniform $C>0$ by \cite[Lemma~3.3]{DZ13}. 

\medskip
\noindent
(c) Since the $U^k_j$ are vertically aligned, we use \eqref{eq:close piece} to bound the difference in curve lengths.
Statement (c) follows immediately using the second cone
condition \eqref{eq:cone 2} to estimate the integral since $|U^1_j| \le 2\delta$.
\end{proof}


\appendix
\section{Complexification}\label{app:complexification}

\begin{proof}[\bf \em Proof of Lemma \ref{lem:complexification}]
We follow \cite[Section 1 (b)]{Sing}. We consider the real vector space $\cB_k^2$ and, for $(x,y)\in\cB_k^2$, define the multiplication by a complex number $a+ib\in\bC$ by
\begin{equation}\label{eq:complex}
(a+ib)(x,y)=(ax-by,ay+bx).
\end{equation}
One can check directly that this defines a complex vector space. Next, we define the norms
\[
\begin{split}
&\|(x,y)\|_{r,k}=\sqrt{\|x\|_k^2+\|y\|_k^2}\\
&\|(x,y)\|_{c,k}=\sup_{\theta\in [0,2\pi]}\|e^{i\theta}(x,y)\|_{r,k}.
\end{split}
\]
Note that, since\footnote{ The second inequality follows  from
\[
\begin{split}
\sqrt{(\|x\|_k^2+\|y\|_k^2)(\|w\|_k^2+\|z\|_k^2)}&\geq \sqrt{\|x\|_k^2\|w\|_k^2+\|y\|_k^2\|z\|_k^2+2\|y\|_k\|w\|_k\|x\|_k\|z\|_k}\\
&=\|x\|_k\|w\|_k+\|y\|_k\|z\|_k.
\end{split}
\]
}
\[
\begin{split}
\|(x,y)+(w,z)\|_{r,k}&\leq \sqrt{\|x\|_k^2+\|w\|_k^2+\|z\|_k^2+\|y\|_k^2+2\|x\|_k\|w\|_k+2\|y\|_k\|z\|_k}\\
&\leq \|(x,y)\|_{r,k}+\|(w,z)\|_{r,k}
\end{split}
\]
$\|\cdot\|_{r,k}$ is a norm for $\cB_k^2$ and with such a norm $\cB_k^2$ is obviously a real Banach space. On the other hand, $\|e^{i\theta}(x,y)\|_{r,k}$ is a continuous function of $\theta$, hence there exists $\theta_*$ such that
\[
\begin{split}
\|(x,y)+(w,z)\|_{c,k}&=\sup_{\theta\in [0,2\pi]}\|e^{i\theta}[(x,y)+(w,z)]\|_{r,k}=\|e^{i\theta_*}[(x,y)+(w,z)]\|_{r,k}\\
&\leq\|e^{i\theta_*}(x,y)\|_{r,k}+\|e^{i\theta_*}(w,z)]\|_{r,k}\leq \|(x,y)\|_{c,k}+\|(w,z)\|_{c,k}.
\end{split}
\]
In addition, for all $\rho\in\bR_+,\vf\in\bR$, we have
\[
\|\rho e^{i\vf}(x,y)\|_{c,k}=\sup_{\theta}\|\rho e^{i\theta}(x,y)\|_{r,k}=\rho\|(x,y)\|_{c,k}=|\rho e^{i\vf}|\|(x,y)\|_{c,k}.
\]
That is $\|\cdot\|_{c,k}$ is a norm for the complex vector space $\bB_k$. Note that, by equation \eqref{eq:complex}, it is natural to write the elements of $\bB_k$ as $x+iy$, $x,y\in \cB_k$.

Since the norm $\|\cdot\|_{c,k}$ dominates the real norm, we have that $\bB_k$ is a complex Banach space. We have then the canonical extension of $\cL_k$ defined as $\cL_k(x+iy)=\cL_k x+i\cL_k y$.
In addition,
\[
\| x + iy \|_{c,k} \le \| x \|_{c, k} + \| y \|_{c,k} = \| x \|_{r, k} + \| y \|_{r, k} \le \sqrt{2} \| x + iy \|_{r,k} \, ,
\]
which implies
\begin{equation}\label{eq:norms_equiv}
\|h\|_{c,k}\leq \sqrt 2\|h\|_{r,k}.
\end{equation}
To conclude the proof note that
\[
\begin{split}
\|\cL_{k+m}\cdots\cL_k(x+iy)\|_{c,k+m+1}&=\|\cL_{k+m}\cdots\cL_k x+i\cL_{k+m}\cdots\cL_k y\|_{c,k+m+1}\\
&\leq \|\cL_{k+m}\cdots\cL_k x\|_{c,k+m+1}+\|i\cL_{k+m}\cdots\cL_k y\|_{c,k+m+1}\\
&=\|\cL_{k+m}\cdots\cL_k x\|_{k+m+1}+ \|\cL_{k+m}\cdots\cL_k y\|_{k+m+1}\\
&\leq C_*(\|x\|_k+\|y\|_k)\leq C_*\sqrt 2\|(x,y)\|_{r,k}\\
&\leq C_*\sqrt 2\|x+iy\|_{c,k}.
\end{split}
\]
\end{proof}


\end{document}